\documentclass[11pt]{article}

\usepackage[T1]{fontenc}
\usepackage[utf8]{inputenc}

\usepackage{graphicx}
\usepackage[a4paper,top=3cm,bottom=3cm,left=2.525cm,right=2.525cm]{geometry}
\usepackage{fancyhdr}

\usepackage[inline]{enumitem}
\usepackage{multicol,multirow,booktabs}

\usepackage{times}

\usepackage[small,labelfont=bf,up,textfont=it,up]{caption}
\usepackage{subcaption}
\usepackage{amsmath,amsfonts,amsthm,amssymb,bm,dsfont}	
\setlist{nolistsep}

\usepackage{numprint}
\npthousandsep{,}\npthousandthpartsep{}\npdecimalsign{.}

\newcommand{\titledoc}{Spectral clustering of network time series via the sample covariance matrix}
\newcommand{\titleshort}{Spectral clustering of network time series via the sample covariance matrix}

\providecommand{\keywords}[1]{{\small{\textbf{\textit{Keywords ---}} #1}}}

\usepackage{fontawesome}
\usepackage{natbib}
\usepackage{setspace}
\usepackage{diagbox}
\usepackage{authblk}

\usepackage{mathtools,bigints}
\mathtoolsset{showonlyrefs}

\usepackage{natbib}

\usepackage{graphicx}
\graphicspath{{Fig/}}

\usepackage{amsmath,amssymb,amsthm,amsfonts}
\allowdisplaybreaks
\usepackage{url} 

\newtheorem{theorem}{Theorem}[section]
\newtheorem{remark}{Remark}%
\newtheorem{definition}[theorem]{Definition}
\newtheorem{lemma}[theorem]{Lemma}
\newtheorem{assumption}[theorem]{Assumption}

\usepackage{dsfont}	
\usepackage{bm}
\usepackage{caption, subcaption}
\usepackage{mathtools}
\usepackage{placeins}
\usepackage{pgfplots}
\pgfplotsset{compat=newest}
\usepackage{standalone}
\usepackage{graphicx}
\usetikzlibrary{pgfplots.dateplot}
\mathtoolsset{showonlyrefs}
\usepackage{booktabs}
\usepackage{appendix}

\DeclareMathOperator{\rank}{rk}
\newcommand{\indfun}{\mathds{1}}

\usepackage[colorlinks=true,allcolors=blue]{hyperref}

\usepackage{algorithm}
\usepackage{algorithmic}
\usepackage{enumitem}

\usepackage{comment}

\DeclareMathOperator*{\argmin}{arg\,min}

\def\acknowledgementsname{Acknowledgments}
\newenvironment{acks}[1][\acknowledgementsname]{\section*{#1}}{\par}
 
\if@balayout
    \renewenvironment{acks}[1][\acknowledgementsname]%
        {%
            \vskip0.5\baselineskip
            \small
            {\noindent\normalfont\sffamily\bfseries #1}\par
            \begingroup\parindent 0pt\parskip 0.5\baselineskip
        }%
        {\endgroup}
\fi

\def\fundingname{Funding}

\def\supplementname{Supplementary Material}
\newenvironment{supplement}[1][\supplementname]{\begin{acks}[\supplementname]}{\end{acks}}

\title{Spectral clustering of network time series via the sample covariance matrix}

\author{{\Large Brendan Martin\textsuperscript{1}, Joshua Agterberg\textsuperscript{2}, Mihai Cucuringu\textsuperscript{1}, \\
\vspace*{-0.5cm}
Alessandra Luati\textsuperscript{3}, and Francesco Sanna Passino\textsuperscript{3}} \\

\textsuperscript{1}Department of Mathematics, University of California, Los Angeles, United States \\
\textsuperscript{2}Department of Statistics, University of Illinois Urbana-Champaign, United States \\
\textsuperscript{3}Department of Mathematics, Imperial College London, United Kingdom
}

\date{}

\title{\Huge\textbf{\titledoc}}

\allowdisplaybreaks

\newcommand\blfootnote[1]{%
  \begingroup
  \renewcommand\thefootnote{}\footnote{\noindent #1}%
  \addtocounter{footnote}{-1}%
  \endgroup
}

\usepackage{multibib}
\newcites{SM}{Supplementary references}

\begin{document}

\maketitle


\begin{abstract}
Spectral clustering for community detection is analysed in multivariate time series models whose dependence structure is determined by an unobserved stochastic {blockmodel}.
We establish that spectral clustering of the sample covariance matrix achieves exact {recovery} of the underlying {communities}. The recovery rates depend explicitly on the network size, sample length, block separation, and degree of data dependence. This {demonstrates} that exact community recovery under a stochastic {blockmodel} is possible even when the adjacency matrix is unobserved. Our theory provides {extensions} of both classical and fine-grained matrix perturbation theory to the setting of dependent data, which may be of independent interest.
\end{abstract}




\blfootnote{
\hspace*{-0.65cm}Corresponding author: Brendan Martin -- \faEnvelopeO\ \texttt{brendan@math.ucla.edu}
}


\section{Introduction}
Statistical inference on random graphs is a well established field of research, with applications in areas such as neuroscience, social science, and econometrics \citep[][]{kolaczyk2014statistical}. Often, the nodes of the graph are divided into communities, and a key task in graph inference is community detection. This problem has been well-studied in the case 
when a realisation of the graph adjacency matrix is observed directly. 
More recently, a line of research has emerged in the multivariate time series literature, 
in which the co-movements among series are driven by a latent underlying graph. Thus, one does not directly observe the network in question. 
A major result of this paper is that we demonstrate that 
imposing a community structure on the underlying graph {is} reflected in the covariance matrix of the resulting time series, 
and
spectral methods applied to the sample covariance matrix of the observed {time} series consistently recover this community structure. 

The stochastic blockmodel \citep[SBM;][]{holland1983stochastic} and its variants such as degree-corrected blockmodels \citep{karrer2011stochastic} and mixed membership blockmodels \citep[][]{airoldi2008mixed} are among the most popular models for graphs with community structure. In a multivariate time series setting, one often observes groups of series exhibiting high intra-group co-movements. Modeling the time series as a vector autoregression (VAR), where the $N \times N$ VAR coefficient matrices are proportional to either the adjacency matrix or normalised Laplacian matrix of an SBM, has been explored in the literature as a way to capture such group-wise co-movements \citep[][]{gudhmundsson2025detecting,martin2024nirvar_arxiv, gudhmundsson2021detecting}. 
The NIRVAR model of \citet[][]{martin2024nirvar_arxiv} and the SB-VAR model of \citet[][]{gudhmundsson2021detecting} study VAR processes in which the VAR coefficient matrices are proportional to the adjacency matrix or normalised Laplacian matrix of an SBM, respectively.  Under the NIRVAR and SB-VAR models, recovery of the underlying communities is key to capturing spillover effects and for downstream tasks such as prediction. 

This paper provides  both strong and weak consistency results for spectral clustering of the sample covariance matrix in the presence of dependence under 
network time series models. 
 ``Strong consistency'' is equivalent to ``exact recovery'' \citep[][]{abbe2018community} and we use the terms interchangeably.  
Exact recovery rates 
are derived under model specifications that belong to the NIRVAR and SB-VAR families.  
To our knowledge, this is the first study to show that spectral clustering achieves exact recovery of the blocks of an SBM when the adjacency matrix is not observed \citep[see][for a survey of cluster analysis of SBMs]{agterberg2023overview}. 
As far as weak recovery is concerned, we provide tight (up to logarithmic terms) weak recovery rates that depend explicitly on the signal-to-noise ratio of the data, and highlight the effect of data dependence on cluster recovery.
 This allows us to compare the final rate to the i.i.d. setting. 


A key theoretical advance of this work is also the extension of both classical and fine-grained matrix perturbation theory to the setting of dependent data. Specifically,  we extend the matrix Bernstein inequality for $\beta$-mixing sequences of \citet[][]{banna2016bernstein} to matrices with unbounded largest eigenvalue and which are not necessarily symmetric. This result is useful beyond the setting of this paper, in applying classical and fine-grained spectral analysis techniques to dependent random matrices. In our context, the result is used to control both the $\ell_{2}$ and $\ell_{2,\infty}$ subspace estimation {errors}. The final bounds depend explicitly on the spectral radius of the VAR coefficient matrix, which controls the dependence between observations of the time series.  

Finally, 
we quantify how low-rank structure in the adjacency matrix of a network 
VAR \citep{zhu2017}  is reflected in the covariance matrix, and hence, sample covariance matrix of the time series. For NIRVAR and SB-VAR, the low-rank structure consists in the blocks of the underlying SBM. In comparison, existing theoretical work either assumes directly that the covariance matrix has a block structure \citep[][]{brownlees2022community} or recovers the block structure not from the sample covariance matrix but instead from an OLS estimate of the VAR coefficient matrix \citep[][]{gudhmundsson2021detecting}.

\subsection{Related literature}\label{sec:lit}
This work is related to several different strands of the literature. First, the literature on group recovery under SBMs \citep[][]{chen2021spectral,abbe2018community,lei2015consistency,rohe2011spectral}, where our contribution is to extend group recovery rates to a setting in which the adjacency matrix is not directly observed.  \citet{brownlees2022community} also prove a weak consistency result for a setting in which the adjacency matrix of a SBM is not directly observed. However, unlike the VAR setting of this paper, they assume directly that the precision matrix is a function of the normalised Laplacian, yielding a spatial lag model.
Second, the literature on mixing sequences \citep[][]{banna2016bernstein,merlevede2011bernstein,merlevede2009bernstein}, where we extend \citet[][Theorem 1]{banna2016bernstein} to unbounded matrices that are not necessarily symmetric. Our proof extends the truncation techniques used in i.i.d. settings \citep[see][Proposition A.3]{hopkins2016fast} to dependent data settings, while we use the standard ``dilation trick'' to deal with matrices that are not necessarily symmetric \citep[see, for example][Section 2.1.17]{tropp2015introduction}. 
Third, the literature on network VAR processes \citep[][]{chen2023community,knight2020generalized,zhu2017,kock2015oracle}, and particularly network VAR processes in which the underlying network is a SBM \citep[][]{gudhmundsson2025detecting,martin2024nirvar_arxiv} where our main contribution is showing that spectral methods applied to the sample covariance matrix instead of an estimate of the VAR coefficient matrix can consistently recover the underlying network structure. 
Lastly, this work is related to spiked covariance models \citep{johnstone2001distribution}. Subspace recovery under spiked covariance models has been well-studied \citep[see, for example][]{bao2022statistical,cai2021subspace,wainwright2019high,wang2017asymptotics} and is closely related to PCA and factor models \citep[see][and references therein]{fan2020statistical}. Existing strong consistency results for spiked covariance models focus on independent data \citep[][]{chen2021spectral}. Our contribution to this literature is a strong consistency result for a spiked covariance model with dependent data.

\subsection{Notation}\label{sec:not} 
For $K\in\mathbb{N}$, let $[K]$ denote the set $\{1,\dots,K\}$. The indicator function is 
$\indfun_{\cdot}\{\cdot\}$, where $\indfun\{\mathcal E\} = 1$ if the event $\mathcal{E}$ occurs, and $0$ otherwise. Let $\mathrm{Card}(S)$ denote the cardinality of a set $S$. We use $\mathcal{S}_{K-1} = \{x \in \mathbb{R}^{K} : \sum_{k=1}^{K}x_{k} = 1, x_{k} > 0, k \in [K]\}$ to denote the interior of the unit simplex in $K$ dimensions. 

For a matrix $M$ and index sets $\mathcal{I}, \mathcal{J} \subseteq [N]$, we let $M_{\mathcal{I},:}$ and $M_{:,\mathcal{J}}$ be the submatrices of $M$ consisting of the corresponding rows and columns, respectively. The transpose of $M$ is $M^{\prime}$, its trace is $\mathrm{tr}(M)$, and its rank is $\rank{(M)}$. Given $\bm{z}=(z_{1},\dots,z_{N})^{\prime} \in \mathbb{R}^{N}$, $\mathrm{diag}(\bm{z})$ is the $N \times N$ diagonal matrix with $\bm{z}$ on the main diagonal. 
We use $\lVert\, \cdot\, \rVert_{2}$ to denote both the Euclidean norm of a vector and the spectral norm of a matrix, depending on the context. The Frobenius norm of a matrix $M$ is written $\lVert M \rVert_{F} = \{\mathrm{tr}(M^{\prime} M)\}^{1/2}$. The $\ell_{2,\infty}$ norm of a matrix $M \in \mathbb{R}^{n \times d}$ is $\lVert M \rVert_{2,\infty} = \max_{1 \leq i \leq n} \lVert M_{i,:} \rVert_{2}$.
The $n \times n$ identity matrix is written $I_{n}$, and we let $\mathbb{M}_{N,K}$ be the collection of all $N \times K$ matrices where each row has exactly one entry equal to 1 and $K -1$ entries equal to 0. 

We say $f(N) = O\{g(N)\}$ or $f(N) \lesssim g(N)$ if there exists a universal constant $C>0$ such that $|f(N)| \leq C |g(N)|$ holds for all sufficiently large $N$. Also, we write $f(N) \asymp g(N)$ if there exist universal constants $C_{1}, C_{2} >0$ such that $C_{1}|g(N)| \leq f(N) \leq C_{2} |g(N)|$ holds for sufficiently large $N$. We write $f(N) \gg g(N)$ if there exists some sufficiently large universal constant $C > 0$ such that $|f(N)| \geq C|g(N)|$.  Similarly, we write $f(N) \ll g(N)$ if there exists a sufficiently small universal constant $C > 0$ such that $|f(N)| \le C |g(N)|$. For ease of presentation, we use $C > 0$ to denote various constants, 
allowed to depend on one another.

\section{Model and estimation}\label{sec:model}
The building block of the network time series models considered in this work is the stochastic blockmodel \citep[SBM;][]{holland1983stochastic}. Definition \ref{def:sbm} below is for an unweighted undirected stochastic blockmodel, which is the setting of the theoretical results in Sections~\ref{sec:adjacency-specification} and~\ref{sec:laplacian-specification}. In Section~\ref{subsec:general-graphs}, we consider the setting of weighted directed graphs.

\begin{definition}[Stochastic blockmodel] \label{def:sbm}
Let $B \in [0,1]^{K \times K}$ be a symmetric connectivity matrix between $K \in {\mathbb{N}}$ communities, let $\alpha_{N} \in (0,1]$ be a sparsity factor, and let $\pi=(\pi_1,\dots,\pi_K) \in \mathcal{S}_{K-1}$ represent group proportions. We say that $A$ is the adjacency matrix for a stochastic blockmodel graph with sparsity factor $\alpha_{N}$, and write $A \sim \mathrm{SBM}(\alpha_{N}B,\pi)$, if the block assignment map $z : [N] \to [K]$ satisfies $\mathbb{P}\{z(i) = k\} = \pi_{k}$ for each $k \in [K]$, and, conditional on $z$, the upper triangular entries $A_{ij}$ are independently generated as 
\begin{align}
A_{ij} &= 
\begin{cases} 
\mathrm{Bernoulli}\bigl(\alpha_{N}\,B_{z(i),z(j)}\bigr) & \text{for} \quad i < j,\\
1 & \text{for} \quad i = j,\\
A_{ji} & \text{for } \quad i > j.
\end{cases}
\end{align}
\end{definition}

    Conditional on the community memberships, 
    the population adjacency matrix $P = \mathbb{E}(A)$ is given by $P = \alpha_{N} ZBZ^{\prime}$, where $Z$ is the community membership matrix, a $N \times K$ binary matrix whose entries satisfy $Z_{ik} = 1$ when vertex $i\in[N]$ belongs to community $k \in [K]$, and zero otherwise. We write the skinny eigendecomposition of $P$ as $P = U_{P} \Lambda_{P} U_{P}^{\prime}$. We also define $\zeta_{k} = \zeta_{k}(Z) = \{1\leq i \leq N : z_{i} = k\}$ and $N_{k} = |\zeta_{k}|$ for all $1 \leq k \leq K$.


Network time series models that incorporate a SBM community structure naturally induce  clustering among the components of the multivariate time series. 
Following the spectral clustering literature \citep{tang2018limit}, we consider two related models. 
  The first, defined in Definition~\ref{def:model-adj}, is based on the adjacency matrix, whereas the second, in Definition~\ref{def:model-laplacian}, uses the normalised Laplacian.

\begin{definition}[Adjacency model] \label{def:model-adj}
    Let $(\bm{X}_{t})_{t \in \mathbb{Z}}$ denote a zero mean, second-order stationary stochastic process where $\bm{X}_{t} = (X_{1,t},\dots,X_{N,t})^{\prime}$ and 
    \begin{align}
    \label{eq:simple-nirvar} 
    \bm{X}_{t} &= \Phi \bm{X}_{t-1} + \bm{\epsilon}_{t} \quad \bm{\epsilon}_{t} \sim \mathcal{N}(\bm{0}, \Sigma),
\end{align}
with $\Phi = \rho A / \lVert A \rVert_{2}$ and $A \sim \mathrm{SBM}(\alpha_{N} B , \pi )$. The $N\times N$ matrix $\Sigma$ denotes the covariance of the innovation term, and $\boldsymbol{0}$ is an $N$-dimensional vector of zeros. 
\end{definition}

\begin{definition}[Laplacian model] \label{def:model-laplacian}
Instate the model in Definition~\ref{def:model-adj} with $\Phi$ replaced by $\Phi_{L}$, with $\Phi_{L} \coloneqq \rho \mathcal{L}(A)$,  where 
$\mathcal{L}(\cdot)$ represents the Laplacian operator 
$\mathcal{L}(A) = [\mathrm{diag}(A\bm{1})]^{-1/2}A[\mathrm{diag}(A\bm{1})]^{-1/2}
$,
and $\bm{1}$ is a vector of $1$s.
\end{definition}


Since $\lVert\Phi\rVert_{2} = \rho$ and $\lVert\Phi_{L}\rVert_{2} = \rho$, both models are stationary if and only if $\rho \in (0,1)$ \citep[][]{lutkepohl2005new}. 
The adjacency model in Definition \ref{def:model-adj} is an instance of the NIRVAR model of \citet[][]{martin2024nirvar_arxiv}. This can be seen by writing the adjacency model as $\bm{X}_{t} = (A \odot \tilde{\Phi})\bm{X}_{t-1} + \bm{\epsilon}_{t}$ with $\tilde{\Phi} = (\rho/\lVert A \rVert_{2})\bm{1}\bm{1}^{\prime}$ being a matrix of fixed weights. While the NIRVAR model allows for multiplex, weighted, and directed graphs, for clarity, we focus on unweighted undirected graphs when presenting theoretical results in Section \ref{sec:adjacency-specification}. 

The Laplacian model in Definition \ref{def:model-laplacian} is an instance of the network vector autoregression model of \citet[][]{zhu2017} with equal network effect and momentum effect ($\beta_{1}=\beta_{2}$ in their notation). Unlike \citet[][]{zhu2017}, we do not assume that the adjacency matrix is observable. \citet[][]{zhu2017} also allow for node-specifice covariates which we do not consider here. The Laplacian model is also an instance of the SB-VAR model of \citet[][]{gudhmundsson2021detecting} with unweighted edges and degree-correction parameters set to $1$. 


For inference under the network time series models in Definition~\ref{def:model-adj} and Definition~\ref{def:model-laplacian}, the main parameter of interest is the community membership matrix $Z$. Under the stochastic blockmodel in Definition~\ref{def:sbm}, the leading $K$ eigenvectors of the population adjacency matrix $P$ reveal the community memberships in the sense that $U_{P} = ZR$ for some matrix $R\in\mathbb{R}^{K \times K}$ \citep[][Lemma 2.1]{lei2015consistency}. 
As such, a substantial portion of this paper is focused on subspace estimation. 
Consider two subspaces $\mathcal{U}_{1}$ and $\mathcal{U}_{2}$, represented by matrices $U_{1} \in \mathbb{R}^{d_{1} \times d_{2}}$ and $U_{2} \in \mathbb{R}^{d_{1} \times d_{2}}$, respectively. 
As in \citet[][]{chen2021spectral}, the optimal rotation distance is defined as 
\begin{align}
    \mathrm{dist}(U_{1},U_{2}) = \min_{R \in \mathcal{O}(d_{2})} \lVert U_{1} - U_{2}R\rVert_{2},
\end{align}
where $\mathcal{O}(d_{2})$ is the 
orthogonal group in $d_2$ dimensions. 

The adjacency matrix $A$ can be thought of as a noisy version of $P$ and it is well-known that $A$ will concentrate around $P$ as long as $N\alpha_{N} \gtrsim \log N$ \citep[][Theorem 5.2]{lei2015consistency}. 
The Davis-Kahan sin$\Theta$ theorem \citep{davis1970rotation} can then be used to control the distance between the leading $K$-dimensional subspaces $\mathcal{U}_{A}$ and $\mathcal{U}_{P}$ of $A$ and $P$, represented by the matrices $U$ and $U_{P}$, respectively, where $U$ is the matrix of the $K$ leading eigenvectors of $A$.
For example, combining Lemma 5.1 and Theorem 5.2 of \citet[][]{lei2015consistency} gives a high probability bound on $\mathrm{dist}(U,U_{P})$. Clustering the rows of $U$ (using, say, $k$-means) will then lead to consistent recovery of $Z$. 
Analogous results to adjacency spectral embedding under an SBM exist for Laplacian spectral embedding: \citet[][Lemma 3.1]{rohe2011spectral} show that the leading $K$ eigenvectors of the population normalised Laplacian $\mathcal{L}(P)$ reveal the community memberships; \citet[][Theorem 4.1]{ke2025optimal} show that $\mathcal{L}(A)$ concentrates around $\mathcal{L}(P)$. 

For the network time series models in Definition \ref{def:model-adj} and Definition \ref{def:model-laplacian}, spectral embedding of the matrices $A$ or $\mathcal{L}(A)$ is not viable, since the network adjacency matrix $A$ is unobserved. Therefore, for estimation of the community membership matrix $Z$, we instead consider spectral embedding of the sample covariance matrix $S_{T} = \sum_{t=1}^{T} \bm{X}_{t}\bm{X}_{t}^{\prime}/T$. Algorithm~\ref{alg:spectral} specifies one possible spectral clustering procedure under this framework. 

\begin{algorithm}[ht]
\caption{A representative spectral clustering algorithm.}
\label{alg:spectral} 
\begin{algorithmic}[1]
\normalsize
\REQUIRE Sample covariance matrix~$S_{T}$, embedding dimension~$d$, 
number of 
groups~$K$. 
\\
\STATE 
Compute the spectral embedding of $S_T$ in $\mathbb{R}^d$ by forming the $N \times d$ matrix $U_S$ whose orthonormal columns are the eigenvectors of $S_{T}$ associated with the $d$ eigenvalues of largest magnitude, collected in the diagonal matrix $\Lambda_S$. \\
\STATE 
Cluster the rows of $U_{S}$ in $K$ groups using a clustering algorithm (for example, $k$-means).\\ 
\STATE Return the vector of 
estimated community memberships $\hat{\bm{z}}\in \{1,\dots,K\}^{N}$, or 
the estimated membership matrix $\hat{Z}\in \{0,1\}^{N \times K}$.
\end{algorithmic}
\end{algorithm}


\section{Theoretical results} \label{sec:theory} 


The rates derived in this section rely on a novel lemma (Lemma~\ref{lemma:truncated-Bernstein}) 
which extends the Bernstein-type inequality for dependent random matrices of \citet{banna2016bernstein} to matrices whose largest eigenvalue is not bounded almost surely and which are not necessarily symmetric. 
This result may be of independent interest to those working with sequences of $\beta$-mixing random matrices. The proof of Lemma~\ref{lemma:truncated-Bernstein} is in Appendix~\ref{sec:appendix}. The proofs of Lemma~\ref{lemma:eigs}, which relates the eigenvalues of the adjacency matrix $A$ to the eigenvalues of the block probability matrix $B$, and Lemma~\ref{lemma:S_T-Gamma}, which bounds the spectral norm covariance estimation error, are also given in Appendix~\ref{sec:appendix}. These two novel lemmas are crucial for subsequent results. The proofs of all other theoretical results are given in the online supplementary material. 

\subsection{A truncated matrix Bernstein inequality for dependent data}\label{sec:trunc-bern}

The sequences we consider are geometrically $\beta$-mixing, also referred to as geometrically absolutely regular \citep[see][for a survey of mixing conditions]{bradley2005basic}. For the remainder of the paper, let $c_{\beta} > 0$ control the geometric decay of the $\beta$-mixing coefficients.

\begin{lemma}(Truncated matrix Bernstein for $\beta$-mixing sequences)\label{lemma:truncated-Bernstein}
    Let $(M_{t})_{t \in [T]}$ be a family of $d_{1} \times d_{2}$ matrices whose entries belong to $\mathbb{R}$. Assume that $(M_{t})_{t \in \mathbb{Z}}$ is a geometrically $\beta$-mixing sequence with  mixing coefficient $c_{\beta}$. Let 
    \begin{align}\label{def:mv}
        v^{2} = \sup_{\tau \subseteq [T]} \frac{1}{\mathrm{Card}(\tau)} \max\Bigg{\{} &\left\lVert \sum_{t,s \in \tau} \mathbb{E}[(M_{t} -\mathbb{E}[M_{t}])(M_{s} -\mathbb{E}[M_{s}])^{\prime}]\right\rVert_{2}, \\ &\left\lVert \sum_{t,s \in \tau} \mathbb{E}[(M_{t} -\mathbb{E}[M_{t}])^{\prime}(M_{s} -\mathbb{E}[M_{s}])]\right\rVert_{2}  \Bigg{\}}.
    \end{align}
    Suppose that for all $0 \leq t \leq T$, 
    \begin{gather}
        \mathbb{P}\{\lVert M_{t} - \mathbb{E}(M_{t}) \rVert_{2} \geq L \} \leq q_{0}, \\
        \lVert \mathbb{E}(M_{t}) - \mathbb{E}(M_{t}\indfun\{\lVert M_{t} - \mathbb{E}(M_{t})\rVert_{2} \geq L\}) \lVert_{2} \leq q_{1}
    \end{gather}
    hold for some $q_{0} \lesssim 1/\{\mathbb{E}(\lVert M_{t} - \mathbb{E}[M_{t}] \rVert_{2}^{4}) (\mathrm{Card}(\tau))^{4}\}$ and $q_{1} \geq 0$. Then, for all $x \geq Tq_{1}$ with $T \geq 2$, there exists $C > 0$ such that
    \begin{equation}
        \mathbb{P}\left( \left\lVert \sum_{t =1}^{T} M_{t} - \mathbb{E}(M_{t}) \right\rVert_{2} \geq x \right) \leq  Tq_{0} + (d_{1} + d_{2})\exp\left(- \frac{C(x-Tq_{1})^{2}}{v^{2}T + c_{\beta}^{-1}L^{2} + (x-Tq_{1})L\gamma(c_{\beta},T)} \right),
    \end{equation}
    where 
    \begin{equation}\label{eq:gamma}
                \gamma(c_{\beta},T) \coloneqq \frac{\log T}{\log 2} \max \left(2, \frac{32 \log T}{c_{\beta} \log 2}\right).
    \end{equation}
\end{lemma}
The proof is in Appendix~\ref{sec:trunc-bern-proof}. Operationally,
if $\mathbb{E}(\lVert M_{t} - \mathbb{E}[M_{t}] \rVert_{2}^{4}) \lesssim d^{4}$, where $d = \max\{d_{1},d_{2}\}$, then  $q_{0} = (d+T)^{-\nu}$ for  $\nu \geq 8$ will satisfy the requirement in Lemma~\ref{lemma:truncated-Bernstein}. 

\subsection{Clustering results for the adjacency model}
\label{sec:adjacency-specification}
This section details theoretical results for the adjacency model in Definition~\ref{def:model-adj}. 

\begin{assumption}[Network regularity conditions] \label{ass:regularity}
    Under the SBM model in Definition~\ref{def:sbm}, we assume that $\mathrm{rank}(B)=K$ 
    where $\lambda_{1}(B),\dots,\lambda_{K}(B)$ are the ordered eigenvalues of $B$. Furthermore, the communities are balanced: $|\zeta_{r}(Z)| \asymp |\zeta_{s}(Z)|$ for $r,s \in [K]$. The sparsity parameter satisfies $\alpha_{N} \geq c_{0}\log N /N$ for some $c_{0} > 0$. 
\end{assumption}

The assumption of balanced communities is standard in the literature \citep{su2019strong} but can be relaxed by keeping track of constants. It is well-known that exact recovery under an SBM is possible only in the sparsity regime $\alpha_{N} \geq c_{0}\log N /N$ \citep[see][for a review of the fundamental limits for
community detection in SBMs]{abbe2018community}. 

For subspace estimation under spiked covariance models, if the noise matrix $\Sigma$ is heteroskedastic, spectral clustering via Algorithm~\ref{alg:spectral} leads to biased estimates. Indeed, unless additional structure, such as homoskedasticity, is imposed on $\Sigma$, the leading eigenvectors of the covariance matrix will be different from the eigenvectors of the low-rank spiked matrix. In the literature, ideas such as diagonal deletion \citep[see, for example][]{yan2024inference,abbe2022lp,cai2021subspace} and rescaling diagonal entries \citep[see, for example][]{zhang2022heteroskedastic} have been recommended to mitigate the bias introduced when there is heteroskedastic noise. For Sections~\ref{sec:adjacency-specification} and~\ref{sec:laplacian-specification}, we assume homoskedastic noise. This assumption will be relaxed in Section \ref{subsec:general-graphs}. 

\begin{assumption}[Homoscedasticity] \label{ass:homoscedasticity}
    The covariance matrix of the innovation term in~\eqref{eq:simple-nirvar} is proportional to the identity: $\Sigma 
    \propto I_{N}$. 
\end{assumption}

Our goal is to recover the SBM community structure, encoded in $U_{P}$, from the leading $K$ eigenvectors of $S_{T}$, denoted $U_{S}$. \citet[][Proposition 4.1]{martin2024nirvar_arxiv} show that, assuming homoskedastic errors, the eigenvectors of $\Gamma$ and $\Phi$ are the same. This allows us to control the distance between $\mathcal{U}_{S}$ and $\mathcal{U}_{P}$ by splitting the problem using the triangle inequality: $\mathrm{dist}(U_{S},U_{P}) \leq \mathrm{dist}(U_{S},U) + \mathrm{dist}(U,U_{P})$, where $U$ contains the simultaneous $K$ leading eigenvectors of $\Gamma$ and $A$. 
    We do not use subscripts for the eigenvector matrix $U$ since its columns are simultaneously the $K$ leading eigenvectors of $\Gamma$, $\Phi$, and $A$.
The term $\mathrm{dist}(U_{S},U)$ arises from the statistical estimation error incurred from using the sample covariance matrix $S_{T}$ as an estimate of the covariance matrix $\Gamma$. In contrast, $\mathrm{dist}(U,U_{P})$ is a structural error arising from the difference between the adjacency matrix $A$ and its expected value $P$. As such, this term does not depend on the number of observed time points $T$, but rather vanishes as the number of nodes in the network $N$ increases. 

We control $\mathrm{dist}(U_{S},U)$ via the Davis-Kahan sin$\Theta$ theorem. First, we show that $S_{T}$ concentrates around $\Gamma$. Lemma~\ref{lemma:S_T-Gamma}, whose proof is deferred to Appendix~\ref{sec:proof3_7}, bounds the size of the perturbation $\lVert S_{T} - \Gamma \rVert_{2}$ and is a main result of this paper. To express the bound in Lemma~\ref{lemma:S_T-Gamma} in terms of the eigenvalues of $B$, we relate the eigenvalues of $\Phi$, denoted $\lambda_{r}(\Phi)$, to the eigenvalues of $B$ through Lemma \ref{lemma:eigs}, whose proof is given in Appendix~\ref{sec:proof3_4}.


\begin{lemma}\label{lemma:eigs}
    Let $\lambda_{1}(\Phi),\dots,\lambda_{N}(\Phi)$ and $\lambda_{1}(B),\dots,\lambda_{K}(B)$ be the ordered eigenvalues of the matrices $\Phi=\rho A/\lVert A \rVert_{2}$ and $B$, respectively. Let Assumption~\ref{ass:regularity} hold. Then for any $s>0$, with probability at least $1 - N^{-s}$, 
    \begin{align}
        \lambda_{r}(\Phi) \asymp \rho\lambda_{r}(B),\ r = 1,\dots,K, & & 
        \lambda_{r}(\Phi) = o(1),\ r = K+1,\dots,N.
    \end{align}
\end{lemma}

\begin{lemma}\label{lemma:S_T-Gamma}
Suppose that 
    \begin{equation}\label{eq:lemma_T_condition}
        T > C \frac{\gamma(c_{\beta},T)^{2}}{1 - \rho^{2}}\log^{3} (T+N)
    \end{equation}
for some sufficiently large constant $C >0$ where $\gamma(c_{\beta},T)$ is defined in \eqref{eq:gamma}.
Then for any $\nu \geq 8$, with probability exceeding $1 - (T+N)^{-\nu}$, 
    \begin{align}\label{eq:S-Gamma-rate}
        \left\lVert S_{T} - \Gamma \right\rVert_{2} &\lesssim \Bigg{[} \frac{1}{(1 - \rho^{2})^{3/2}} \sqrt{\frac{K}{T}} +  \left\{1 + \frac{1}{\sqrt{ 1 - \rho^{2}}}\right\}\sqrt{\frac{N}{T}} \Bigg{]} \log^{1/2}(T+N).
    \end{align}
\end{lemma}

\begin{remark}\label{remark:cov-est-rho}
    The covariance estimation error in \eqref{eq:S-Gamma-rate} increases with $\rho$, since a larger $\rho$ corresponds to a higher degree of data dependence, thus reducing the information contained in the data for fixed $T$.
\end{remark}
A high probability bound on $\mathrm{dist}(U_{S},U)$ is obtained by combining Lemma~\ref{lemma:S_T-Gamma} with a variant of the Davis-Kahan sin$\Theta$ theorem 
in \citet[][Corollary 2.8]{chen2021spectral}: if $\lVert S_{T} - \Gamma\rVert_{2} \leq (1-1/\sqrt{2})\Delta_{\Gamma}$, then 
    $\mathrm{dist}(U_{S},U) \leq 2\lVert S_{T} - \Gamma\rVert_{2} / \Delta_{\Gamma}$,
where $\Delta_{\Gamma} \coloneqq \lambda_{K}(\Gamma) - \lambda_{K+1}(\Gamma)$ is the eigengap. 
We now state a high probability bound on $\mathrm{dist}(U_{S},U_{P})$, given by Theorem~\ref{thm:dist_U_S_U_P}. 

\begin{theorem}\label{thm:dist_U_S_U_P}
Let Assumption~\ref{ass:regularity} hold and let $\gamma(c_{\beta},T)$ be defined as in \eqref{eq:gamma}. 
Suppose  
    \begin{align}
        T &\geq C \left[\left\{ \frac{1}{\rho^{2} \lambda_{K}^{2}(B) }\right\}^{2}\left\{ 
        \frac{K}{(1 - \rho^{2})^{3}} 
        +  \left(1 + \frac{1}{\sqrt{ 1 - \rho^{2}}}\right)^{2}N\right\} + \frac{\gamma(c_{\beta},T)^{2}}{1 - \rho^{2}} \log^{2} T\right]\log^{3} T
        \label{eq:T_cond_combined}
    \end{align}
    for some sufficiently large constant $C>0$.
    Then, for any $\nu \geq 8$, with probability at least $1 - (T+N)^{-\nu}$,
    \begin{multline}\label{eq:dist_U_P}
        \mathrm{dist}(U_{S},U_{P}) \lesssim \\
        \frac{K\sqrt{K}}{\sqrt{N \alpha_{N}}\lambda_{K}(B)} 
         + \left\{ \frac{1}{\rho^{2} \lambda_{K}^{2}(B) }\right\}\left\{ \frac{1}{(1 - \rho^{2})^{3/2}} \sqrt{\frac{K}{T}} +  \left(1 + \frac{1}{ \sqrt{1 - \rho^{2}}}\right)\sqrt{\frac{N}{T}}\right\} \log^{1/2}(T+N).
    \end{multline}
\end{theorem}

    The rate in~\eqref{eq:dist_U_P} consists of two terms: the first arises from the error in estimating $U_{P}$ via $U$, while the second arises from the error in estimating $U$ via $U_{S}$. The first term is the standard error incurred in estimating the blocks of an SBM from a single realised adjacency matrix \citep[][]{lei2015consistency}, while the second term goes to zero as $T \to \infty$ (for fixed $N$).

\begin{remark}\label{remark:SNR}
    From~\eqref{eq:dist_U_P}, $\mathrm{dist}(U_{S},U_{P}) \xrightarrow[]{} \infty$ as $\rho \to 0$. This is expected since if $\rho = 0$, then $\bm{X}_{t} \sim \mathcal{N}(\bm{0},\Sigma)$, and $\Gamma$ does not inherit any low-rank 
    structure from $\Phi$, i.e. there is no signal. 
\end{remark}

\begin{remark}
    The second term of the bound in~\eqref{eq:dist_U_P} can be expressed in terms of the noise $\sigma_{\Gamma}^{2} \coloneqq \lambda_{K+1}(\Gamma) $, the condition number of $\Gamma$, labelled $\kappa_{\Gamma} \coloneqq \lambda_{1}(\Gamma)/\lambda_{K}(\Gamma)$, and the spectral radius $\rho$ which determines the strength of the data dependence:
    \begin{equation}\label{eq:noise-form-distUS}
        \mathrm{dist}(U_{S},U) \lesssim \, \frac{1}{\sqrt{1-\rho^{2}}}\left\{ \kappa_{\Gamma} \sqrt{\frac{K}{T}} + \frac{\sigma_{\Gamma}}{\sqrt{\lambda_{K}(\Gamma)}}\sqrt{\frac{\kappa_{\Gamma}N}{T}} + \frac{\sigma_{\Gamma}^{2}}{\lambda_{K}(\Gamma)}\sqrt{\frac{N}{T}}\right\} \log^{1/2}(T+N).
    \end{equation}
    Writing the bound in this form highlights the linear and quadratic dependence on the noise level $\sigma_{\Gamma}$ that is typical of spectral methods applied to the sample covariance matrix. In particular, we compare to the subspace estimation rate in \citet[][Theorem 3.6]{chen2021spectral} who consider a spiked covariance model with i.i.d. data. The rate in \eqref{eq:noise-form-distUS} is the same as \citet[][Theorem 3.6]{chen2021spectral} except for the $1/\sqrt{1-\rho^{2}}$ factor which is due to data dependence. 
\end{remark}

Two notions of community recovery are considered: \emph{weak recovery}, namely recovering all but a vanishing fraction of vertex memberships with probability tending to one, and \emph{exact recovery}, namely recovering all vertex memberships exactly with probability tending to one. Community recovery under the spectral clustering procedure in Algorithm~\ref{alg:spectral} will be possible if the leading $K$-dimensional eigenspace of $S_{T}$ is close to that of $P$. For weak recovery,  closeness is measured in terms of $\mathrm{dist}(U_{S},U_{P})$, which we have bounded in Theorem \ref{thm:dist_U_S_U_P}, whereas for exact recovery the appropriate metric is $\min_{R \in \mathcal{O}(K)} \lVert U_{S} - U_{P}R\rVert_{2, \infty}$ \citep[][]{chen2021spectral,fan2020statistical}.  

\subsubsection{Weak recovery}
Let us consider estimation of the group membership matrix via $k$-means clustering, defined as 
\begin{align}\label{eq:k-means}
    (\hat{Z}, \hat{\Theta}) = \argmin_{Z \in \mathbb{M}_{N,K}, \Theta \in \mathbb{R}^{K \times K}} \lVert Z \Theta - U_{S}  \rVert_{F}.
\end{align}
Finding a global minimiser of \eqref{eq:k-means} is known to be non-deterministic polynomial-time hard \citep{aloise2009np}. We find an approximate solution whose value is within a constant fraction of the optimal value. In particular, there exist polynomial-time algorithms \citep[see][]{kumar2010clustering} that find $(\hat{Z}, \hat{\Theta}) \in \mathbb{M}_{N,K} \times \mathbb{R}^{K \times K}$ such that 
\begin{align}
    \lVert \hat{Z} \hat{\Theta} - U_{S}  \rVert_{F} \leq (1+\varepsilon)\min_{Z \in \mathbb{M}_{N,K}, \Theta \in \mathbb{R}^{K \times K}} \lVert Z \Theta - U_{S}  \rVert_{F}.
\end{align}
Under this approximate $k$-means procedure, Theorem~\ref{thm:kmeans} states a recovery rate for weak clustering.

\begin{theorem}\label{thm:kmeans}
    Let $\hat{Z}$ be the output of spectral clustering using $(1+\varepsilon)$-approximate $k$-means. Then, under the assumptions of Theorem~\ref{thm:dist_U_S_U_P}, with probability at least $1 - (T+N)^{-\nu}$ for $\nu \geq 8$, there exist subsets $S_{k} \subset \zeta_{k}$ for $k = 1, \dots, K$ and a $K \times K$ permutation matrix $J$ such that $\hat{Z}_{\zeta,:}J = Z_{\zeta,:}$ where $\zeta = \cup_{k=1}^{K} (\zeta_{k} \setminus S_{k})$ and 
    \begin{multline}
        \sum_{k=1}^{K} \frac{|S_{k}|}{N_{k}} \lesssim \\
        \frac{K^{4} }{N \alpha_{N}\lambda_{K}^{2}(B)} + 
        \left\{\frac{1}{\rho^{2} \lambda_{K}^{2}(B) }\right\}^{2} \left\{ \frac{1}{(1 - \rho^{2})^{3}} \frac{K^{2}}{T} +  \left(1 + \frac{1}{ \sqrt{1 - \rho^{2}}}\right)^{2}\frac{KN}{T}\right\} \log(T+N).
        \label{eq:weak-cluster-rate}
    \end{multline}
\end{theorem}

Theorem~\ref{thm:kmeans} establishes that the fraction of SBM vertices which are misclustered by the spectral method in Algorithm \ref{alg:spectral} is negligible with high probability as long as $N/T$ is small and the SBM is sufficiently dense. The second term in \eqref{eq:weak-cluster-rate} decreases with $\rho$ so that higher data dependence leads to a smaller fraction of misclustered vertices. This might appear in conflict with Remark~\ref{remark:cov-est-rho}, which notes that the covariance estimation error in~\eqref{eq:S-Gamma-rate} increases with $\rho$. However, since the SNR (eigengap) increases at a faster rate in $\rho$, the overall misclustering rate decreases with $\rho$. 

\subsubsection{Exact recovery}
In order to establish exact recovery results for the spectral clustering algorithm, controlling the entry-wise behaviour of $U_{S}$ is crucial. Theorem~\ref{thm:strong-cov-estimation} bounds the $\ell_{2,\infty}$-error between the eigenvectors of the sample covariance matrix and their population counterparts.

\begin{theorem}[Covariance estimation] \label{thm:strong-cov-estimation}
    Instate the assumptions of Theorem~\ref{thm:dist_U_S_U_P}.  
    There exists an orthogonal matrix $W_{U}$ such that for any $\nu \geq 8$ , with probability $1 - (T+N)^{-\nu}$,
    \begin{align}\label{eq:strong-cov-bound}
    \|U_{S} - U W_{U}\|_{2 ,\infty} \lesssim &\  \frac{1}{\rho^{2} \lambda_{K}^{2}(B) } \Bigg{[} \left\{ \frac{1}{(1 - \rho^{2})^{3/2}} \sqrt{\frac{K^{3}}{NT}} +  \frac{1}{\sqrt{1 - \rho^{2}}}  \sqrt{\frac{K^{2}}{T}} \right\} \log^{1/2}T \\ 
    &+  \left\{ \frac{1}{(1 - \rho^{2})^{3/2}} \frac{\sqrt{NK}}{T} +   \frac{1}{ \sqrt{1 - \rho^{2}}}\frac{N}{T}\right\} \log T \log^{1/2}(T+N) \\
    &+  \left\{ \frac{1}{(1 - \rho^{2})^{3/2}} \sqrt{\frac{K}{T}} +   \frac{1}{ \sqrt{1 - \rho^{2}}}\sqrt{\frac{N}{T}}\right\} \log^{1/2}(T+N) \\
    &+ \left\{ \frac{1}{(1 - \rho^{2})^{3}} \frac{K}{T} + \frac{1}{(1 - \rho^{2})^{2}} \frac{\sqrt{NK}}{T} +   \frac{1}{ 1 - \rho^{2}}\frac{N}{T}\right\} \sqrt{\frac{K}{N}}\log(T+N) \Bigg{]}.
\end{align}
\end{theorem}

As with the spectral norm bound, by the triangle inequality  
\begin{align}
    \min_{R \in \mathcal{O}(K)}\lVert U_{S} - U_{P}R \rVert_{2,\infty} \leq \min_{R \in \mathcal{O}(K)}\lVert U_{S} - UR\rVert_{2,\infty} + \min_{R \in \mathcal{O}(K)}\lVert U - U_{P}R\rVert_{2,\infty}.
\end{align}
Under Assumption \ref{ass:regularity}, by \citet[][Theorem 1]{agterberg2023overview}, with high probability,
\begin{align}
    \min_{R \in \mathcal{O}(K)}\lVert U - U_{P}R\rVert_{2,\infty} = O\left( \frac{\sqrt{K \log N}}{N\alpha_{N}^{1/2}} \right).
\end{align}
Combining the above equality with Theorem~\ref{thm:strong-cov-estimation} yields a bound on $\min_{R \in \mathcal{O}(K)}\lVert U_{S} - U_{P}R \rVert_{2,\infty}$ and allows us to prove that the $k$-means algorithm gives exact recovery in the regime $T \gtrsim N^{2}\log^{2}N$, stated in Theorem~\ref{thm:strong-clustering}. If the network were observable, Assumption~\ref{ass:regularity} would be sufficient for exact recovery \citep[][Theorem 6]{lyzinski2014perfect}; the additional conditions in Theorem~\ref{thm:strong-clustering} are due to the network being latent. 

\begin{theorem}[Strong clustering] \label{thm:strong-clustering}
Instate the assumptions of Theorem~\ref{thm:strong-cov-estimation}. If $T > CN^{2}\log^{2}N$ for some $C>0$, then for any $\nu \geq 8$, with probability at least $1 - (N+T)^{-\nu}$,
    $\sup_{1 \leq i \leq N} \indfun\{ \hat{z}_{i} \neq z_{i} \} = 0.$
\end{theorem}
Given Theorem~\ref{thm:strong-cov-estimation}, the proof of Theorem~\ref{thm:strong-clustering}  
follows a similar strategy to that of \citet[][Theorem 6]{lyzinski2014perfect}. 

    We conclude this section by observing that a substantial portion of the literature over the past twenty years has been dedicated to $\ell_{2}$ and $\ell_{2,\infty}$ perturbation theory in the setting of i.i.d. data \citep[see, for example][]{cai2021subspace,abbe2020entrywise,cape2019two,chen2019spectral,el2013robust}. The core technical tool used to obtain such fine-grained error control is ``leave-one-out'' analysis \citep[see][Chapter~4 for an overview of these proof techniques]{chen2021spectral}. The assumption of i.i.d. data is of crucial importance in leave-one-out analysis; as noted by \citet[][Section~5]{chen2021spectral}, handling dependency structure ``might require ideas beyond the current leave-one-out framework''. In our setting of both temporal and cross-sectional dependence, standard leave-one-out analysis does not apply. Instead, to control the entry-wise covariance estimation error in Theorem \ref{thm:strong-cov-estimation}, we use \citet[][Theorem~3.7]{cape2019two}. The bound in Theorem~\ref{thm:strong-cov-estimation} implies that strong clustering will hold if $T \gtrsim N^{2}\log^{2}N$, which is assumed in Theorem~\ref{thm:strong-clustering}.

\subsection{Clustering results for the Laplacian model}\label{sec:laplacian-specification}
A row-wise large deviation bound on the eigenvectors of the normalized graph Laplacian has been 
established by \citet[][Theorem~4.1]{ke2025optimal}. Assumption~\ref{ass:regularity-laplacian} imposes mild regularity conditions on the network, allowing the use of \citet[][Theorem~4.1]{ke2025optimal}. 

\begin{assumption}[Regularity conditions for Laplacian spectral embedding] \label{ass:regularity-laplacian}
    There exist constants $c_{2}>0$ and $c_{3}>0$ such that 
    \begin{enumerate}[label=(\alph*)]
    \item $\max_{k \neq 1}\lambda_{K}(B) \leq (1-c_{2})\lambda_{1}(B)$ and $\min_{1 \leq k \leq K}\{\sum_{1 \leq l \leq K} B_{kl}\} \geq c_{2}K$,
    \item $\min_{1 \leq k \leq K} \{(U_{B})_{k1}\} \geq c_{3} \max_{1 \leq k \leq K}\{(U_{B})_{k1}\}$, where $(U_{B})_{:1}$ is the leading eigenvector of $B$, with eigenvalue $\lambda_{1}(B)$.
    \end{enumerate}
\end{assumption}

Assumption~\ref{ass:regularity-laplacian} is equivalent to \citet[][Condition 2.1(b)-(c)]{ke2025optimal}. The first part of Assumption~\ref{ass:regularity-laplacian}(a) is an eigengap condition, 
strengthening 
Perron's theorem  \citep{horn2012matrix}. The second part of (a) requires that 
each row sum
of $B$ is $O(K)$. By Perron's theorem  \citep{horn2012matrix}, $(U_{B})_{:1}$ is a positive vector when $B$ is irreducible; hence, Assumption~\ref{ass:regularity-laplacian}(b) is only slightly stronger than requiring 
irreducibility of $B$. 


As in Section~\ref{sec:adjacency-specification}, we relate $\lambda_{r}(\Gamma)$ to $\lambda_{r}(\Phi_{L})$ via \citet[][Proposition 4.1]{martin2024nirvar_arxiv} and also assume homoskedastic errors (Assumption \ref{ass:homoscedasticity}).
The eigenvalues $\lambda_{r}(\Phi_{L})$ are further related to $\lambda_{r}(B)$ through the following lemma. 
\begin{lemma}\label{lemma:eigs-laplacian}
    Let Assumptions~\ref{ass:regularity} and \ref{ass:regularity-laplacian} hold. Then the results of Lemma~\ref{lemma:eigs} hold with $\Phi$ replaced by $\Phi_{L}$.
\end{lemma}

We now bound the distance between the leading $K$-dimensional subspaces $\mathcal{U}_{S}$ and $\tilde{\mathcal{U}}_{P}$ of $S_{T}$ and $\mathcal{L}(P)$, represented by the matrices $U_{S}$ and $\tilde{U}_{P}$, respectively.

\begin{theorem}\label{thm:dist_U_S_U_P-lap}
Let Assumptions~\ref{ass:regularity} and~\ref{ass:regularity-laplacian} hold, and let $\gamma(c_{\beta},T)$ be defined as in Lemma~\ref{lemma:S_T-Gamma}. 
Suppose 
    \begin{align}
        T &\geq C \left[\left\{ \frac{1}{\rho^{2} \lambda_{K}^{2}(B) }\right\}^{2}\left\{ \frac{K}{(1 - \rho^{2})^{3}} 
        +  \left(1 + \frac{1}{\sqrt{ 1 - \rho^{2}}}\right)^{2}N\right\} + \frac{\gamma(c_{\beta},T)^{2}}{1 - \rho^{2}} \log^{2} T\right]\log^{3} T 
    \end{align}
    for some sufficiently large constant $C>0$.
    Then, for any $\nu \geq 8$, with probability at least $1 - (T+N)^{-\nu}$,
    \begin{multline} 
        \mathrm{dist}(\tilde{U}_{S},\tilde{U}_{P}) \lesssim \\
        \frac{\log^{1/2}N}{\sqrt{N} \alpha_{N}^{3/4}\lambda_{K}(B)} + \left\{ \frac{1}{\rho^{2} \lambda_{K}^{2}(B) }\right\}\left\{ \frac{1}{(1 - \rho^{2})^{3/2}} \sqrt{\frac{K}{T}} +  \left(1 + \frac{1}{ \sqrt{1 - \rho^{2}}}\right)\sqrt{\frac{N}{T}}\right\} \log^{1/2}(T+N).
    \end{multline}
\end{theorem}

\subsubsection{Weak recovery}
Theorem~\ref{thm:dist_U_S_U_P-lap} gives the same rate as in Theorem~\ref{thm:dist_U_S_U_P}, therefore allowing us to prove weak recovery of the vertex memberships under the Laplacian model.

\begin{theorem}\label{thm:kmeans-lap}
    Let $\hat{Z}$ be the output of spectral clustering using $(1+\varepsilon)$-approximate $k$-means. Then under the assumptions of Theorem~\ref{thm:dist_U_S_U_P-lap}, with probability at least $1 - (T+N)^{-\nu}$ for $\nu \geq 8$, there exist subsets $S_{k} \subset \zeta_{k}$ for $k = 1, \dots, K$ and a $K \times K$ permutation matrix $J$ such that $\hat{Z}_{\zeta,:}J = Z_{\zeta,:}$ where $\zeta = \cup_{k=1}^{K} (\zeta_{k} \setminus S_{k})$ and 
    \begin{multline}\label{eq:weak-cluster-rate-lap}
        \sum_{k=1}^{K} \frac{|S_{k}|}{N_{k}} \lesssim \\
        \frac{\log N }{N \alpha_{N}^{3/2}\lambda_{K}^{2}(B)}  +
         \left\{\frac{1}{\rho^{2} \lambda_{K}^{2}(B) }\right\}^{2} \left\{ \frac{1}{(1 - \rho^{2})^{3}} \frac{K^{2}}{T} +  \left(1 + \frac{1}{ \sqrt{1 - \rho^{2}}}\right)^{2}\frac{KN}{T}\right\} \log(T+N). 
    \end{multline}
\end{theorem}

\subsubsection{Exact recovery}
 Theorem~\ref{thm:strong-cov-estimation} holds in the Laplacian setting; since the proof follows the same technique as the adjacency setting, we omit it for conciseness. 
 Note that in the Laplacian setting, $U \in \mathbb{R}^{N \times K}$ is the matrix of $K$ leading eigenvectors of $\Gamma$ and $\Phi_{L}$, simultaneously.  By the triangle inequality,
\begin{align}
    \min_{R \in \mathcal{O}(K)}\lVert U_{S} - \tilde{U}_{P}R \rVert_{2,\infty} \leq \min_{R \in \mathcal{O}(K)}\lVert U_{S} - UR\rVert_{2,\infty} + \min_{R \in \mathcal{O}(K)}\lVert U - \tilde{U}_{P}R\rVert_{2,\infty}.
\end{align}
Under Assumptions~\ref{ass:regularity} and~\ref{ass:regularity-laplacian}, by \citet[][Theorem~4.1]{ke2025optimal}, with probability at least $1-N^{-4}$,
\begin{align}
    \min_{R \in \mathcal{O}(K)}\lVert U - \tilde{U}_{P}R\rVert_{2,\infty} \lesssim  \frac{\log^{1/2}N}{N\alpha_{N}^{3/4}\lambda_{K}(B)}.
\end{align}
Combining this with Theorem~\ref{thm:strong-cov-estimation} yields a high probability bound on $\min_{R \in \mathcal{O}(K)}\lVert U_{S} - \tilde{U}_{P}R \rVert_{2,\infty}$. This allows us to prove that the $k$-means algorithm gives exact clustering in the regime $T \gtrsim N^{2}\log^{2}N$, stated in Theorem~\ref{thm:strong-clustering-lap}.

\begin{theorem}[Strong clustering for Laplacian model] \label{thm:strong-clustering-lap}
Instate the assumptions of Theorem~\ref{thm:dist_U_S_U_P-lap} and Theorem~\ref{thm:strong-cov-estimation}. If $T > CN^{2}\log^{2}N$ for some $C>0$, then for any $\nu \geq 8$, with probability at least $1 - (N+T)^{-\nu}$, 
    $\sup_{1 \leq i \leq N} \indfun\{ \hat{z}_{i} \neq z_{i} \} = 0.$
\end{theorem}

\subsection{Relation to existing results}
    The weak clustering rate in Theorem~\ref{thm:kmeans-lap} is of the same order as the weak clustering result given by \citet[][Theorem~1]{gudhmundsson2021detecting}. Their ``blockbuster'' algorithm for subspace recovery 
    utilises a different approach 
    compared to this work: 
    instead of applying spectral methods directly to the sample covariance matrix, \citet{gudhmundsson2021detecting} first estimate $\Phi$ using ordinary least squares (OLS) and then apply spectral methods to the OLS estimate $\hat{\Phi}$. One downside of this approach is the need to estimate the $N^{2}$ parameters of $\Phi$, which may be prohibitive in high-dimensional settings. To remedy this, one might consider regularised VAR methods \citep[see][for example]{basu2015regularized}; however these methods require the tuning of hyperparameters. 

    The weak clustering result of \citet[][Theorem~2]{brownlees2022community} is also of the same order in $N$ and $T$ as Theorem~\ref{thm:kmeans-lap}. Similar to our work, \citet{brownlees2022community} also have two layers of randomness: 
    the first arising from the difference between the realisation of the graph and its population parameters; the second arising from the difference between the sample covariance matrix and the population covariance matrix due to sampling uncertainty in the data $\bm{X}_{1},\dots,\bm{X}_{T}$. However, \citet{brownlees2022community} directly assume a block structure within the precision matrix, which they show is equivalent to a contemporaneous regression model, in contrast to the network time series models considered in this work.

     The rate in Theorem~\ref{thm:kmeans-lap} depends explicitly on the spectral radius $\rho$ of the VAR coefficient matrix. This is a direct consequence of using the matrix Bernstein inequality of \citet{banna2016bernstein} as opposed to using the scalar Bernstein inequality of \citet{merlevede2011bernstein} combined with an epsilon-net and union bound argument as done in \citet{gudhmundsson2021detecting} and \citet{brownlees2022community}. 

\subsection{General error covariance, weighted directed graphs}\label{subsec:general-graphs}

This section considers $\ell_{2}$ and $\ell_{2,\infty}$-norm bounds (used in  weak and exact recovery, respectively) in the cases of weighted directed graphs and general noise covariances $\Sigma$. For ease of presentation, we focus on the adjacency model, although analogue results hold for the Laplacian model. We refer to \citet[][Definition 3]{gallagher2023spectral} for a definition of a weighted SBM, \citet[][Definition 1]{rohe2016co} for a definition of a directed SBM, and \citet[][Definition 1]{gudhmundsson2025detecting} for a definition of a weighted directed SBM. In the case of directed SBMs, the left and right singular vectors of $P$ correspond to two types of block membership, commonly referred to as giver and receiver groups, respectively. 

\citet[][Proposition 4.1]{martin2024nirvar_arxiv} shows that for homoskedastic noise and undirected graphs, the eigenvectors of $\Gamma$ and $\Phi$ are the same. Relaxing these two assumptions means that 
\citet[][Proposition 4.1]{martin2024nirvar_arxiv} 
may not hold. 
As such, the eigenvectors of $\Gamma$ are not necessarily equal to the left singular vectors of $\Phi$. For this section, we write the top $K$ eigenvectors of $\Gamma$ as $U_{\Gamma}$ and the singular value decomposition of $\Phi$ as $\Phi = U_{\Phi} \Sigma_{\Phi} V_{\Phi}^{\prime} + U_{\Phi,\perp} \Sigma_{\Phi,\perp} V_{\Phi,\perp}^{\prime}$ where $\Sigma_{\Phi} = \mathrm{diag}\{\sigma_{1}(\Phi),\dots,\sigma_{K}(\Phi)\}$ is the diagonal matrix containing the top $K$ largest in magnitude singular values of $\Phi$. The triangle inequality implies 
\begin{align}
    \mathrm{dist}(U_{S},U_{P}) \leq \mathrm{dist}(U_{S},U_{\Gamma}) + \mathrm{dist}(U_{\Gamma},U_{\Phi}) + \mathrm{dist}(U_{\Phi},U_{P}).
\end{align}
A high probability bound on $\mathrm{dist}(U_{S},U_{\Gamma})$ can be obtained by modifying the proof of Lemma~\ref{lemma:S_T-Gamma}. In particular, the variance parameters in the matrix Bernstein inequality can be bound solely in terms of the lag-$h$ autocovariances $\Gamma(h) = \mathbb{E}(\bm{X}_{t}\bm{X}_{t+h}^{\prime})$. 
Lemma~\ref{lemma:S_T-Gamma-general-errors} below controls the $\ell_{2}$ covariance estimation error in the setting of weighted directed graphs and general noise covariance $\Sigma$ (note that $\bm{\epsilon}_{t}$ is still assumed to be Gaussian).

\begin{lemma}\label{lemma:S_T-Gamma-general-errors}
Let Assumption \ref{ass:regularity} hold and suppose that 
    \begin{equation}\label{eq:lemma_T_condition-ge}
        T > C \frac{\gamma(c_{\beta},T)^{2}}{1 - \rho^{2}}\log^{3} (T+N)
    \end{equation}
for some sufficiently large constant $C >0$ where $\gamma(c_{\beta},T)$ is defined in \eqref{eq:gamma}.
Then for any $\nu \geq 8$, with probability exceeding $1 - (T+N)^{-\nu}$, 
    \begin{align} 
        \left\lVert S_{T} - \Gamma \right\rVert_{2} &\lesssim \frac{\lVert \Sigma \rVert_{2}}{(1 - \rho^{2})^{3/2}} \left( \sqrt{\frac{K}{T}}  +\sqrt{\frac{N}{T}} \right) \log^{1/2}(T+N).
    \end{align}
\end{lemma}

Lemma~\ref{lemma:S_T-Gamma-general-errors} has the same rate in $N$ and $T$ as in Lemma~\ref{lemma:S_T-Gamma}. Note however, that the constants multiplying the $\sqrt{N/T}$ term in Lemma~\ref{lemma:S_T-Gamma} are smaller than those in Lemma~\ref{lemma:S_T-Gamma-general-errors}. This is because a general noise covariance $\Sigma$ causes the leading $K$ eigenvectors of $\Gamma(h)$ to be different from $U_{\Gamma}$. Therefore, an unbiased estimator of $U_{\Gamma}$ would need to account for this misalignment. 

Given Lemma~\ref{lemma:S_T-Gamma-general-errors}, a high probability bound on $\mathrm{dist}(U_{S},U_{\Gamma})$ can be obtained via the Davis-Kahan sin$\Theta$ theorem. Since this is the exact same proof technique as Theorem~\ref{thm:dist_U_S_U_P}, we omit an explicit derivation for conciseness. A high probability bound on $\mathrm{dist}(U_{\Phi},U_{P})$ 
is obtained by combining \citet[][Theorem OA.1]{gudhmundsson2025detecting} with the Davis-Kahan sin$\Theta$ theorem. In particular, \citet[][Theorem~OA.1]{gudhmundsson2025detecting} proves that for a weighted directed SBM, $\lVert A - P \rVert_{2} \lesssim \sqrt{\alpha_{N}N\log N }$ with high probability. Therefore, by the Davis-Kahan sin$\Theta$ theorem and Lemma~\ref{lemma:eigs}, $\mathrm{dist}(U_{\Phi},U_{P}) \lesssim \log^{1/2}N/\sqrt{\alpha_{N}N}$ with high probability. 

It remains to control $\mathrm{dist}(U_{\Gamma},U_{\Phi})$. Under the assumption of homoskedastic noise, we can bound $\mathrm{dist}(U_{\Gamma},U_{\Phi})$ in terms of $\lVert \Phi \Phi^{\prime} - \Phi^{\prime}\Phi \lVert_{2}$, which is a measure of how far $\Phi$ is from being a normal matrix. Note that Lemma~\ref{lemma:singular-angle} below holds for weighted directed graphs but, unlike all other results in this section, not for heteroskedastic noise.   

\begin{lemma}\label{lemma:singular-angle}
    Assume $\Sigma \propto I_{N}$. For any $s>0$, with probability at least $1 - N^{-s}$, 
    \begin{align}
        \mathrm{dist}(U_{\Gamma},U_{\Phi}) \lesssim \frac{\lVert \Phi \Phi^{\prime} - \Phi^{\prime}\Phi \lVert_{2}}{(1-\rho^{2})^{3} \lambda_{K}^{2}(B)}. 
    \end{align}
\end{lemma}

Lemma~\ref{lemma:singular-angle} shows that if $\Phi$ is a normal matrix, the subspace estimation error $\mathrm{dist}(U_{\Gamma},U_{\Phi})$ is zero. Also, $\mathrm{dist}(U_{\Gamma},U_{\Phi})$ is inversely proportional to $\lambda_{K}(B)$, which measures the separability of the underlying communities. 

The entry-wise $\ell_{2,\infty}$ error between $U_{S}$ and $U_{\Gamma}$ can be controlled in the case of weighted directed graphs and general noise covariance via Theorem~\ref{thm:hetero-strong-cov} below. The proof of Theorem~\ref{thm:hetero-strong-cov} follows that of Theorem~\ref{thm:strong-cov-estimation}, with a modified bound on the Bernstein variance parameter. The bound in Theorem~\ref{thm:hetero-strong-cov} achieves the same rate in $N$ and $T$ as 
Theorem~\ref{thm:strong-cov-estimation}, but has larger constants due to the bias incurred when considering heteroskedastic noise.

\begin{theorem}[Covariance estimation under general noise and weighted directed graphs] \label{thm:hetero-strong-cov}
    Let Assumption~\ref{ass:regularity} hold. 
    There exists an orthogonal matrix $W_{U}$ such that, for any $\nu \geq 8$, with probability $1 - (T+N)^{-\nu}$,
    \begin{align}\label{eq:general-strong-cov-bound}
    \|U_{S} - U_{\Gamma} W_{U}\|_{2 ,\infty} \lesssim &\   \frac{1}{\rho^{2} \lambda_{K}^{2}(B) } \Bigg{[} \frac{\lVert \Sigma \rVert_{2}}{(1 - \rho^{2})^{3/2}} \left\{ \sqrt{\frac{K^{3}}{NT}} +   \sqrt{\frac{K^{2}}{T}} \right\} \log^{1/2}T \\ 
    &+  \frac{\lVert \Sigma \rVert_{2}}{(1 - \rho^{2})^{3/2}} \left\{\frac{\sqrt{NK}}{T} +   \frac{N}{T}\right\} \log T \log^{1/2}(T+N) \\
    &+  \frac{\lVert \Sigma \rVert_{2}}{(1 - \rho^{2})^{3/2}} \left\{ \sqrt{\frac{K}{T}} +   \sqrt{\frac{N}{T}}\right\} \log^{1/2}(T+N) \\
    &+ \frac{\lVert \Sigma \rVert_{2}^{2}}{(1 - \rho^{2})^{3}} \left\{  \frac{K}{T} +  \frac{\sqrt{NK}}{T} +  \frac{N}{T}\right\} \sqrt{\frac{K}{N}}\log(T+N) \Bigg{]}.
\end{align}
\end{theorem}
For weighted, undirected graphs, the $\ell_{2,\infty}$-norm error between $U_{\Phi}$ and $U_{P}$ can be controlled by \citet[][Theorem~1]{gallagher2023spectral}. As far as we are aware, fine-grained analysis of adjacency spectral embedding for weighted directed SBMs has not been explicitely addressed in the literature. Extending \citet[][Theorem~1]{gallagher2023spectral} to the case of weighted directed SBMs would be feasible by following a similar procedure to the proof of \citet[][Theorem~1]{corneck2026spectral}.

\section{Numerical studies}\label{sec:sim}

\begin{figure}[t]
\centering
\begin{subfigure}{.49\textwidth}
  \centering
  \includegraphics[width=.95\linewidth]{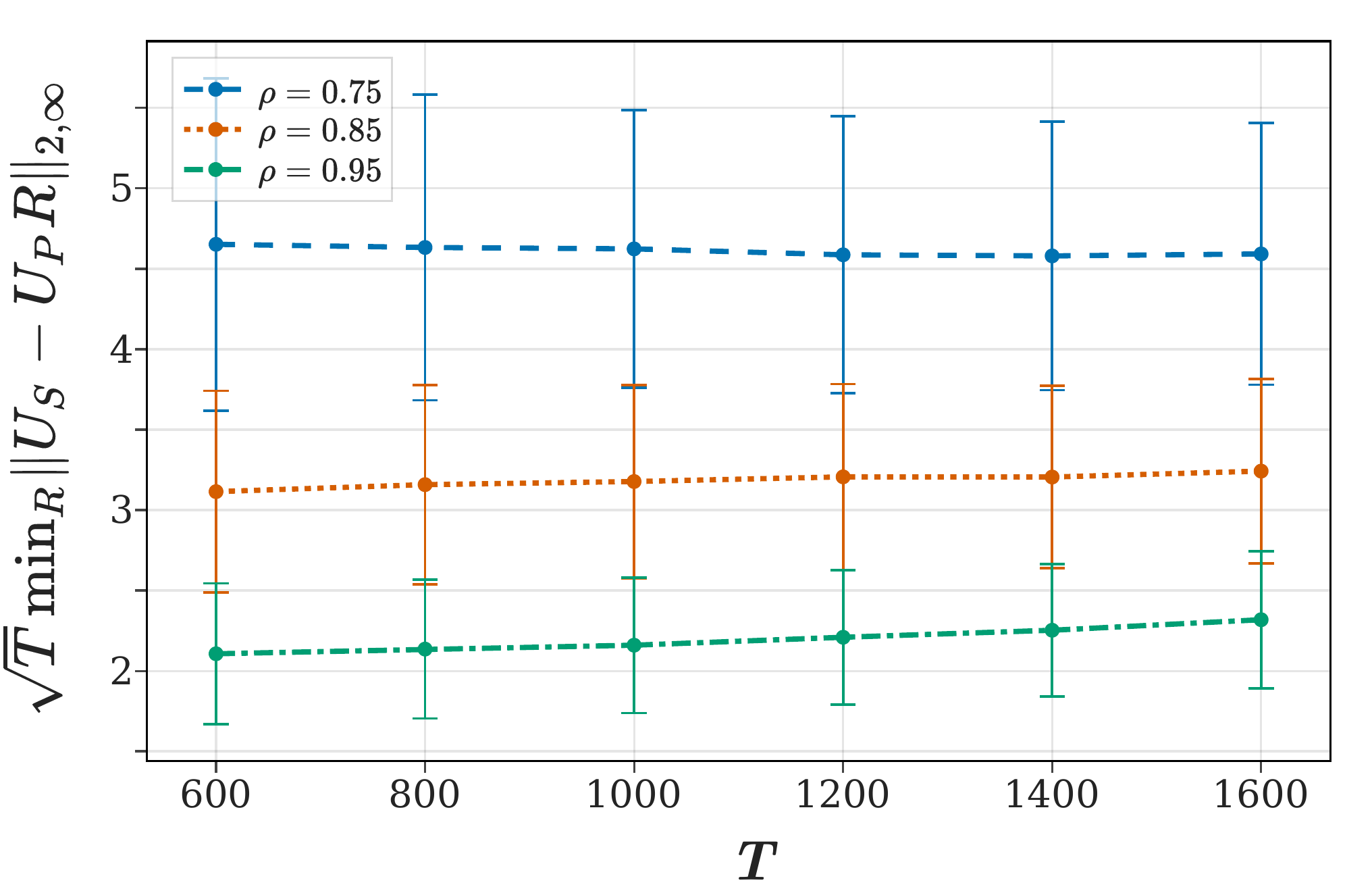}
  \caption{Homoscedastic $\Sigma$}
  \label{fig:varyT-sub1}
\end{subfigure}
\begin{subfigure}{.49\textwidth}
  \centering
  \includegraphics[width=.95\linewidth]{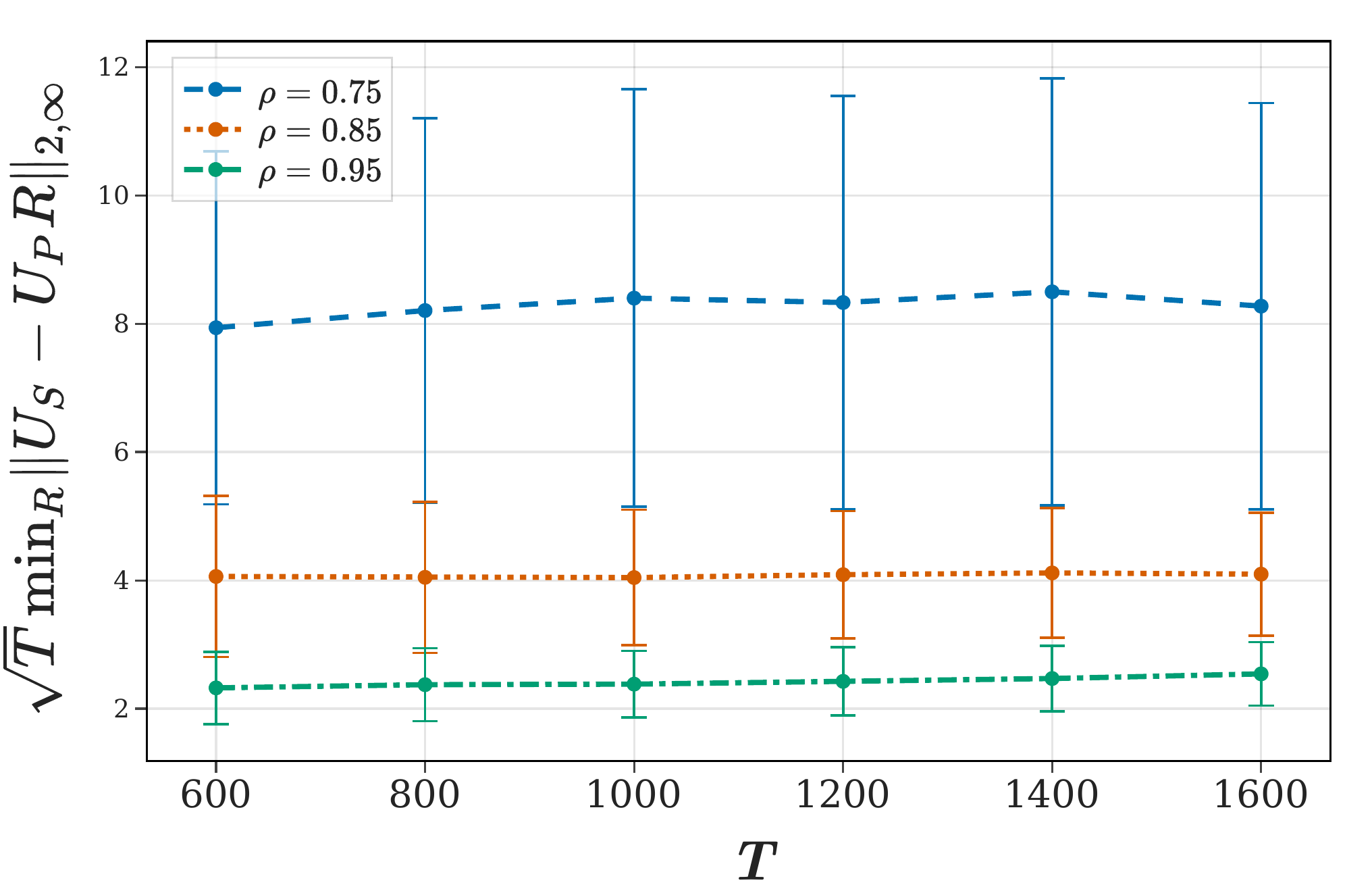}
  \caption{Heteroscedastic $\Sigma$}
  \label{fig:varyT-sub2}
\end{subfigure}
\begin{subfigure}{.49\textwidth}
  \centering
  \includegraphics[width=.95\linewidth]{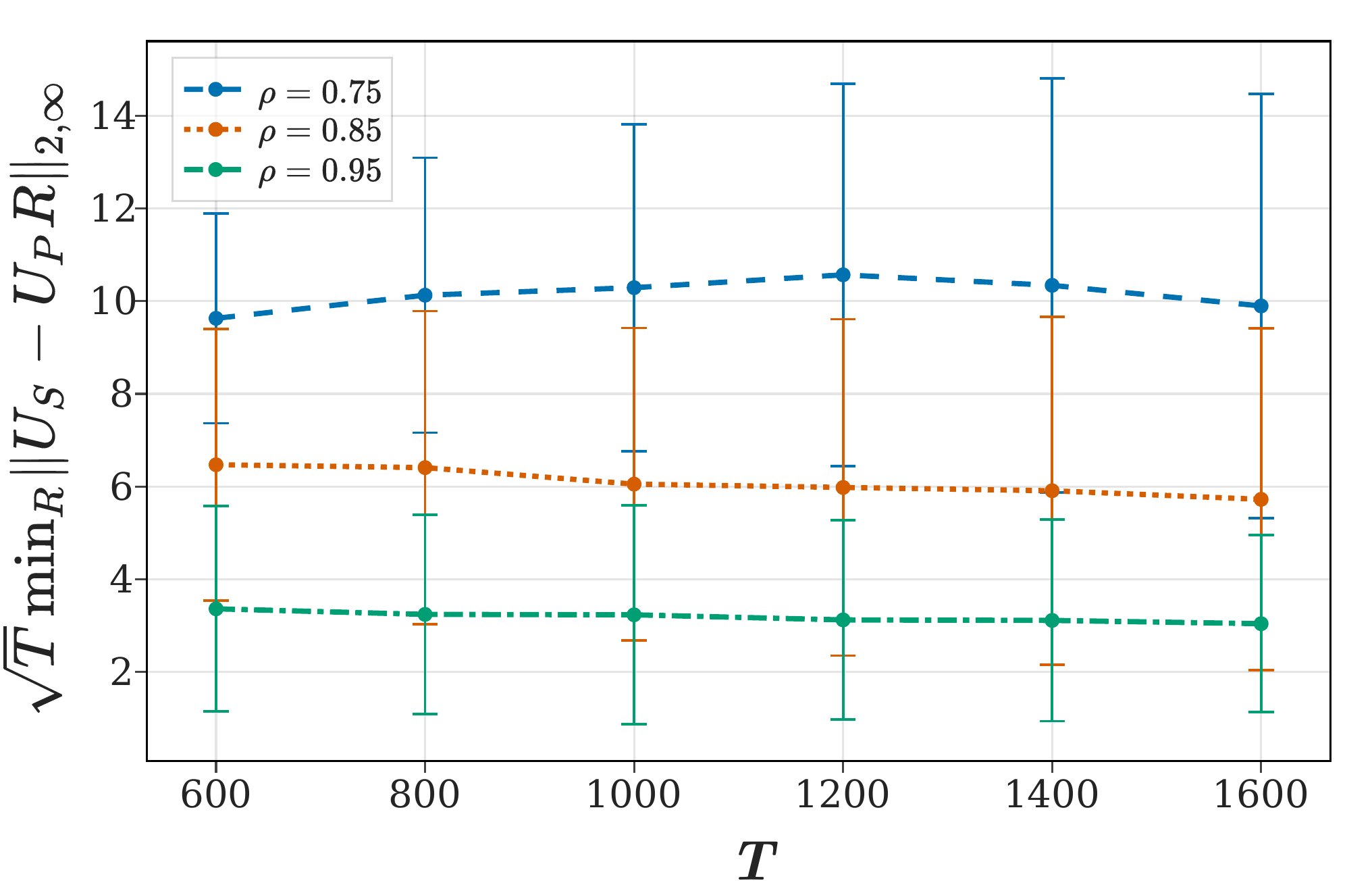}
  \caption{Heavy-tailed $\Sigma$}
  \label{fig:varyT-sub3}
\end{subfigure}
\begin{subfigure}{.49\textwidth}
  \centering
  \includegraphics[width=.95\linewidth]{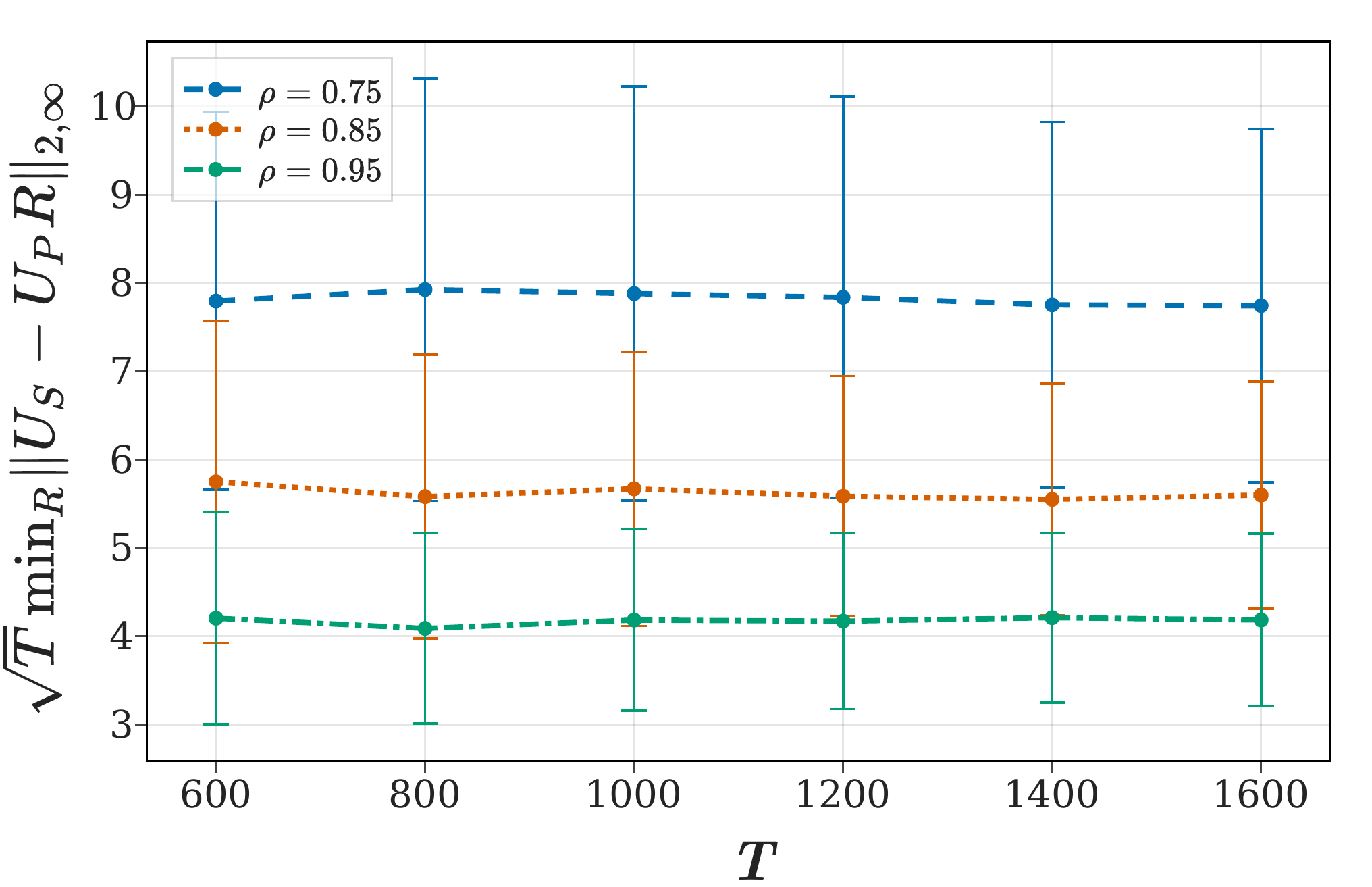}
  \caption{Weighted directed $\Phi$}
  \label{fig:varyT-sub4}
\end{subfigure}
\caption{Subspace recovery rate as a function of the number of time series observations $T$. The rate for different values of the VAR spectral radius $\rho$, which controls the SNR, are shown.}
\label{fig:varyT}
\end{figure}

We conducted simulation studies assessing the finite sample properties of the spectral estimator considered in this work, examining its sensitivity to different sparsity and data dependence regimes, and demonstrating its robustness to weighted, directed graphs and heavy-tailed noise. 
First, we simulated data from the adjacency model in Definition~\ref{def:model-adj}, with $N = 50$, $K=2$, $N_{1} = N_{2}=25$, $\alpha_{N} = 1$, $\Sigma = I_{N}$, and $B = [[0.9, 0.05],[0.05,0.9]]$. Following Algorithm~\ref{alg:spectral}, we computed $U_{S}$ and $\sqrt{T}\min_{R \in \mathcal{O}(K)}\lVert U_{S} - U_{P}R \rVert_{2,\infty}$, calculated using Procrustes alignment. This was repeated $1500$ times, with the mean and standard deviation over the 1500 replications being recorded. The results are reported in Figure~\ref{fig:varyT}. Figure~\ref{fig:varyT}(\subref{fig:varyT-sub1}) shows that $\sqrt{T}\min_{R \in \mathcal{O}(K)}\lVert U_{S} - U_{P}R \rVert_{2,\infty}$ is constant in $T$, matching the rate given by Theorem~\ref{thm:strong-cov-estimation}. Note that the error bars quantify the standard deviation of the distribution of $\sqrt{T}\min_{R \in \mathcal{O}(K)}\lVert U_{S} - U_{P}R \rVert_{2,\infty}$, not a standard error on a mean.  Figure~\ref{fig:varyT}(\subref{fig:varyT-sub1}) also shows that the value of the constant increases with decreasing $\rho$, and the error bars also increase, since $\rho$ determines the SNR (see Remark~\ref{remark:SNR}). 

Next, we repeated the above simulation under two different regimes for the error matrix: \emph{(i)} heteroscedastic $\Sigma$, with $\Sigma_{ij} = \sigma_{i}^{2}\indfun\{i=j\}$ where $\sigma_{i} \sim \text{Uniform}(0.5,1)$, and \emph{(ii)} heavy-tailed $\Sigma$, where the error matrix is sampled from a multivariate Student’s-$t$ distribution with 5 degrees of freedom.
Figure~\ref{fig:varyT}(\subref{fig:varyT-sub2}) and (\subref{fig:varyT-sub3}) show that the rate given by Theorem \ref{thm:strong-cov-estimation} continues to hold in this simulation study for both heteroscedastic and heavy-tailed error processes, respectively. 
Additionally, we repeated the simulation study using a weighted and directed SBM adjacency matrix, where the weight on edge $(i,j)$ is sampled from $\mathcal{N}(\mu_{i},1)$ with $\mu_{i} = 1$ if $z_{i}=1$ and $\mu_{i}=-1$ if $z_{i} = 2$. Figure~\ref{fig:varyT}(\subref{fig:varyT-sub4}) shows that the rate given by Theorem~\ref{thm:strong-cov-estimation} also holds for this choice of adjacency matrix. 

\begin{figure}[t]
\begin{center}
\begin{subfigure}{.49\textwidth}
  \centering
  \includegraphics[width=.95\linewidth]{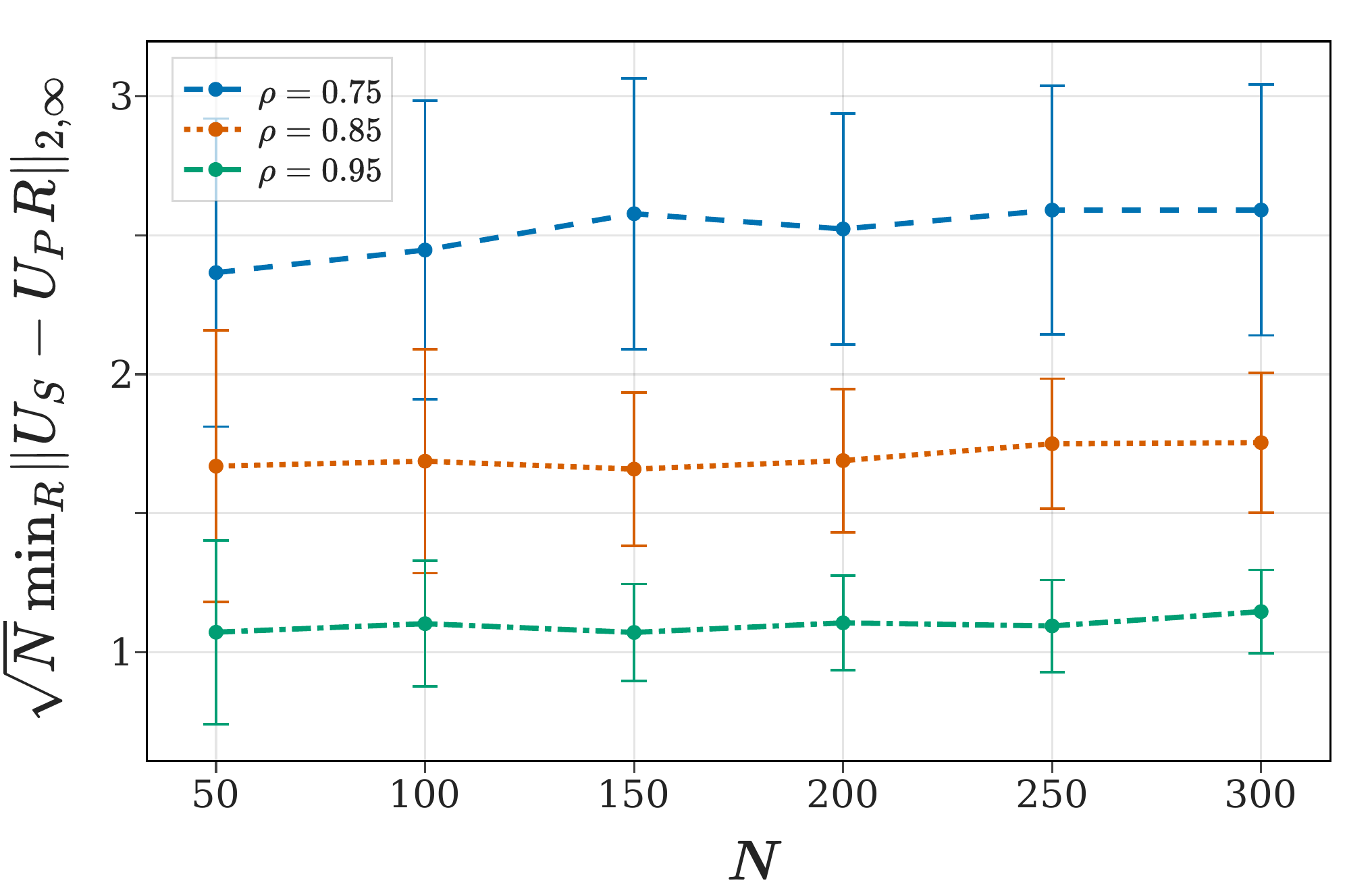}
  \caption{Varying SNR}
  \label{fig:varyN-sub1}
\end{subfigure}
\begin{subfigure}{.49\textwidth}
      \centering
      \includegraphics[width=0.95\linewidth]{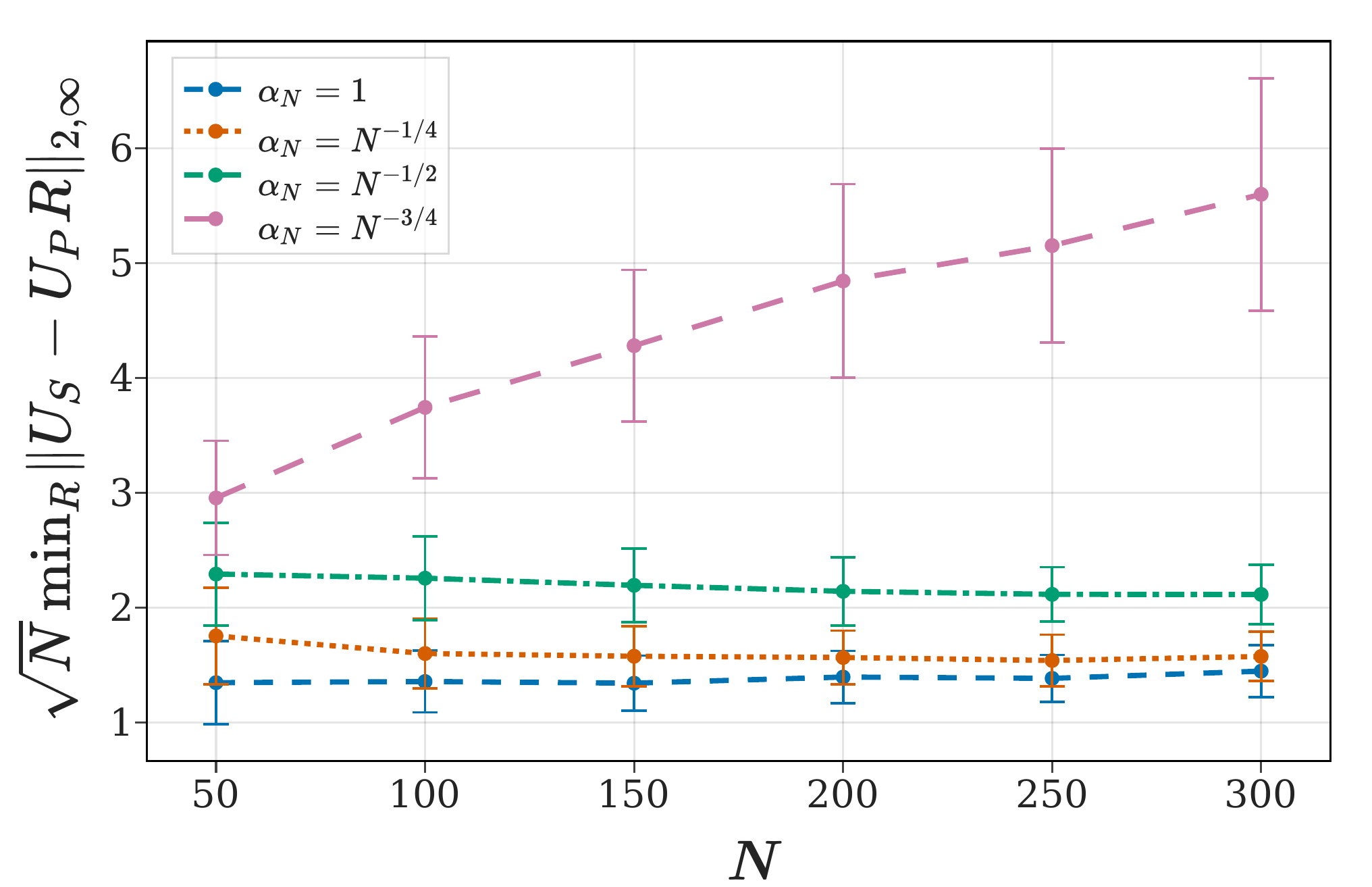}
  \caption{Varying graph sparsity $\alpha_{n}$}
  \label{fig:varyN-sub2}
\end{subfigure}
\caption{Subspace recovery rate as a function of the number of graph vertices $N$. (a) The rate for different values of the VAR spectral radius $\rho$, which controls the SNR. (b) The rate for fixed $\rho$ and variable graph sparsity $\alpha_{N}$.} 
\label{fig:varyN} 
\end{center}
\end{figure} 

We performed further simulations to study $\min_{R \in \mathcal{O}(K)}\lVert U_{S} - U_{P}R \rVert_{2,\infty}$ as a function of $N$. Again, we set $N_{1} = N_{2} = N/2$, $\alpha_{N} = 1$, $\Sigma = I_{N}$, and $B = [[0.9, 0.05],[0.05,0.9]]$. We chose $T = 4N$ to align with the requirement given by \eqref{eq:T_cond_combined}. Figure~\ref{fig:varyN}(\subref{fig:varyN-sub1}) shows that $\sqrt{N}\min_{R \in \mathcal{O}(K)}\lVert U_{S} - U_{P}R \rVert_{2,\infty}$ is constant as a function of $N$. This indicates that the bound in Theorem~\ref{thm:strong-cov-estimation} may be tightened by a factor of $\sqrt{N}$, specifically the terms on the second and third lines of \eqref{eq:strong-cov-bound}. In the i.i.d. setting, these terms could be tightened by a factor of $\sqrt{N}$ using leave-one-out analysis (see the discussion at the end of Section~\ref{sec:adjacency-specification}). We repeated the simulation with varying levels of graph sparsity $\alpha_{N}$ (and fixed $\rho = 0.9$). Figure~\ref{fig:varyN}(\subref{fig:varyN-sub2}) shows that for $\alpha_{N} = N^{-3/4}$, $\sqrt{N}\min_{R \in \mathcal{O}(K)}\lVert U_{S} - U_{P}R \rVert_{2,\infty}$ is no longer constant in $N$ but rather increases with $N$. This is expected since $N^{-3/4} < \log(N)/N$ for $N \in [100,300]$, which violates the minimum sparsity condition in Assumption~\ref{ass:regularity}. 

\section{Discussion}

There are a number of avenues for future work. Firstly, obtaining tight $\ell_{2,\infty}$ subspace estimation bounds in the setting of both temporal and cross-sectional dependence remains an open problem. \citet[][]{fan2024can} tackled the case of contemporaneous cross-sectional dependence by factorising out the cross-sectional covariance matrix and employing leave-one-out techniques on the residual matrix. However, when there is cross-sectional dependence at all temporal lags (as in our setting), adapting the factorisation approach of \citet[][]{fan2024can} appears fruitless to us. Rather, a new proof technique other than leave-one-out type analysis might be required. 

Secondly, obtaining tight subspace estimation error bounds in the presence of heteroskedastic noise might require modifying the spectral clustering algorithm used in this paper via techniques similar to heteroskedastic PCA \citep[][]{yan2024inference}, for example. Thirdly, extensions to degree-corrected and mixed membership SBMs would increase the flexibility of the models considered in this paper.

Finally, in the setting of i.i.d. data and observable adjacency matrices, there exists a multitude of distributional results for spectral methods applied to SBMs \citep[see][for an overview]{agterberg2023overview}. The consistency rates proved in this paper are a first step towards distributional results in the setting of dependent data and an unobserved SBM adjacency matrix. A key motivation for obtaining such distributional results is to develop tests for whether two time series belong to the same block of the underlying SBM. 



\begin{appendix}

\section{Proofs of main results}\label{sec:appendix}
\subsection{Proof of Lemma~\ref{lemma:truncated-Bernstein}}\label{sec:trunc-bern-proof}

For any matrix $X \in \mathbb{R}^{d_{1} \times d_{2}}$, consider the dilation \citep[see, for example] [Section 2.1.17]{tropp2015introduction}
\begin{align}
    S(X) &=
\begin{pmatrix}
0 & X \\
X^{\prime} & 0
\end{pmatrix} \in \mathbb{R}^{(d_{1} + d_{2}) \times (d_{1} + d_{2})},
\end{align} 
and note that $\lVert S(X) \rVert_{2} = \lVert X \rVert_{2}$. Let $M_{t}^{\circ} = M_{t} - \mathbb{E}(M_{t})$ and note that   $S(M_{t}^{\circ})$ is symmetric and zero-mean, and equation \eqref{def:mv} 
can be written as
\begin{align}
    v^{2} = \sup_{\tau \subseteq [T]} \frac{1}{\mathrm{Card}(\tau)} \lambda_{1}\left\{ \mathbb{E}\left( \sum_{t \in \tau } S(M_{t}^{\circ}) \right)^{2} \right\}.
\end{align}
Define the truncated counterpart to $M_{t}^{\circ}$ as
\begin{align}
    \tilde{M}_{t} \coloneqq \begin{cases}
        M_{t}^{\circ} & \text{if}\ \lVert M_{t}^{\circ} \rVert_{2} \leq L, \\
        0 & \text{otherwise}.
    \end{cases}
\end{align}
As $\tilde{M}_{t}$ is a measurable function of $M_{t}$, the sequence $(\tilde{M}_{t})_{t \in [T]}$ is also geometrically $\beta$-mixing with constant $c_{\beta}$.
Since
\begin{align}
    \mathbb{E}\left( \sum_{t \in \tau } S(\tilde{M}_{t} - \mathbb{E}[\tilde{M}_{t}])\right)^{2} \preceq \mathbb{E}\left( \sum_{t \in \tau } S(\tilde{M}_{t})\right)^{2},
\end{align}
we have
\begin{align}
    \left\lVert \mathbb{E}\left( \sum_{t \in \tau } S(\tilde{M}_{t} - \mathbb{E}[\tilde{M}_{t}])\right)^{2} \right\rVert_{2} \leq \left\lVert \mathbb{E}\left( \sum_{t \in \tau } S(\tilde{M}_{t})\right)^{2} \right\rVert_{2}.
\end{align}
Define the events $\mathcal{E}_{t,s} \coloneqq \{\lVert M_{t}^{\circ} \rVert_{2} \leq L\} \cap \{ \lVert M_{s}^{\circ} \rVert_{2} \leq L\}$ and   $\mathcal{E}_{t,s}^{\prime}$,  its complement. Then 
\begin{align}\label{eq:trunc-variance-two-terms}
    \left\lVert \mathbb{E}\left( \sum_{t \in \tau } S(\tilde{M}_{t})\right)^{2} \right\rVert_{2} 
    &=\left\lVert \sum_{t, s \in \tau } \mathbb{E}[ S(M_{t}^{\circ}) S(M_{s}^{\circ})] - \mathbb{E}[ S(M_{t}^{\circ}) S(M_{s}^{\circ}) \indfun\{\mathcal{E}_{t,s}^{\prime}\}]  \right\rVert_{2} \\ 
    &\leq \left\lVert \sum_{t, s \in \tau } \mathbb{E}[ S(M_{t}^{\circ}) S(M_{s}^{\circ})] \right\rVert + \sum_{t, s \in \tau} \left\lVert \mathbb{E}[ S(M_{t}^{\circ}) S(M_{s}^{\circ}) \indfun\{\mathcal{E}_{t,s}^{\prime}\}]  \right\rVert_{2}.
\end{align}
By Jensen's inequality and the Cauchy-Schwarz inequality,
\begin{align}\label{eq:jensen-CS}
    \left\lVert \mathbb{E}[ S(M_{t}^{\circ}) S(M_{s}^{\circ}) \indfun\{\mathcal{E}_{t,s}^{\prime}\}]  \right\rVert_{2} &\leq \mathbb{E}\left(\left\lVert S(M_{t}^{\circ}) S(M_{s}^{\circ}) \right\rVert_{2} \indfun\{\mathcal{E}_{t,s}^{\prime}\}\right) \\
    &\leq \sqrt{\mathbb{E}\left( \left\lVert S(M_{t}^{\circ}) \right\rVert_{2}^{4} \right) \mathbb{P}(\mathcal{E}_{t,s}^{\prime})}.
\end{align}
The assumption $\mathbb{P}\{\lVert M_{t} - \mathbb{E}(M_{t}) \rVert_{2} \geq L \} \leq q_{0}$ implies $\mathbb{P}(\mathcal{E}_{t,s}^{\prime}) \leq q_{0}$. Then the assumption $q_{0} \lesssim 1/\{\mathbb{E}(\lVert M_{t} - \mathbb{E}[M_{t}] \rVert_{2}^{4}) (\mathrm{Card}(\tau))^{4}\}$ and \eqref{eq:jensen-CS} imply that $\left\lVert \mathbb{E}[ S(M_{t}^{\circ}) S(M_{s}^{\circ}) \indfun\{\mathcal{E}_{t,s}^{\prime}\}]  \right\rVert_{2} \lesssim 1/(\mathrm{Card}(\tau))^{2}$. Since the summation in the second term on the last line of \eqref{eq:trunc-variance-two-terms} contains $(\mathrm{Card}(\tau))^{2}$ such terms, it follows that the second term on the last line of \eqref{eq:trunc-variance-two-terms} is $O(1)$ and therefore
\begin{align}
    \tilde{v}^{2} \coloneqq \sup_{\tau \subseteq [T]} \frac{1}{\mathrm{Card}(\tau)} \lambda_{1}\left\{ \mathbb{E}\left( \sum_{t \in \tau } S(\tilde{M}_{t} - \mathbb{E}[\tilde{M}_{t}]) \right)^{2} \right\} \lesssim v^{2}.
\end{align}
The matrices $S(\Tilde{M}_{t} - \mathbb{E}[\tilde{M}_{t}])$ are self-adjoint, zero-mean, and have a maximum eigenvalue which is bounded almost surely by $L$. Therefore, by \citet[][Theorem 1]{banna2016bernstein} we have 
\begin{equation}
    \label{eq:Bernstein-B-lemma} 
    \mathbb{P}\left( \left\lVert \sum_{t = 1}^{T} S(\Tilde{M}_{t} - \mathbb{E}[\tilde{M}_{t}]) \right\rVert_{2} \geq \alpha \right) \leq (d_{1} + d_{2})\exp\left(-\frac{C \alpha^{2}}{v^{2}T + c_{\beta}^{-1}L^{2} + \alpha L_{1}\gamma(c_{\beta},T)}\right).
\end{equation}
Letting $M^{\circ} \coloneqq \sum_{t = 1}^{T} M_{t}^{\circ}$ and $\tilde{M} \coloneqq \sum_{t = 1}^{T} \tilde{M}_{t}$ gives
\begin{align}
    \mathbb{P}\left( \left\lVert  S(M^{\circ}) \right\rVert_{2} \geq x \right) = &\ \mathbb{P}\left( \left\lVert S(M^{\circ}) \right\rVert_{2} \geq x \middle| S(M^{\circ}) = S(\tilde{M})  \right) \mathbb{P}\left(  S(M^{\circ}) = S(\tilde{M})  \right) \\ &+ \mathbb{P}\left( \left\lVert S(M^{\circ}) \right\rVert_{2} \geq x \middle| S(M^{\circ}) \neq S(\tilde{M})  \right) \mathbb{P}\left(  S(M^{\circ}) \neq S(\tilde{M})  \right). 
\end{align}
By assumption, $\mathbb{P}(S(M_{t}^{\circ}) \neq S(\tilde{M}_{t})) \leq q_{0}$. Applying a union bound gives $\mathbb{P}(S(M^{\circ}) \neq S(\tilde{M})) \leq T q_{0}$. Therefore
\begin{align}
    \mathbb{P}\left( \left\lVert  S(M^{\circ}) \right\rVert_{2} \geq x \right) \leq \mathbb{P}\left( \left\lVert S(M^{\circ}) \right\rVert_{2} \geq x \middle| S(M^{\circ}) = S(\tilde{M})  \right) + T q_{0}.
\end{align}
By assumption, $\lVert \mathbb{E}(\tilde{M}_{t})\rVert_{2} \leq q_{1}$ and by a union bound, $\lVert \mathbb{E}(S(\tilde{M}))\rVert_{2} \leq T q_{1}$. By the triangle inequality, 
 \begin{align}
     \left\lVert S(M^{\circ}) \right\rVert_{2} \leq \left\lVert S(M^{\circ}) - \mathbb{E}(S(\tilde{M}))  \right\rVert_{2} + \left\lVert \mathbb{E}(S(\tilde{M})) \right\rVert_{2} \leq \left\lVert S(M^{\circ}) - \mathbb{E}(S(\tilde{M}))  \right\rVert_{2} + T q_{1}.
 \end{align}
Thus,
\begin{align}
\label{eq:A-prob-equal-lemma}
    &\mathbb{P}\left( \left\lVert S(M^{\circ}) \right\rVert_{2} \geq x\ \middle| S(M^{\circ}) = S(\tilde{M}) \right) \\ &\qquad\qquad\qquad\leq \mathbb{P}\left( \left\lVert S(M^{\circ}) - \mathbb{E}(S(\tilde{M})) \right\rVert_{2} + T q_{1} \geq x\ \middle|\ S(M^{\circ}) = S(\tilde{M})  \right) \\ 
    &\qquad\qquad\qquad= \mathbb{P}\left( \left\lVert S(\tilde{M}) - \mathbb{E}(S(\tilde{M})) \right\rVert_{2}  \geq x - T q_{1}\ \middle|\ S(M^{\circ}) = S(\tilde{M})  \right).
\end{align}
The last line of \eqref{eq:A-prob-equal-lemma} is bounded by \eqref{eq:Bernstein-B-lemma}. Recalling that $\lVert S( M^{\circ}) \rVert_{2} = \lVert \sum_{t = 1}^{T} M_{t} - \mathbb{E}(M_{t}) \rVert_{2}$ and setting $\alpha = x - T q_{1}$,
\begin{equation}
        \mathbb{P}\left( \left\lVert \sum_{t = 1}^{T} M_{t} - \mathbb{E}(M_{t}) \right\rVert_{2} \geq x \right) \leq  Tq_{0} + (d_{1} + d_{2})\exp\left(-\frac{C(x-Tq_{1})^{2}}{v^{2}T + c_{\beta}^{-1}L_{1}^{2} + (x-Tq_{1})L_{1}\gamma(c_{\beta},T)}\right),
\end{equation}
which proves the claim. \qed

\subsection{Proof of Lemma~\ref{lemma:eigs}} \label{sec:proof3_4}

    By Ostrowski's theorem \citep[see][Theorem 4.5.9]{horn2012matrix}, for $r \in [K]$,
    \begin{align}
        \lambda_{r}(P) &= \alpha_{N}\lambda_{r}(ZBZ^{\prime}) = \alpha_{N} \theta_{r} \lambda_{r}(B),
    \end{align}
    where $\theta_{r} \in [N_{\text{min}},N_{\max}]$ with $N_{\text{min}}$ and $N_{\text{max}}$ being the smallest and largest block size, respectively. Under Assumption \ref{ass:regularity}, the blocks are balanced and therefore, 
    \begin{align}
        \lambda_{r}(P) \asymp \frac{\alpha_{N} N \lambda_{r}(B)}{K}.
    \end{align}
    By Weyl's inequality \citep[see][Lemma 2.2]{chen2021spectral},
    \begin{align}
        |\lambda_{r}(A) - \lambda_{r}(P)| \leq \lVert A - P \rVert_{2}. 
    \end{align}
    Then, by \citet[][Theorem 5.2]{lei2015consistency}, for any $s>0$, with probability at least $1 - N^{-s}$,
    \begin{align}
        |\lambda_{r}(A) - \lambda_{r}(P)| &\lesssim \sqrt{N \alpha_{N}}.
    \end{align}
    Thus, since $\lambda_{r}(B)$ is bounded away from zero for $r \in [K]$, with probability at least $1 - N^{-s}$,
    \begin{align}
        \lambda_{r}(A) \asymp \frac{\alpha_{N} N \lambda_{r}(B)}{K}.
    \end{align}
    Furthermore, since under Assumption \ref{ass:regularity}, $\lambda_{1}(B) \lesssim K$, then with probability at least $1 - N^{-s}$,
    \begin{align}
        \lambda_{r}(\Phi) = \rho \frac{\lambda_{r}(A)}{\lambda_{1}(A)}  
        \asymp \rho \lambda_{r}(B). 
    \end{align}
    For $r = K+1,\dots,N$, since $P$ is rank $K$, $\lambda_{r}(P) = 0$. Thus with probability at least $1 - N^{-s}$, $|\lambda_{r}(A)| \lesssim \sqrt{N \alpha_{N}}$ and 
    \begin{align}
        \lambda_{r}(\Phi) &\asymp \frac{\rho}{\sqrt{N \alpha_{N}}} 
        \xrightarrow[]{N \to \infty} 0. 
    \end{align}
    Therefore, with probability at least $1 - N^{-s}$ $\lambda_{r}(\Phi) = o(1)$ for $r = K+1,\dots,N$. \qed

\subsection{Proof of Lemma~\ref{lemma:S_T-Gamma}} \label{sec:proof3_7}
In order to prove Lemma~\ref{lemma:S_T-Gamma}, we will bound $\lVert S_{T} - \Gamma\rVert_{2}$ by the sum of three terms, each of which can be bounded in turn using the truncated matrix Bernstein-type inequality given in Lemma \ref{lemma:truncated-Bernstein}. We start by defining $\bm{Y}_{t} \coloneqq U^{\prime}\bm{X}_{t}$ and $\bm{Y}_{t}^{\perp} \coloneqq U_{\perp}^{\prime}\bm{X}_{t}$. Note that since $U^{\prime}\Gamma U = \Lambda_{\Gamma}$ and $U_{\perp}^{\prime}\Gamma U_{\perp} = \Lambda_{\Gamma, \perp}$, then $\bm{Y}_{t} \sim \mathcal{N}(0, \Lambda_{\Gamma})$ and $\bm{Y}_{t}^{\perp} \sim \mathcal{N}(0, \Lambda_{\Gamma,\perp})$. Furthermore, since $U^{\prime} \Phi = \Lambda_{\Phi} U^{\prime}$, the autoregressive equations of $\bm{Y}_{t}$ and $\bm{Y}_{t}^{\perp}$ are $\bm{Y}_{t+1} = \Lambda_{\Phi} \bm{Y}_{t} + \Tilde{\bm{\epsilon}}_{t+1}$ with $\Tilde{\bm{\epsilon}}_{t+1} \sim \mathcal{N}(0,U^{\prime} \Sigma U )$, and $\bm{Y}_{t+1}^{\perp} = \Lambda_{\Phi,\perp} \bm{Y}_{t}^{\perp} + \Tilde{\bm{\epsilon}}_{t+1}^{\perp}$ with $\Tilde{\bm{\epsilon}}_{t+1}^{\perp} \sim \mathcal{N}(0,U_{\perp}^{\prime} \Sigma U_{\perp} )$.
Define $W_{t} = \bm{Y}_{t} \bm{Y}_{t}^{\prime} - \Lambda_{\Gamma}$ and $W_{t}^{\perp} = \bm{Y}_{t}^{\perp} (\bm{Y}_{t}^{\perp})^{\prime} - \Lambda_{\Gamma,\perp}.$ Note that $\bm{Y}_{t} \bm{Y}_{t}^{\prime}  \sim W_{K}(1,\Lambda_{\Gamma})$, where $W_{K}(1,\Lambda_{\Gamma})$ is the $K$-dimensional Wishart distribution with $1$ degree of freedom and scale parameter $\Lambda_{\Gamma}$. Similarly, $\bm{Y}_{t}^{\perp} (\bm{Y}_{t}^{\perp})^{\prime} \sim W_{N-K}(1,\Lambda_{\Gamma,\perp})$.

Let $\bm{X}_{t}^{U} \coloneqq UU^{\prime} \bm{X}_{t}$ and $\bm{X}_{t}^{U_{\perp}} \coloneqq U_{\perp}U_{\perp}^{\prime} \bm{X}_{t}$ be  projections of $\bm{X}_{t}$ onto the subspaces corresponding to $U$ and $U_{\perp}$, respectively. Then $\bm{X}_{t}^{U} \sim \mathcal{N}(0,U\Lambda_{\Gamma} U^{\prime})$ and $\bm{X}_{t}^{U_{\perp}} \sim \mathcal{N}(0,U_{\perp}\Lambda_{\Gamma, \perp} U_{\perp}^{\prime})$.
Since $\text{Cov}(\bm{X}_{t}^{U},\bm{X}_{t}^{U_{\perp}}) = 0$, then $\bm{X}_{t}^{U}$ and $\bm{X}_{t}^{U_{\perp}}$ are independent. Further note that
\begin{align}
\label{eq:J-K-AR}
    \bm{X}_{t}^{U} &= \Phi \bm{X}_{t-1}^{U} + \bm{\epsilon}_{t}^{U} \quad \quad \bm{\epsilon}_{t}^{U} \sim \mathcal{N}(0,UU^{\prime}\Sigma UU^{\prime}) \nonumber \\
    \bm{X}_{t}^{U_{\perp}} &= \Phi \bm{X}_{t-1}^{U_{\perp}} + \bm{\epsilon}_{t}^{U_{\perp}} \quad \quad \bm{\epsilon}_{t}^{U_{\perp}} \sim \mathcal{N}(0,U_{\perp}U_{\perp}^{\prime}\Sigma U_{\perp}U_{\perp}^{\prime}).
\end{align}
Define 
    $G_{t} \coloneqq \bm{X}_{t}^{U}(\bm{X}_{t}^{U_{\perp}})^{\prime} + \bm{X}_{t}^{U_{\perp}}(\bm{X}_{t}^{U})^{\prime}$.
Then 
    $UU^{\prime}XX^{\prime}U_{ \perp}U_{ \perp}^{\prime} \nonumber + U_{ \perp}U_{ \perp}^{\prime}XX^{\prime}UU^{\prime} = \sum_{t=1}^{T} G_{t}$,
with $\mathbb{E}(G_{t}) = 0$.
We can write 
\begin{align}
\label{eq:S-Gamma}
     S_{T} - \Gamma  &= \frac{1}{T}\left( \sum_{t=1}^{T} \bm{X}_{t} \bm{X}_{t}^{\prime} - \Gamma \right) \nonumber \\
     &= \frac{1}{T}\left[ \sum_{t=1}^{T} \left\{ (UU^{\prime} + U_{\perp}U_{\perp}^{\prime})\bm{X}_{t} \bm{X}_{t}^{\prime}(UU^{\prime} + U_{\perp}U_{\perp}^{\prime}) - U \Lambda_{\Gamma} U^{\prime} - U_{\perp} \Lambda_{\Gamma , \perp} U_{\perp}^{\prime} \right\} \right] \nonumber \\  
     &= U \frac{1}{T} \sum_{t=1}^{T} W_{t} U^{\prime} + \frac{1}{T}  \sum_{t=1}^{T} G_{t} + U_{\perp} \frac{1}{T} \sum_{t=1}^{T} W_{t}^{\perp} U_{\perp}^{\prime}.
\end{align}
Using the triangle inequality, we can bound each of the terms of \eqref{eq:S-Gamma} in turn: 
\begin{align}
    \label{eq:triangle-S-Gamma}
    \lVert S_{T} - \Gamma\rVert_{2} \leq \frac{1}{T}\left\lVert \sum_{t=1}^{T} W_{t} \right\rVert_{2} + \frac{1}{T}\left\lVert \sum_{t=1}^{T} G_{t} \right\rVert_{2} + \frac{1}{T}\left\lVert \sum_{t=1}^{T} W_{t}^{\perp} \right\rVert_{2}.
\end{align}

We start with $\sum_{t=1}^{T} W_{t}$. 
For the truncated version of the matrix Bernstein inequality, we first choose an appropriate truncation level, $L_{1}$. We  show that $\mathbb{P}(\lVert  W_{t} \rVert_{2} \geq L_{1}) \leq q_{0}$ where 
    \begin{align}
    \label{eq:L_{1}-bound1} \nonumber
    L_{1} = C \left[  \sqrt{K} \left\{ \lambda_{1}(\Gamma) + \sqrt{\lambda_{1}(\Gamma) \log\left( \frac{1}{q_{0}} \right)} \right\} + \log\left( \frac{1}{q_{0}} \right) + \sqrt{\lambda_{1}(\Gamma)\log\left( \frac{1}{q_{0}} \right)} + \lambda_{1}(\Gamma)  \right].
\end{align}
This result follows from  Theorem 3.8 of \citet[][]{cai2022non} which gives
\begin{align}
\mathbb{P}\left( \lVert  W_{t} \rVert_{2} \gtrsim  \lambda_{1}(\Gamma) + x^{2} + 2\sqrt{\lambda_{1}(\Gamma) \sum_{i=1}^{K} \lambda_{i}(\Gamma)} + 2x \sqrt{\sum_{i=1}^{K}\lambda_{i}(\Gamma)} +2x\sqrt{\lambda_{1}(\Gamma)} \right) \leq e^{-x^{2}}.
\end{align}
Using $\sum_{i=1}^{K} \lambda_{i}(\Gamma) \leq K \lambda_{1}(\Gamma)$ and setting $q_{0} = \exp (-x^{2})$ then gives the desired high probability bound. Note that if we choose $q_{0} = T^{-\nu }$, then $L_{1} = O(\log T)$. 

Next, we show that $\lVert\mathbb{E}(W_{t}) - \mathbb{E}(W_{t} \indfun\{\lVert W_{t} \rVert_{2} < L_{1}\}) \rVert_{2} \leq q_{1}$ 
with 
$q_{1} = \lambda_{1}(\Gamma)\sqrt{(K^{2} +4K+1)q_{0}}$. Note that 
\begin{align}
    \mathbb{E}(W_{t}) = \mathbb{E}(W_{t} \indfun\{\lVert W_{t} \rVert_{2} < L_{1}\}) + \mathbb{E}(W_{t} \indfun\{\lVert W_{t} \rVert_{2} > L_{1}\}),
\end{align}
and since $\mathbb{E}(W_{t}) = 0$, then $\lVert\mathbb{E}(W_{t} \indfun\{\lVert W_{t} \rVert_{2} < L_{1}\}) \rVert_{2} = \lVert\mathbb{E}(W_{t} \indfun\{\lVert W_{t} \rVert_{2} > L_{1}\}) \rVert_{2}$.
By Jensen's inequality,
    \begin{align}
        \lVert&\mathbb{E}(W_{t}) - \mathbb{E}(W_{t} \indfun\{\lVert W_{t} \rVert_{2} < L_{1}\}) \rVert_{2} = \lVert\mathbb{E}(W_{t} \indfun\{\lVert W_{t} \rVert_{2} < L_{1}\}) \rVert_{2} 
        \nonumber \\ 
        &= \lVert\mathbb{E}(W_{t} \indfun\{\lVert W_{t} \rVert_{2} > L_{1}\}) \rVert_{2} \leq \mathbb{E}(\left\lVert W_{t} \indfun\{\lVert W_{t} \rVert_{2} > L_{1}\} \right\rVert_{2}) 
        = \mathbb{E}(\left\lVert W_{t}\right\rVert_{2} \indfun\{\lVert W_{t} \rVert_{2} > L_{1}\} ) \nonumber \\ 
        &\leq \sqrt{\mathbb{E}(\left\lVert W_{t}\right\rVert_{2}^{2}) \mathbb{E}(\indfun\{\lVert W_{t} \rVert_{2} > L_{1}\}) } 
        = \sqrt{\mathbb{E}(\left\lVert W_{t}\right\rVert_{2}^{2}) \mathbb{P}(\lVert W_{t} \rVert_{2} > L_{1}) },
    \end{align}
    where the second inequality follows from the Cauchy-Schwarz inequality. Since $W_{t} = \bm{Y}_{t}\bm{Y}_{t}^{\prime} - \Lambda_{\Gamma}$, we have that $\lVert W_{t}\rVert_{2} \leq \lVert \bm{Y}_{t}\bm{Y}_{t}^{\prime}\rVert_{2} + \lambda_{1}(\Gamma)$ and thus
    \begin{align}
        \mathbb{E}(\lVert W_{t}\rVert_{2}^{2}) \leq \mathbb{E}(\lVert \bm{Y}_{t}\bm{Y}_{t}^{\prime}\rVert_{2}^{2}) + 2 \lambda_{1}(\Gamma)\mathbb{E}(\lVert\bm{Y}_{t}\bm{Y}_{t}^{\prime}\rVert_{2})    + \lambda_{1}^{2}(\Gamma).
    \end{align}
    Now $\bm{Y}_{t}\bm{Y}_{t}^{\prime}$ is a rank one matrix with eigenvector $\bm{Y}_{t}$ and eigenvalue $\bm{Y}_{t}^{\prime}\bm{Y}_{t} = \sum_{i = 1}^{K} (\bm{Y}_{t})_{i}^{2}$. But $(\bm{Y}_{t})_{i}^{2} \sim \lambda_{i}(\Gamma) \chi_{1}^{2}$, and therefore
\begin{align}
\label{eq:chi-squared-Z}
    \lVert \bm{Y}_{t}\bm{Y}_{t}^{\prime}\rVert_{2} = \sum_{i = 1}^{K} (\bm{Y}_{t})_{i}^{2} 
    = \sum_{i = 1}^{K} \lambda_{i}(\Gamma) Z_{i}^{2}
    \leq \lambda_{1}(\Gamma) \sum_{i = 1}^{K}  Z_{i}^{2},
\end{align}
where $Z_{i} \overset{iid}{\sim} \mathcal{N}(0,1)$. Therefore, $\lVert \bm{Y}_{t}\bm{Y}_{t}^{\prime}\rVert_{2} \leq \lambda_{1}(\Gamma) \Upsilon$ where $\Upsilon \sim \chi_{K}^{2}$. Using $\mathbb{E}(\Upsilon) = K$ and $\mathrm{Var}(\Upsilon) = 2K$, 
    \begin{align}
        \mathbb{E}(\lVert W_{t}\rVert_{2}^{2}) &\leq \lambda_{1}^{2}(\Gamma) \left\{\mathbb{E}(\Upsilon^{2}) + 2\mathbb{E}(\Upsilon)    + 1\right\} 
        = \lambda_{1}^{2}(\Gamma) \left\{K^{2} +4K+1\right\}.
    \end{align}
    Since $\mathbb{P}(\lVert  W_{t} \rVert_{2} \geq L_{1}) \leq q_{0}$, 
    \begin{align}
        \lVert\mathbb{E}(W_{t}) - \mathbb{E}(W_{t} \indfun\{\lVert W_{t} \rVert_{2} < L_{1}\}) \rVert_{2} \leq \lambda_{1}(\Gamma)\sqrt{(K^{2} +4K+1)q_{0}}. 
    \end{align}
Note that if $q_{0} = T^{-\nu}$, then $q_{1} = O\left(T^{-\nu/2}\right)$.

Next, we compute 
\begin{align}
\label{eq:variance-param}
        v^{2} = \sup_{\vartheta \subseteq \{1,\dots,T\}} \frac{1}{\mathrm{Card}(\vartheta)} \lambda_{\text{max}}\left[\mathbb{E} \left( \sum_{t \in \vartheta} W_{t}  \right)^{2} \right] = \frac{1}{T}\lambda_{\text{max}}\left[\mathbb{E} \left( \sum_{t = 1}^{T} W_{t}  \right)^{2} \right].
\end{align}
We start by writing 
\begin{align}
\label{eq:variance-expand}
    \mathbb{E} \left( \sum_{t = 1}^{T}  W_{t}  \right)^{2} = \sum_{t=1}^{T} \mathbb{E}(W_{t}^{2} ) + \sum_{t=1}^{T-1} \sum_{h=1}^{T-t} \mathbb{E}(W_{t} W_{t+h}) +  \mathbb{E}(W_{t+h} W_{t}).
\end{align}
First, $\mathbb{E}(W_{t}^{2} ) = \mathbb{E}\{ (\bm{Y}_{t} \bm{Y}_{t}^{\prime} - \Lambda_{\Gamma} )^{2}\} = \mathbb{E} (\bm{Y}_{t} \bm{Y}_{t}^{\prime} )^{2} - \Lambda_{\Gamma}^{2}$. The second moments of $\bm{Y}_{t} \bm{Y}_{t}^{\prime} \sim W_{K}(1,\Lambda_{\Gamma})$ are 
$
    \mathbb{E} (\bm{Y}_{t} \bm{Y}_{t}^{\prime} )^{2} = \Lambda_{\Gamma} \mathrm{tr}(\Lambda_{\Gamma}) + 2 \Lambda_{\Gamma}^{2}
$ \citep[see e.g.][]{letac2004all}.
Hence: 
\begin{equation}
\label{eq:diag-term}
    \sum_{t=1}^{T} \mathbb{E}(W_{t}^{2} ) = T \{\Lambda_{\Gamma} \mathrm{tr} (\Lambda_{\Gamma}) + \Lambda_{\Gamma}^{2}\}. 
\end{equation}
Next, we look at the cross terms in \eqref{eq:variance-expand}. 
Now 
\begin{align}
    \mathbb{E}(W_{t} W_{t+h}) = \mathbb{E}\{(\bm{Y}_{t} \bm{Y}_{t}^{\prime} - \Lambda_{\Gamma}) (\bm{Y}_{t+h} \bm{Y}_{t+h}^{\prime} - \Lambda_{\Gamma}) \} 
    = \mathbb{E}(\bm{Y}_{t} \bm{Y}_{t}^{\prime} \bm{Y}_{t+h} \bm{Y}_{t+h}^{\prime} ) - \Lambda_{\Gamma}^{2}.
\end{align}
We will compute the $(ij)$-th element of $\mathbb{E}(\bm{Y}_{t} \bm{Y}_{t}^{\prime} \bm{Y}_{t+h} \bm{Y}_{t+h}^{\prime} )$ using Isserlis' theorem \citep[see, for example][]{vignat2012generalized}. Noting that $\bm{Y}_{t}^{\prime} \bm{Y}_{t+h}$ is a scalar, we have
\begin{align}
    \label{eq:isserlis}
    &\mathbb{E}\{ (\bm{Y}_{t}^{\prime} \bm{Y}_{t+h}) (\bm{Y}_{t})_{i} (\bm{Y}_{t+h})_{j} \} = \sum_{k=1}^{K} \mathbb{E}\{ (\bm{Y}_{t})_{k} (\bm{Y}_{t+h})_{k} (\bm{Y}_{t})_{i} (\bm{Y}_{t+h})_{j} \} \\ 
    & = \sum_{k=1}^{K} \Big[\mathbb{E}\{ (\bm{Y}_{t})_{k} (\bm{Y}_{t+h})_{k} \} \mathbb{E}\{ (\bm{Y}_{t})_{i} (\bm{Y}_{t+h})_{j} \} + \mathbb{E}\{ (\bm{Y}_{t})_{k} (\bm{Y}_{t})_{i} \} \mathbb{E}\{ (\bm{Y}_{t+h})_{k} (\bm{Y}_{t+h})_{j} \} \nonumber \\ &\quad \quad \quad + \mathbb{E}\{ (\bm{Y}_{t})_{k} (\bm{Y}_{t+h})_{j} \} \mathbb{E}\{ (\bm{Y}_{t+h})_{k} (\bm{Y}_{t})_{i} \}\Big], 
\end{align}
where the last line used Isserlis' theorem. Now, 
\begin{align}
    \mathbb{E}\{ (\bm{Y}_{t})_{i} (\bm{Y}_{t+h})_{j} \} &= \mathbb{E}\left\{ (\bm{Y}_{t})_{i} \left(\Lambda_{\Phi}^{h}\bm{Y}_{t} + \sum_{l=0}^{h-1} \Lambda_{\Phi}^{l} \Tilde{\bm{\epsilon}}_{t+h-l}\right)_{j} \right\} 
    = \lambda_{j}^{h}(\Phi) \lambda_{i}(\Gamma) \delta_{ij},
\end{align}
where $\delta_{ij}$ is the Kronecker delta. Thus, 
\begin{align}
    \mathbb{E}\{ (\bm{Y}_{t}^{\prime} \bm{Y}_{t+h}) (\bm{Y}_{t})_{i} (\bm{Y}_{t+h})_{j} \} = &\, \sum_{k=1}^{K} \lambda_{k}(\Phi)^{h} \lambda_{k}(\Gamma) \lambda_{j}^{h}(\Phi) \lambda_{i}(\Gamma) \delta_{ij} + \lambda_{i}(\Gamma) \delta_{ik} \lambda_{j}(\Gamma) \delta_{kj} \nonumber \\ &+ \lambda_{j}^{h}(\Phi) \lambda_{k}(\Gamma) \delta_{kj} \lambda_{k}(\Phi)^{h} \lambda_{i}(\Gamma) \delta_{ik} \nonumber \\ 
    =&\, \lambda_{j}^{h}(\Phi) \lambda_{i}(\Gamma) \delta_{ij} \mathrm{tr}(\Lambda_{\Phi}^{h} \Lambda_{\Gamma}) + \lambda_{i}^{2}(\Gamma) \delta_{ij} + \{\lambda_{i}^{h}(\Phi) \lambda_{i}(\Gamma) \}^{2} \delta_{ij}.
\end{align}
Therefore $\mathbb{E}(\bm{Y}_{t} \bm{Y}_{t}^{\prime} \bm{Y}_{t+h} \bm{Y}_{t+h}^{\prime} ) = (\Lambda_{\Phi}^{h} \Lambda_{\Gamma}) \mathrm{tr}(\Lambda_{\Phi}^{h} \Lambda_{\Gamma}) + \Lambda_{\Gamma}^{2} + (\Lambda_{\Phi}^{h} \Lambda_{\Gamma})^{2}$ and, finally, 
\begin{align}
\label{eq:cross-terms}
    \mathbb{E}(W_{t} W_{t+h}) = (\Lambda_{\Phi}^{h} \Lambda_{\Gamma}) \mathrm{tr}(\Lambda_{\Phi}^{h} \Lambda_{\Gamma}) +  (\Lambda_{\Phi}^{h} \Lambda_{\Gamma})^{2}.
\end{align}
Similarly, $\mathbb{E}(W_{t+h} W_{t})$ is given by the right hand side of \eqref{eq:cross-terms}.

Noting that 
$\mathrm{tr}(\Lambda_{\Gamma}) \leq K/(1 - \rho^{2})$
and  
$\mathrm{tr}(\Lambda_{\Phi}^{h} \Lambda_{\Gamma}) \leq K \rho^{h}/(1 - \rho^{2})$,
we have that, on combining \eqref{eq:diag-term} and \eqref{eq:cross-terms},
\begin{align}
    \lambda_{\text{max}}\left[\mathbb{E} \left( \sum_{t =1}^{T} W_{t}  \right)^{2} \right]  \leq T\left( \frac{K+1}{(1 - \rho^{2})^{2}}\right) + 2 \frac{K+1}{(1 - \rho^{2})^{2}} \sum_{h=1}^{T-1} (T-h) \rho^{2h} \lesssim \frac{KT}{(1-\rho^{2})^{3}}.
\end{align}
Therefore, $v^{2} \lesssim K/(1-\rho^{2})^{3}$.
We now show that $\mathbb{E}(\lVert W_{t} \rVert_{2}^{4}) = O(1)$, which will allow us to use Lemma \ref{lemma:truncated-Bernstein}. Since $\lVert W_{t} \rVert_{2} \leq \lVert \bm{Y}_{t}\bm{Y}_{t}^{\prime} \rVert_{2} + \lambda_{1}(\Gamma)$, then $\mathbb{E}\lVert W_{t} \rVert_{2}^{4} \lesssim \mathbb{E}\lVert \bm{Y}_{t} \bm{Y}_{t}^{\prime} \rVert_{2}^{4} = \mathbb{E}\lVert \bm{Y}_{t} \rVert_{2}^{8} \lesssim K^{4}$. Since $K$ is fixed, this implies that  $\mathbb{E}(\lVert W_{t} \rVert_{2}^{4}) = O(1)$. Therefore, setting $q_{0} \lesssim T^{-\nu}$ for any $\nu \geq 4$ allows us to use the truncated matrix Bernstein inequality in Lemma \ref{lemma:truncated-Bernstein}. Indeed, by Lemma \ref{lemma:truncated-Bernstein},
\begin{align}
        \mathbb{P}\left( \left\lVert \sum_{t = 1}^{T} W_{t} \right\rVert_{2} \geq x \right) \leq Tq_{0} + K\exp\left(-\frac{C(x-Tq_{1})^{2}}{v^{2}T + c_{\beta}^{-1}L_{1}^{2} + (x-Tq_{1})L_{1}\gamma(c_{\beta},T)}\right).
    \end{align}
The requirement on $T$ given by \eqref{eq:lemma_T_condition} then implies 
\begin{align}
    \frac{1}{T} \left\lVert \sum_{t = 1}^{T} W_{t} \right\rVert_{2} \lesssim \frac{1}{\{1 - \rho^{2}\}^{3/2}} \sqrt{\frac{K}{T}} \log^{1/2} T
\end{align}
with probability exceeding $1 - T^{-\nu}$.   

Next, we bound $\sum_{t=1}^{T} W_{t}^{\perp}$. 
Following the same proof as for bounding $\sum_{t=1}^{T} W_{t}$, but with $K$ replace with $N-K$ and $\lambda_{1}(\Phi)$ replaced with $\lambda_{K+1}(\Phi)$, the variance parameter satisfies $v_{\perp} \lesssim N/\{1-\lambda_{K+1}^{2}(\Phi)\}^{3}$. By Lemma~\ref{lemma:eigs}, $\lambda_{K+1}(\Phi) = o(1)$. Therefore, if 
        $T \geq C \gamma(c_{\beta},T)^{2} \log^{3}(T+N)$
for some sufficiently large constant $C >0$, then
\begin{align}\label{eq:Q_t_bound}
    \frac{1}{T} \left\lVert \sum_{t = 1}^{T} W_{t}^{\perp} \right\rVert_{2} \lesssim  \sqrt{\frac{N}{T}} \log^{1/2}(T+N)
\end{align}
with probability exceeding $1 - (T+N)^{-\nu}$ for any $\nu \geq 8$. Note that we choose $\nu \geq 8$ here since $\mathbb{E}\lVert W_{t}^{\perp}\rVert_{2}^{4} \lesssim N^{4}$ in this case.

Next, we bound the cross term, $\sum_{t=1}^{T} G_{t}$. 
We start by choosing the truncation level as the high probability bound on $\lVert G_{t} \rVert_{2}$. We claim that $\mathbb{P}(\lVert  G_{t} \rVert_{2} \geq L_{G}) \leq T^{-\nu}$ where $$L_{G} = C\sqrt{K(N-K)\lambda_{1}(\Gamma)\lambda_{K+1}(\Gamma)}\log T.$$
   We first note that $\lVert \bm{X}_{t}^{U}\rVert_{2} = \lVert U \bm{Y}_{t} \rVert_{2} = \lVert \bm{Y}_{t}\rVert_{2} $ and $\lVert \bm{X}_{t}^{U_{\perp}}\rVert_{2} = \lVert U_{\perp} \bm{Y}_{t}^{\perp} \rVert_{2} = \lVert \bm{Y}_{t}^{\perp}\rVert_{2} $ where $\bm{Y}_{t} \sim \mathcal{N}(0,\Lambda_{\Gamma})$ and $\bm{Y}_{t}^{\perp} \sim \mathcal{N}(0,\Lambda_{\Gamma, \perp})$. The Chernoff bound for Gaussian random variables gives 
    \begin{align}
        \lVert \bm{Y}_{t}\lVert _{\infty } \leq C\sqrt{\lambda_{1}(\Gamma)\log T}, & & 
        \lVert \bm{Y}_{t}^{\perp}\lVert _{\infty } \leq C\sqrt{ \lambda_{K+1}(\Gamma) \log T},
    \end{align}
    each with probability greater or equal to $1 - T^{-\nu /2}$. Therefore 
    \begin{align}
        \lVert G_{t}\rVert_{2} &= \lVert \bm{X}_{t}^{U}(\bm{X}_{t}^{U_{\perp}})^{\prime} + \bm{X}_{t}^{U_{\perp}}(\bm{X}_{t}^{U})^{\prime} \rVert_{2} 
        \leq 2\lVert \bm{X}_{t}^{U}\rVert_{2} \lVert \bm{X}_{t}^{U_{\perp}}\rVert_{2} \nonumber \\ 
        &\leq 2 \sqrt{K} \lVert \bm{Y}_{t}\lVert _{\infty } \sqrt{N-K} \lVert \bm{Y}_{t}^{\perp}\lVert _{\infty }
        \leq C\sqrt{K(N-K)\lambda_{1}(\Gamma)\lambda_{K+1}(\Gamma)}\log T 
    \end{align}
    with probability greater or equal to $1 - T^{-\nu}$.
    
Since $\bm{X}_{t}^{U}$ and $\bm{X}_{t}^{U_{\perp}}$ are independent Gaussian random vectors with $\mathbb{E}(\bm{X}_{t}^{U}) = 0$ and $\mathbb{E}(\bm{X}_{t}^{U_{\perp}}) =0$, then the distribution of $\bm{X}_{t}^{U}(\bm{X}_{t}^{U_{\perp}})^{\prime}$ is symmetric around 0. Therefore the distribution of the truncated random variable $G_{t} \indfun\{\lVert G_{t} \rVert_{2} < L\}$ is also symmetric around 0, implying that $\lVert  \mathbb{E}(G_{t} \indfun\{\lVert G_{t} \rVert_{2} < L_{G}\}) \rVert_{2} = 0$ for any $L_{G} > 0$.

Next, we compute 
$v_{G}^{2} = \lambda_{\text{max}}[\mathbb{E} ( \sum_{t =1}^{T} G_{t}  )^{2} ]/T$.
Write 
\begin{align}
\label{eq:variance-expand-G}
    \mathbb{E} \left( \sum_{t = 1}^{T}  G_{t}  \right)^{2} = \sum_{t=1}^{T} \mathbb{E}(G_{t}^{2} ) + 2\sum_{t=1}^{T-1} \sum_{h=1}^{T-t} \mathbb{E}(G_{t} G_{t+h}).
\end{align}
Repeated application of \eqref{eq:J-K-AR} gives 
\begin{align}
    \bm{X}_{t+h}^{U} = \Phi^{h} \bm{X}_{t}^{U} + \sum_{s=0}^{h-1} \Phi^{s} \bm{\epsilon}_{t+h-j}^{U} & & 
    \bm{X}_{t+h}^{U_{\perp}} = \Phi^{h} \bm{X}_{t}^{U_{\perp}} + \sum_{s=0}^{h-1} \Phi^{s} \bm{\epsilon}_{t+h-j}^{U_{\perp}},
\end{align}
from which it is seen that $\mathbb{E}((\bm{X}_{t}^{U_{\perp}})^{\prime}\bm{X}_{t+h}^{U}) = 0$ and $\mathbb{E}((\bm{X}_{t}^{U})^{\prime}\bm{X}_{t+h}^{U_{\perp}}) = 0$. Thus,
\begin{align}
    \mathbb{E}(G_{t}G_{t+h}) &= \mathbb{E}\{ (\bm{X}_{t}^{U} (\bm{X}_{t}^{U_{\perp}})^{\prime} + \bm{X}_{t}^{U_{\perp}} (\bm{X}_{t}^{U})^{\prime})(\bm{X}_{t+h}^{U} (\bm{X}_{t+h}^{U_{\perp}})^{\prime} + \bm{X}_{t+h}^{U_{\perp}} (\bm{X}_{t+h}^{U})^{\prime})\} \nonumber \\ 
    &= \mathbb{E}( (\bm{X}_{t}^{U_{\perp}})^{\prime}\bm{X}_{t+h}^{U_{\perp}}) \mathbb{E}(\bm{X}_{t}^{U} (\bm{X}_{t+h}^{U})^{\prime}) + \mathbb{E}( (\bm{X}_{t}^{U})^{\prime}\bm{X}_{t+h}^{U}) \mathbb{E}(\bm{X}_{t}^{U_{\perp}} (\bm{X}_{t+h}^{U_{\perp}})^{\prime}).
\end{align}
Now, 
\begin{align}
    \mathbb{E}( (\bm{X}_{t}^{U})^{\prime}\bm{X}_{t+h}^{U}) &= \mathbb{E}\left\{(\bm{X}_{t}^{U})^{\prime} (\Phi^{h} \bm{X}_{t}^{U} + \sum_{s=0}^{h-1} \Phi^{s} \bm{\epsilon}_{t+h-j}^{U}) \right\} 
    = \mathbb{E}((\bm{X}_{t}^{U})^{\prime} \Phi^{h} \bm{X}_{t}^{U}) \\ &= \mathrm{tr}(\Phi^{h} U \Lambda_{\Gamma} U^{\prime}) 
    = \mathrm{tr}(U \Lambda_{\Phi}^{h} \Lambda_{\Gamma} U^{\prime}) 
    = \mathrm{tr}( \Lambda_{\Phi}^{h} \Lambda_{\Gamma}),
\end{align}
where we used $\Phi^{h}U = U \Lambda_{\Phi}^{h}$ and $\mathbb{E}((\bm{X}_{t}^{U})^{\prime} \bm{\epsilon}_{s}^{U}) =0$ for $t \neq s$. Similarly, 
    $\mathbb{E}( (\bm{X}_{t}^{U_{\perp}})^{\prime}\bm{X}_{t+h}^{U_{\perp}}) = \mathrm{tr}( \Lambda_{\Phi,\perp}^{h} \Lambda_{\Gamma,\perp})$,  
    $\mathbb{E}( \bm{X}_{t}^{U}(\bm{X}_{t+h}^{U})^{\prime}) = U \Lambda_{\Gamma} \Lambda_{\Phi}^{h} U^{\prime}$, and
    $\mathbb{E}( \bm{X}_{t}^{U_{\perp}}(\bm{X}_{t+h}^{U_{\perp}})^{\prime}) = U_{\perp} \Lambda_{\Gamma,\perp} \Lambda_{\Phi,\perp}^{h} U_{\perp}^{\prime}$.
Therefore, 
\begin{align}
\label{eq:G-variance-cross}
    \mathbb{E}(G_{t}G_{t+h}) &= U \Lambda_{\Gamma} \Lambda_{\Phi}^{h} U^{\prime}\mathrm{tr}( \Lambda_{\Phi,\perp}^{h} \Lambda_{\Gamma,\perp}) + U_{\perp} \Lambda_{\Gamma,\perp} \Lambda_{\Phi,\perp}^{h} U_{\perp}^{\prime} \mathrm{tr}( \Lambda_{\Phi}^{h} \Lambda_{\Gamma}).
\end{align}
Plugging \eqref{eq:G-variance-cross} into \eqref{eq:variance-expand-G} yields
\begin{align}
    \lambda_{\text{max}}\left[\mathbb{E} \left( \sum_{t =1}^{T} G_{t}  \right)^{2} \right] 
    &\leq \frac{2N}{\{1-\lambda_{1}^{2}(\Phi)\}\{1-\lambda_{K+1}^{2}(\Phi)\}} \sum_{h=0}^{T-1} (T-h) \{\lambda_{1}(\Phi)\lambda_{K+1}(\Phi)\}^{h} \nonumber \\ 
    &\leq T\frac{2N}{\{1-\lambda_{1}^{2}(\Phi)\}\{1-\lambda_{K+1}^{2}(\Phi)\}\{1 - \lambda_{1}(\Phi)\lambda_{K+1}(\Phi)\}}.
\end{align}
Recalling that $\lambda_{1}(\Phi) = 1/(1-\rho^{2})$ and $\lambda_{K+1}(\Lambda_{\Phi}) = o(1)$ gives $v_{G}^{2} \lesssim N/{(1 - \rho^{2})}$.

Next we analyse $\mathbb{E}\lVert G_{t}\rVert_{2}^{4}$. Notice that $\lVert G_{t} \rVert_{2}$ is the product of two subgaussian random variables. This follows from 
$\lambda_{1}(G_{t}^{2}) = \lVert \bm{X}_{t}^{U} \rVert_{2}^{2} \lVert \bm{X}_{t}^{U_{\perp}} \rVert_{2}^{2}$ and thus $\lVert G_{t} \rVert_{2} = \sqrt{\lambda_{1}(G_{t}^{2})} = \lVert \bm{X}_{t}^{U} \rVert_{2} \lVert \bm{X}_{t}^{U_{\perp}} \rVert_{2}$. Since $\lVert \bm{X}_{t}^{U} \rVert_{2} $ and $\lVert \bm{X}_{t}^{U_{\perp}} \rVert_{2}$ are independent, we have
\begin{align}
    \mathbb{E}\left(\lVert G_{t} \rVert_{2}^{4} \right) = \mathbb{E}\left(\lVert \bm{X}_{t}^{U} \rVert_{2}^{4} \right) \mathbb{E}\left(\lVert \bm{X}_{t}^{U_{\perp}} \rVert_{2}^{4} \right) \lesssim N^{2},
\end{align}
where we recall from earlier that $\mathbb{E}\left(\lVert \bm{X}_{t}^{U} \rVert_{2}^{4} \right) \lesssim \lambda_{1}^{2}(\Lambda_{\Gamma})K^{2} = O(1)$ and $\mathbb{E}(\lVert \bm{X}_{t}^{U_{\perp}} \rVert_{2}^{4} ) \lesssim \lambda_{K+1}^{2}(\Lambda_{\Gamma})N^{2}$. Thus, for any $q_{0} \lesssim (T+N)^{-\nu}$ with $\nu \geq 6$, by Lemma \ref{lemma:truncated-Bernstein},
\begin{align}
    \label{eq:bernstein-G}
        \mathbb{P}\left( \left\lVert \sum_{t = 1}^{T} G_{t} \right\rVert_{2} \geq x \right) \leq Tq_{0} + N\exp\left(-\frac{Cx^{2}}{v_{G}^{2}T + c_{\beta}^{-1}L_{G}^{2} + xL_{G}\gamma(c_{\beta},T)}\right),
    \end{align}
The condition 
$
    T > C \gamma(c_{\beta},T)^{2} \log^{3} (T+N)
$
for some sufficiently large $C > 0$ then implies that
\begin{align}\label{eq:G_t_bound}
   \frac{1}{T}\left\lVert \sum_{t = 1}^{T} G_{t} \right\rVert_{2} \lesssim \frac{1}{\sqrt{1 - \rho^{2}}  }\sqrt{\frac{N}{T}} \log^{1/2}(T+N)
\end{align}
with probability exceeding $1 - (T+N)^{-\nu}$. 

Combining the bounds for each of the three terms in \eqref{eq:triangle-S-Gamma} yields the result. \qed


\end{appendix}

\begin{acks}[Acknowledgments]
We thank Marcelo Medeiros for insightful discussions. 

Brendan Martin acknowledges funding from the Engineering and Physical Sciences Research
Council (EPSRC), grant number EP/S023151/1.
Francesco Sanna Passino acknowledges funding from the EPSRC, grant number EP/Y002113/1.

Brendan Martin is also affiliated with the Department of Mathematics, Imperial College London (United Kingdom). Mihai Cucuringu is also affiliated with the Department of Statistics and the Oxford-Man Institute of Quantitative Finance at the University of Oxford (United Kingdom), and with The Alan Turing Institute (United Kingdom). Alessandra Luati is also affiliated with the Department of Statistical Sciences at the University of Bologna (Italy). 
\end{acks}


\begin{supplement}
Proofs of 
theorems and technical lemmas are reported in the supplementary material. 
Code for the numerical experiments 
is available at \url{https://github.com/bmartin9/nirvar_theory}.
\end{supplement}




\bibliographystyle{rss}
\singlespacing
\bibliography{references.bib}


\appendix
\newpage


\setcounter{section}{0}
\renewcommand{\thesection}{S\arabic{section}}

\setcounter{equation}{0}
\renewcommand{\theequation}{S.\arabic{equation}}

\setcounter{figure}{0}
\renewcommand{\thefigure}{S.\arabic{figure}}

\setcounter{table}{0}
\renewcommand{\thetable}{S.\arabic{table}}

\begin{center}
{\LARGE\textbf{SUPPLEMENTARY MATERIAL}}
\end{center}

\section{Adjacency model}\label{sec:appendix-adjacency}

\subsection{Proof of  Theorem~\ref{thm:dist_U_S_U_P}}  \label{proof:dist_U_S_U_P}
    The Davis-Kahan sin$\Theta$ theorem \citepSM{davis1970rotation} implies that if 
    \begin{align}\label{eq:davis-kahan-assumption}
        \left\lVert S_{T} - \Gamma \right\rVert_{2} < (1 - 1/\sqrt{2})\Delta_{\Gamma},
    \end{align}
    then 
    \begin{equation}\label{eq:D-K}
        \text{dist} (U_{S}, U) \leq \frac{2\left\lVert S_{T} - \Gamma \right\rVert_{2} }{\Delta_{\Gamma}}. 
    \end{equation}
    Recall $\Delta_{\Gamma} = \lambda_{K}(\Gamma) - \lambda_{K+1}(\Gamma)$. The assumption in~\eqref{eq:T_cond_combined} on $T$ implies that 
    \begin{align}\label{eq:num_obs_rate}
        \frac{2\left\lVert S_{T} - \Gamma \right\rVert_{2} }{\Delta_{\Gamma}} &< \frac{1}{\log T} 
    \end{align}
    and therefore \eqref{eq:davis-kahan-assumption} holds under the assumptions of the theorem. By Lemma \ref{lemma:S_T-Gamma}, 
    \begin{align}
        \left\lVert S_{T} - \Gamma \right\rVert_{2} &\lesssim \Bigg{[} \frac{1}{\{1 - \rho^{2}\}^{3/2}} \sqrt{\frac{K}{T}} +  \left\{1 + \frac{1}{ \sqrt{1 - \rho^{2}}}\right\}\sqrt{\frac{N}{T}} \Bigg{]} \log^{\frac{1}{2}}(T+N)
    \end{align}
    with probability exceeding $1 - (T+N)^{-\nu}$ as long as \eqref{eq:num_obs_rate} and \eqref{eq:lemma_T_condition} hold. 
    These conditions hold by the assumption given by \eqref{eq:T_cond_combined}.  By Assumption \ref{ass:regularity}, Lemma \ref{lemma:eigs}, and \citetSM[][Proposition 4.1]{martin2024nirvar_arxiv}, for any $s>0$, with probability at least $1-N^{-s}$,
     $\lambda_{r}(\Gamma) \asymp 1/\{1 - \rho^{2}\lambda_{r}^{2}(B)/\lambda_{1}^{2}(B)\}$ for $r = 1,\dots, K$ and $\lambda_{r}(\Gamma)  \asymp 1$ for $r=K+1,\dots,N$. Note that $\lambda_{1}(\Phi) = \rho$, and thus $\lambda_{1}(\Gamma) = 1/(1 - \rho^{2})$. We also have $\lambda_{K}(\Gamma) \asymp \kappa^{2}/(\kappa^{2} - \rho^{2})$ where we define the condition number of $B$ as $\kappa \coloneqq \lambda_{1}(B)/\lambda_{K}(B)$. 
Therefore, with probability at least $1-N^{-s}$,
\begin{align}
    \frac{1}{\Delta_{\Gamma}} \asymp \frac{\kappa^{2} - \rho^{2} }{\rho^{2} } \leq \frac{\kappa^{2}}{\rho^{2}} \asymp \frac{1}{\rho^{2}\lambda_{K}^{2}(B)},
\end{align}
which, from \eqref{eq:D-K}, implies that, with probability exceeding $1 - (T+N)^{-\nu}$,
\begin{align}\label{eq:dist_U_S}
        \text{dist}(U_{S},U) &\lesssim \left\{ \frac{1}{\rho^{2} \lambda_{K}^{2}(B) }\right\}\left\{ \frac{1}{(1 - \rho^{2})^{3/2}} \sqrt{\frac{K}{T}} +  \left(1 + \frac{1}{ \sqrt{1 - \rho^{2}}}\right)\sqrt{\frac{N}{T}}\right\} \log^{\frac{1}{2}}(T+N).
    \end{align}
By \citetSM[][Lemma 5.1]{lei2015consistency}, there exists a $K \times K$ orthogonal matrix $R$ such that 
\begin{align}
    \lVert U - U_{P}R \rVert_{F} \leq \frac{2\sqrt{2K}}{\lambda_{K}(P)} \lVert A - P \rVert_{2}.
\end{align}
Under the assumptions of the theorem, by \citetSM[][Theorem 5.2]{lei2015consistency}, $\lVert A - P \rVert_{2} \lesssim\sqrt{N \alpha_{N}}$ with probability at least $1 - N^{-\nu}$. Therefore, using $\lVert U - U_{P}R \rVert_{2} \leq \lVert U - U_{P}R \rVert_{F}$, and $\lambda_{K}(P) \asymp N\alpha_{N}\lambda_{K}(B)/K$, we have 
\begin{align}\label{eq:dist_U_UP}
    \text{dist}(U,U_{P}) \lesssim \frac{K\sqrt{K}}{\sqrt{N \alpha_{N}}\lambda_{K}(B)}
\end{align}
with probability at least $1 - N^{-\nu}$. By the triangle inequality,
\begin{align}\label{eq:triangle}
    \text{dist}(U_{S},U_{P}) \leq \text{dist}(U_{S},U) + \text{dist}(U,U_{P}).
\end{align}
The result follows from substituting \eqref{eq:dist_U_S} and \eqref{eq:dist_U_UP} into \eqref{eq:triangle}. \qed

\subsection{Proof of Theorem~\ref{thm:kmeans}} \label{proof:kmeans}

    We apply \citetSM[][Lemma 5.3]{lei2015consistency} to $U_{S}$ and $U_{P}R^{*}$ where $R^{*} = \argmin_{R \in \mathbb{R}^{K \times K}} \lVert U_{S} - U_{P}R \rVert_{2}$. By \citepSM[][Lemma 2.1]{lei2015consistency}, $U_{P}R = Z\Theta R = Z \Tilde{\Theta}$ where $\lVert \Tilde{\Theta}_{k,:} - \Tilde{\Theta}_{l,:} \rVert_{2} = \sqrt{1/N_{k} + 1/N_{l}}$. We  choose $\delta_{k} \leq \min_{l \neq k} \lVert \Tilde{\Theta}_{k,:} - \Tilde{\Theta}_{l,:} \rVert_{2} = \sqrt{1/N_{k} + 1/\max \{N_{l} :l \neq k}\}$ and require 
    \begin{align}\label{eq:kmeans-requirement}
        (16 + 8 \varepsilon)\lVert U_{S} - U_{P} \rVert_{F}^{2} \leq N_{k} \delta_{k}^{2}
    \end{align}
    for all $k\in[K]$. 
    A sufficient condition for \eqref{eq:kmeans-requirement} to hold is 
    \begin{align}\label{eq:kmeans-condition}
        (16 + 8 \varepsilon)K\lVert U_{S} - U_{P} \rVert_{2}^{2} \leq 1.
    \end{align}
    Condition \eqref{eq:kmeans-condition} is implied by the assumption of Theorem \ref{thm:dist_U_S_U_P} on the number of observations given by \eqref{eq:T_cond_combined} since $K$ is a fixed constant. Therefore \citetSM[][Lemma 5.3]{lei2015consistency} implies there exists subsets $S_{k} \subset \zeta_{k}$ for $k = 1, \dots, K$ and a $K \times K$ permutation matrix $J$ such that $\hat{Z}_{\zeta,:}J = Z_{\zeta,:}$ where $\zeta = \cup_{k=1}^{K} (\zeta_{k} \setminus S_{k})$ and, upon choosing $\delta_{k} = \sqrt{1/N_{k}}$, Theorem \ref{thm:dist_U_S_U_P} implies 
    \begin{align}
        \sum_{k=1}^{K} \frac{|S_{k}|}{N_{k}} \leq&\, (16 + 8 \varepsilon)K\lVert U_{S} - U_{P} \rVert_{2}^{2} \nonumber \\ \lesssim&\,  \left(\frac{1}{\rho^{2} \lambda_{K}^{2}(B) }\right)^{2} \left\{ \frac{1}{(1 - \rho^{2})^{3}} \frac{K^{2}}{T} +  \left(1 + \frac{1}{ \sqrt{1 - \rho^{2}}}\right)^{2}\frac{KN}{T}\right\} \log^{\frac{1}{2}}(T+N) \\ &+   \frac{K^{4} }{N \alpha_{N}\lambda_{K}^{2}(B)} .
    \end{align}
    with probability at least $1 - (T+N)^{-\nu}$. 
    \qed

\subsection{Proof of Theorem~\ref{thm:strong-cov-estimation}} \label{proof:strong-cov-estimation}

 The assumption in~\eqref{eq:T_cond_combined} on $T$ implies that 
    \begin{align}\label{eq:strong-num_obs_rate}
        \frac{2\left\lVert S_{T} - \Gamma \right\rVert_{2} }{\Delta_{\Gamma}} &< \frac{1}{\log T}, 
    \end{align}
and therefore $\lambda_{K}(\Gamma) \geq \lVert S_{T} - \Gamma \rVert_{2}$. Let $E_{T} = S_{T} - \Gamma$. \citetSM[][Theorem 3.7]{cape2019two} thus yields
\begin{multline}
\|U_{S} - U W_{U}\|_{2, \infty}
\;\le\;
C\,\frac{\|(U_{\perp} U_{\perp}^{\prime})\,E_{T}(UU^{\prime})\|_{2, \infty}}{\lambda_{K}(\Gamma)}
+\,C\,\frac{\|(U_{\perp} U_{\perp}^{\prime})\,E_{T}(U_{\perp} U_{\perp}^{\prime})\|_{2, \infty}\;\|\sin\Theta(U_{S},U)\|_{2}}{\lambda_{K}(\Gamma)} \\
+\,C\,\frac{\|(U_{\perp} U_{\perp}^{\prime})\,\Gamma\,(U_{\perp} U_{\perp}^{\prime})\|_{2, \infty}\;\|\sin\Theta(U_{S},U)\|_{2}}{\lambda_{K}(\Gamma)}
+\,\|\sin\Theta(U_{S},U)\|_{2}^{2}\;\|U\|_{2, \infty}.
\label{eq:theorem3.7}
\end{multline}

We analyse each of the four terms in \eqref{eq:theorem3.7} in turn. By \citetSM[][Proposition 6.5]{cape2019two},
\begin{align}
\|(U U^{\prime})\,E_{T}(UU^{\prime})\|_{2, \infty}
&\;\le\;
\|U U^{\prime}\|_{\infty}\;\|E_{T} U\|_{2, \infty}.
\end{align}

The assumption of balanced communities implies bounded coherence, that is, $\lVert U \rVert_{2,\infty} \leq C\sqrt{K/N}$. Therefore, 
\begin{align}
\|U U^{\prime}\|_{\infty}
=\|I - U U^{\prime}\|_{\infty}
\le 1 + \sqrt{N}\,\|U U^{\prime}\|_{2, \infty}
< 1 + C\,\sqrt{K}\,.
\end{align}
We bound $E_{T}U \in \mathbb{R}^{N \times K}$ in two-to-infinity norm in the following manner:
\begin{align}
\|E_{T} U\|_{2, \infty}
\;\le\;
\sqrt{K}\,\|E_{T} U\|_{\max}
=\;
\sqrt{K}\,\max_{i \in [N],\,j \in [K]}\,|(E_{T}U)_{ij}|.
\end{align}
We can write 
\begin{align}
    E_{T}U &= \frac{1}{T}\sum_{t=1}^{T} (UU^{\prime} +U_{\perp}U_{\perp}^{\prime})\bm{X}_{t}\bm{X}_{t}^{\prime}U -U\Lambda_{\Gamma} 
    = \frac{1}{T}\sum_{t=1}^{T} \bm{X}_{t}^{U}\bm{Y}_{t}^{\prime} -U\Lambda_{\Gamma} + \frac{1}{T}\sum_{t=1}^{T} \bm{X}_{t}^{U_{\perp}}\bm{Y}_{t}^{\prime},
\end{align}
where we recall that $\bm{Y}_{t} \coloneqq U^{\prime}\bm{X}_{t}$, $\bm{X}_{t}^{U} \coloneqq UU^{\prime} \bm{X}_{t}$ and $\bm{X}_{t}^{U_{\perp}} \coloneqq U_{\perp}U_{\perp}^{\prime} \bm{X}_{t}$ with 
\begin{align}
    \bm{Y}_{t} \sim \mathcal{N}(0, \Lambda_{\Gamma}), & & 
    \bm{X}_{t}^{U} \sim \mathcal{N}(0,U\Lambda_{\Gamma} U^{\prime}), & & 
    \bm{X}_{t}^{U_{\perp}} \sim \mathcal{N}(0,U_{\perp}\Lambda_{\Gamma, \perp} U_{\perp}^{\prime}).
\end{align}
Since the product of subgaussian random variables is subexponential, we bound the following sum of dependent, centered subexponential random variables:
\begin{align}\label{eq:EU-sum}
    (E_{T}U)_{ij} = \frac{1}{T}\sum_{t=1}^{T} X_{i,t}^{U} Y_{j,t} -(U\Lambda_{\Gamma})_{ij} + \frac{1}{T}\sum_{t=1}^{T} X_{i,t}^{U_{\perp}} Y_{j,t}.
\end{align}
To bound the sums in \eqref{eq:EU-sum}, we will use a truncated version of the Bernstein-type inequality for geometrically absolutely regular random variables given by \citetSM[][Theorem~2]{merlevede2009bernstein}. Since $X_{i,t}^{U} Y_{j,t} -(U\Lambda_{\Gamma})_{ij}$ and $X_{i,t}^{U_{\perp}} Y_{j,t}$ are measureable functions of the geometrically absolutely regular sequence $(\bm{X}_{t})_{t \geq 1}$, then the sequences $(X_{i,t}^{U} Y_{j,t} -(U\Lambda_{\Gamma})_{ij})_{t \geq 1}$ and $(X_{i,t}^{U_{\perp}} Y_{j,t})_{t \geq 1}$ are also geometrically absolutely regular.

Define $Z_{t}^{(ij)} \coloneqq X_{i,t}^{U} Y_{j,t} -(U\Lambda_{\Gamma})_{ij}$. For ease of presentation, we suppress the superscript on $Z_{t}^{(ij)}$ in the remainder of the proof and write $Z_{t} = X_{i,t}^{U} Y_{j,t} -(U\Lambda_{\Gamma})_{ij}$. Note that $\mathbb{E}(Z_{t}) = 0$. Also, using Isserlis' theorem, we have 
\begin{align}
    \mathbb{E}(Z_{t}^{2}) &= (U\Lambda_{\Gamma}U^{\prime})_{ii}(\Lambda_{\Gamma})_{jj} + (U\Lambda_{\Gamma})_{ij}^{2} 
    \leq 2 \lambda_{1}^{2}(\Gamma) \lVert U \rVert_{2, \infty}^{2} 
    \leq C \lambda_{1}^{2}(\Gamma) \frac{K}{N},
\end{align}
where we used the assumption of balanced communities, which implies bounded coherence, in the last line.

Next, we find a high probability bound on $\sum_{t=1}^{T}Z_{t}$. The (univariate) subexponential norm and the (univariate) subgaussian norm of a random variable $Y$ are, respectively,
\begin{align}
\|Y^2\|_{\psi_1}
\;:=\;
\sup_{p \ge 1}\,
p^{-1}\,\Bigl(\mathbf{E}\bigl[(Y^2)^p\bigr]\Bigr)^{\tfrac{1}{p}},
& & 
\|Y\|_{\psi_2}
\;:=\;
\sup_{p \ge 1}\,
p^{-\tfrac12}\,\Bigl(\mathbf{E}\bigl[\lvert Y\rvert^p\bigr]\Bigr)^{\tfrac{1}{p}}.
\end{align}
By \citetSM[][Lemma 2.7.7]{vershynin2018high}, 
\begin{equation}
    \lVert X_{i,t}^{U} Y_{j,t} -( U\Lambda_{\Gamma})_{ij}  \rVert_{\psi_1} \leq 2 \lVert X_{i,t}^{U} Y_{j,t} \rVert_{\psi_1} \leq \lVert X_{i,t}^{U} \rVert_{\psi_2} \lVert Y_{j,t} \rVert_{\psi_2}. 
\end{equation}
Now 
\begin{align}
    \lVert X_{i,t}^{U} \rVert_{\psi_2} &\leq C \max_{1 \leq i \leq N} \sqrt{\text{Var}(X_{i,t}^{U})} 
    \leq C \sqrt{\lambda_{1}(\Gamma) \lVert U \rVert_{2, \infty}^{2}} 
    \leq C \sqrt{\lambda_{1}(\Gamma)} \sqrt{\frac{K}{N}}.
\end{align}
Since Var$(Y_{j,t}) = \lambda_{j}(\Gamma)$, we have $\lVert Y_{j,t} \rVert_{\psi_2} \leq C \sqrt{\lambda_{1}(\Gamma)}$. Therefore,
\begin{align}\label{eq:Z-psi1}
    \lVert X_{i,t}^{U} Y_{j,t} -( U\Lambda_{\Gamma})_{ij}  \rVert_{\psi_1} \leq C \lambda_{1}(\Gamma) \sqrt{\frac{K}{N}}.
\end{align}
This implies \citepSM[see, for example][Proposition 2.7.1]{vershynin2018high} that 
\begin{align}
    \mathbb{P}\left(|X_{i,t}^{U} Y_{j,t} -( U\Lambda_{\Gamma})_{ij}| > C \lambda_{1}(\Gamma) \sqrt{\frac{K}{N}} \log T \right) < q_{0},
\end{align}
where $q_{0}$ is small, for example, $q_{0} = T^{-\nu}$. For the truncated Bernstein-type inequality, we also need a bound on $|\mathbb{E}(Z_{t}) - \mathbb{E}(\tilde{Z}_{t})|$ where $\tilde{Z}_{t} \coloneqq Z_{t} \indfun \{|Z_{t}| \leq L\}$ with $L = C \lambda_{1}(\Gamma) \sqrt{K/N} \log T$. First note that since $\mathbb{E}(Z_{t}) = 0$, then $\mathbb{E}(\tilde{Z}_{t}) = - \mathbb{E}(Z_{t} \indfun \{ |Z_{t}| > L\}$. Therefore, 
\begin{align}
    |\mathbb{E}(Z_{t}) - \mathbb{E}(\tilde{Z}_{t})| &= |\mathbb{E}(Z_{t} \indfun \{ |Z_{t}| > L\}| 
    \leq \mathbb{E}(|Z_{t}| \indfun \{ |Z_{t}| > L\}) \\
    &\leq \sqrt{\mathbb{E}(Z_{t}^{2}) \mathbb{P}( |Z_{t}| > L)} 
    \leq \sqrt{C \lambda_{1}^{2}(\Gamma) \frac{K}{N} q_{0}},
\end{align}
where we used Jensen's inequality in the second line and the Cauchy-Schwarz inequality in the third.

We now compute the variance parameter
\begin{align}
    v^{2} = \mathbb{E}(Z_{t}) + 2\sum_{h=1}^{T-1} |\mathbb{E}(Z_{t}^{2}Z_{t+h})| .
\end{align}
By Isserlis' theorem, 
\begin{align}
    \mathbb{E}(Z_{t}Z_{t+h}) =&\, \mathbb{E}(X_{i,t}^{U}Y_{j,t}X_{i,t+h}^{U}Y_{j,t+h}) - (U\Lambda_{\Gamma})_{ij}^{2} \\
    =&\, \mathbb{E}(X_{i,t}^{U}Y_{j,t})\mathbb{E}(X_{i,t+h}^{U}Y_{j,t+h}) + \mathbb{E}(X_{i,t}^{U}X_{i,t+h}^{U})\mathbb{E}(Y_{j,t}Y_{j,t+h}) \\ &+ \mathbb{E}(X_{i,t}^{U}Y_{j,t+h})\mathbb{E}(Y_{j,t}X_{i,t+h}^{U}) - (U\Lambda_{\Gamma})_{ij}^{2}. 
\end{align}
We note that
\begin{align}
    &\mathbb{E}(X_{i,t}^{U}Y_{j,t}) = (U\Lambda_{\Gamma})_{ij}, 
    &|\mathbb{E}(X_{i,t}^{U}X_{i,t+h}^{U})| \leq \lambda_{1}^{h}(\Phi) \lambda_{1}(\Gamma) \lVert U \rVert_{2, \infty}^{2}, \\
    &|\mathbb{E}(Y_{j,t}Y_{j,t+h})| \leq \lambda_{1}^{h}(\Phi) \lambda_{1}(\Gamma), 
    &|\mathbb{E}(X_{i,t}^{U}Y_{j,t+h})| \leq \lambda_{1}^{h}(\Phi) \lambda_{1}(\Gamma) \lVert U \rVert_{2, \infty}.
\end{align}
Therefore, 
\begin{align}
    |\mathbb{E}(Z_{t}Z_{t+h})| \leq 2\lambda_{1}^{2h}(\Phi) \lambda_{1}^{2}(\Gamma) \lVert U \rVert_{2, \infty}^{2}
    \leq C \lambda_{1}^{2h}(\Phi) \lambda_{1}^{2}(\Gamma) \frac{K}{N}.
\end{align}
The variance parameter is therefore bounded as 
\begin{align}
    v^{2} &\leq C \lambda_{1}^{2}(\Gamma) \frac{K}{N} \sum_{h=0}^{T-1} \lambda_{1}^{2h}(\Phi) \leq C \lambda_{1}^{2}(\Gamma) \frac{K}{N} \frac{1}{1-\lambda_{1}^{2}(\Phi)} = C \frac{1}{(1 - \rho^{2})^{3}} \frac{K}{N}.
\end{align}
We observe that the variance of the truncated sum is of the same order as $v^{2}$ since $|\mathbb{E}(\tilde{Z}_{t}\tilde{Z}_{s})| \leq |\mathbb{E}(Z_{t}Z_{s})| + |\mathbb{E}(Z_{t}Z_{s} \indfun\{\mathcal{E}_{Z}^{\prime} \})|$ where $\mathcal{E}_{Z} \coloneqq \{|Z_{t} |<L\} \cap\{|Z_{s}|<L\}$ and by Jensen's inequality and the Cauchy-Schwarz inequality,
\begin{align}
    |\mathbb{E}(Z_{t}Z_{s} \indfun\{\mathcal{E}_{Z}^{\prime} \})| &\leq \mathbb{E}(Z_{t}Z_{s}) \indfun\{\mathcal{E}_{Z}^{\prime} \})  
    \leq \sqrt{\mathbb{E}(Z_{t}^{2} Z_{s}^{2}) \mathbb{P}(\mathcal{E}_{Z}^{\prime})}  
    \leq \{\mathbb{E}(Z_{t}^{4}) \mathbb{E}( Z_{s}^{4}) \}^{1/4} \sqrt{\mathbb{P}(\mathcal{E}_{Z}^{\prime})} \\ 
    &\leq \lVert Z_{t} \rVert_{\psi_{1}}^{2} \sqrt{\mathbb{P}(\mathcal{E}_{Z}^{\prime})} 
    \leq C \lambda_{1}^{2}(\Gamma) \frac{K}{N} \sqrt{q_{0}},
\end{align}
where we used \eqref{eq:Z-psi1} in the last line.
 Combining a union bound with a Bernstein inequality for the truncated random variables $\tilde{Z}_{t}$ \citepSM[][Theorem 2]{merlevede2009bernstein}, it follows that there exists a constant $C$  depending only on $c_{\beta}$ such that for $T > C (1 - \rho^{2}) \log^{7} T $, 
 \begin{align}\label{eq:parallel-2inf}
     \left|\frac{1}{T} \sum_{t=1}^{T} Z_{t}\right| \leq C \left( \frac{1}{1 - \rho^{2}} \right)^{3/2} \sqrt{\frac{K}{NT}} \log^{1/2}T.
 \end{align}
Next we bound the second term on the right hand side of \eqref{eq:EU-sum}. By \citetSM[][Lemma 2.7.7]{vershynin2018high}, 
\begin{equation}
    \lVert X_{i,t}^{U_{\perp}} Y_{j,t}  \rVert_{\psi_1} \leq \lVert X_{i,t}^{U_{\perp}} \rVert_{\psi_2} \lVert Y_{j,t} \rVert_{\psi_2}. 
\end{equation}
Now 
\begin{align}
    \lVert X_{i,t}^{U_{\perp}} \rVert_{\psi_2} &\leq C \max_{1 \leq i \leq N} \sqrt{\text{Var}(X_{i,t}^{U_{\perp}})} 
    \leq C \sqrt{\lambda_{K+1}(\Gamma) \lVert U_{\perp} \rVert_{2, \infty}^{2}} 
    \leq C \sqrt{\lambda_{K+1}(\Gamma)}.
\end{align}
Since Var$(Y_{j,t}) = \lambda_{j}(\Gamma)$, we have $\lVert Y_{j,t} \rVert_{\psi_2} \leq C \sqrt{\lambda_{1}(\Gamma)}$. Therefore,
\begin{align}
    \lVert X_{i,t}^{U_\perp} Y_{j,t}  \rVert_{\psi_1} \leq C \sqrt{\lambda_{K+1}(\Gamma)} \sqrt{\lambda_{1}(\Gamma)}.
\end{align}
This implies \citepSM[see, for example][Proposition 2.7.1]{vershynin2018high} that 
\begin{align}
    \mathbb{P}\left(|X_{i,t}^{U_{\perp}} Y_{j,t} | > C \sqrt{\lambda_{K+1}(\Gamma)} \sqrt{\lambda_{1}(\Gamma)} \log T \right) < q_{0},
\end{align}
where $q_{0}$ is small. Since $X_{i,t}^{U_{\perp}}$ and $Y_{j,t}$ are independent Gaussian random variables with $\mathbb{E}(X_{i,t}^{U_{\perp}}) = 0$ and $\mathbb{E}(Y_{j,t}) =0$, then the distribution of $X_{i,t}^{U_{\perp}} Y_{j,t}$ is symmetric around 0. Therefore the distribution of the truncated random variable $X_{i,t}^{U_{\perp}} Y_{j,t} \indfun\{\lvert X_{i,t}^{U_{\perp}} Y_{j,t} \rvert < L\}$ is also symmetric around 0. This implies that $\lvert \mathbb{E}(X_{i,t}^{U_{\perp}} Y_{j,t} ) -  \mathbb{E}(X_{i,t}^{U_{\perp}} Y_{j,t} \indfun\{\lvert X_{i,t}^{U_{\perp}} Y_{j,t} \rvert < L\}) \rvert = 0$ for any $L > 0$.

We note that $\mathbb{E}[(X_{i,t}^{U_{\perp}} Y_{j,t} )^{2}] = \mathbb{E}[(X_{i,t}^{U_{\perp}} )^{2}]\mathbb{E}( Y_{j,t}^{2}) \leq \lambda_{K+1}(\Gamma) \lambda_{1}(\Gamma)$. Using Isserlis' theorem and noting that $\mathbb{E}(X_{i,t}^{U_{\perp}} Y_{j,s}) = 0$, we have that 
\begin{align}
    \mathbb{E}(X_{i,t}^{U_{\perp}} Y_{j,t} X_{i,s}^{U_{\perp}} Y_{j,s}) \leq \lambda_{K+1}^{h}(\Phi) \lambda_{K+1}(\Gamma) \lambda_{1}^{h}(\Phi) \lambda_{1}(\Gamma).
\end{align}
Therefore, the variance parameter is bounded as 
\begin{align}
    v_{\perp}^{2} &= \mathbb{E}(X_{i,t}^{U_{\perp}} Y_{j,t}) + 2\sum_{h=1}^{T-1} |\mathbb{E}(X_{i,t}^{U_{\perp}} Y_{j,t} X_{i,t+h}^{U_{\perp}} Y_{j,t+h})| 
    \leq C \lambda_{1}(\Gamma) \lambda_{K+1}(\Gamma) \sum_{h=0}^{T-1} \lambda_{1}^{h}(\Phi) \lambda_{K+1}^{h}(\Phi) \\ 
    &\leq C \frac{\lambda_{1}(\Gamma) \lambda_{K+1}(\Gamma)}{1 - \lambda_{1}(\Phi) \lambda_{K+1}(\Phi)} 
    \lesssim C \frac{1}{1 - \rho^{2}}. 
\end{align}

For convenience, let $L_{\perp} \coloneqq C \sqrt{\lambda_{1}(\Gamma)}\log T$, $Z_{t}^{\perp} \coloneqq X_{i,t}^{U_{\perp}} Y_{j,t} $ and $\tilde{Z}_{t}^{\perp} \coloneqq X_{i,t}^{U_{\perp}} Y_{j,t} \indfun\{\lvert X_{i,t}^{U_{\perp}} Y_{j,t} \rvert < L_{\perp}\}$. The variance of the sum of truncated random variables $\tilde{Z}_{t}^{\perp}$ is of the same order as $v_{\perp}^{2}$ since, upon defining the event $\mathcal{E}_{Z,\perp} \coloneqq \{ |Z_{t, \perp}| < L_{\perp} \} \cap \{ |Z_{s, \perp}| < L_{\perp} \}$ we have 
\begin{align}
    |\mathbb{E}(Z_{t\perp} Z_{s\perp} \indfun \{\mathcal{E}_{Z,\perp}^{\prime} \})| &\leq \sqrt{\mathbb{E}(Z_{t,\perp}^{2} Z_{s,\perp}^{2})\mathbb{P}(\mathcal{E}_{Z,\perp}^{\prime})} 
    \leq \sqrt{\mathbb{E}(Z_{s,\perp}^{4})\mathbb{P}(\mathcal{E}_{Z,\perp}^{\prime})} \\ 
    &= \sqrt{\mathbb{E}\{(X_{i,t}^{U_{\perp}})^{4}\} \mathbb{E}\{Y_{j,t}^{4}\}\mathbb{P}(\mathcal{E}_{Z,\perp}^{\prime})}  
    \leq C \sqrt{q_{0}},
\end{align}
where we used the fact that $X_{i,t}^{U_{\perp}}$ and $Y_{j,t}$ are independent Gaussian random variables with bounded fourth moments. 
Again, combining a union bound with a Bernstein inequality for the truncated random variables $\tilde{Z}_{t}^{\perp}$ \citepSM[][Theorem 2]{merlevede2009bernstein}, it follows that there exists a constant $C$ depending only on $c_{\beta}$ such that for $T > C \log^{7} T $, 
 \begin{align}\label{eq:orth-2inf}
     \left|\frac{1}{T} \sum_{t=1}^{T} X_{i,t}^{U_{\perp}} Y_{j,t}\right| \leq C \left( \frac{1}{1 - \rho^{2}} \right)^{1/2} \sqrt{\frac{1}{T}} \log^{1/2}T.
 \end{align}

Combining \eqref{eq:parallel-2inf} and \eqref{eq:orth-2inf} we have that for any $i \in [N]$ and $j \in [K]$, there exists a constant $C$ such that for $T > C \log^{7} T $, 
\begin{align}
    |(E_{T}U)_{ij}| \lesssim \left\{ \left( \frac{1}{1 - \rho^{2}} \right)^{3/2} \sqrt{\frac{K}{NT}} + \left( \frac{1}{1 - \rho^{2}} \right)^{1/2} \sqrt{\frac{1}{T}} \right\} \log^{1/2}T
\end{align}
with high probability and therefore
\begin{align}
    \|(U_{\perp} U_{\perp}^{\prime})\,E_{T}(UU^{\prime})\|_{2, \infty} \lesssim \left\{ \left( \frac{1}{1 - \rho^{2}} \right)^{3/2} \sqrt{\frac{K^{3}}{NT}} + \left( \frac{1}{1 - \rho^{2}} \right)^{1/2} \sqrt{\frac{K^{2}}{T}} \right\} \log^{1/2}T
\end{align}
with high probability. 

We now bound the remaining terms in \eqref{eq:theorem3.7}. By Lemma~\ref{lemma:S_T-Gamma}, 
\begin{align}
    \lVert U_{\perp} U_{\perp}^{\prime} E_{T} U_{\perp} U_{\perp}^{\prime} \rVert_{2, \infty} \leq \lVert U_{\perp} U_{\perp}^{\prime} E_{T} U_{\perp} U_{\perp}^{\prime} \rVert_{2} \lesssim \sqrt{\frac{N}{T}}\log^{1/2}T.
\end{align}
By Theorem \ref{thm:dist_U_S_U_P}, with probability at least $1-(T+N)^{-\nu}$
\begin{align}
    \lVert \sin \Theta (U_{S},U) \rVert_{2} \lesssim \left\{ \frac{1}{\rho^{2} \lambda_{K}^{2}(B) }\right\}\left\{ \frac{1}{(1 - \rho^{2})^{3/2}} \sqrt{\frac{K}{T}} +  \left(1 + \frac{1}{\sqrt{1 - \rho^{2}}}\right)\sqrt{\frac{N}{T}}\right\} \log^{\frac{1}{2}}(T+N).
\end{align}
Note that $ \lVert U_{\perp} U_{\perp}^{\prime} \Gamma U_{\perp} U_{\perp}^{\prime} \rVert_{2, \infty} = \lVert U_{\perp} \Lambda_{\Gamma} U_{\perp}^{\prime} \rVert_{2, \infty} \leq \lambda_{K+1}(\Gamma) = O(1)$. Together, these observations yield that there exists $W_{U} \in \mathcal{O}(K \times K)$ and $C>0$ depending only on $c_{\beta}$, such that for $T \geq C\log^{7}T$ and $\nu \geq 8$, 
\begin{align}
    \|U_{S} - U W_{U}\|_{2, \infty} \leq &\,  \frac{C}{\rho^{2} \lambda_{K}^{2}(B) } \Bigg{[} \left\{ \frac{1}{(1 - \rho^{2})^{3/2}} \sqrt{\frac{K^{3}}{NT}} +  \frac{1}{\sqrt{1 - \rho^{2}}}  \sqrt{\frac{K^{2}}{T}} \right\} \log^{1/2}T \\ 
    &+   \left\{ \frac{1}{(1 - \rho^{2})^{3/2}} \frac{\sqrt{NK}}{T} +   \frac{1}{ \sqrt{1 - \rho^{2}}}\frac{N}{T}\right\} \log T \log^{\frac{1}{2}}(T+N) \\
    &+  \left\{ \frac{1}{(1 - \rho^{2})^{3/2}} \sqrt{\frac{K}{T}} +   \frac{1}{ \sqrt{1 - \rho^{2}}}\sqrt{\frac{N}{T}}\right\} \log^{\frac{1}{2}}(T+N) \\
    &+ \left\{ \frac{1}{(1 - \rho^{2})^{3}} \frac{K}{T} + \frac{1}{(1 - \rho^{2})^{2}} \frac{\sqrt{NK}}{T} +   \frac{1}{ 1 - \rho^{2}}\frac{N}{T}\right\} \sqrt{\frac{K}{N}}\log(T+N) \Bigg{]},
\end{align}
with probability at least $1 - (T+N)^{-\nu}$. Noting that \eqref{eq:T_cond_combined} implies $T \gtrsim \log^{7}T$ then yields the result. \qed

\subsection{Proof of Theorem~\ref{thm:strong-clustering}} \label{proof:strong-clustering}

We assume the event in Theorem \ref{thm:strong-cov-estimation} occurs and show that this, along with the assumption of $T \gtrsim N^{2}\log^{2}N$, implies the result. The $k$-means clustering of the rows of $U_{S}$ into $K$ blocks returns 
\begin{align}
\hat{C} \coloneqq \min_{C \in \mathcal{C}_{K}} \lVert C - U_{S} \rVert_{F}, \quad \text{where} \quad
\mathcal{C}_{K} = \{C \in \mathbb{R}^{n \times d} : C \text{ has } K \text{ distinct rows}\},
\end{align}
are the optimal cluster centroids. 
Let $W_{P} \coloneqq \argmin_{W \in \mathcal{O}(K \times K)}\lVert U_{P}W - U_{S} \rVert_{F}$. Since $U_{P}$ has $K$ distinct rows, it follows that
\begin{align}
\lVert \hat{C} - U_{S} \rVert_{F} \leq\lVert U_{P}W_{P} - U_{S} \rVert_{F} \leq \sqrt{N}\lVert U_{P}W_{P} - U_{S} \rVert_{2,\infty}.
\end{align}
By Theorem \ref{thm:strong-cov-estimation}, if $T \gtrsim N^{2}\log^{2}N$, with probability at least $1-(T+N)^{-\nu}$, $\min_{R \in \mathcal{O}(K \times K)}\lVert U_{S} - UR \rVert_{2,\infty} \lesssim 1/{\sqrt{N\log N}}$. By \citetSM[][Theorem 1]{agterberg2023overview}, with probability at least $1-N^{-\nu}$, $\min_{R \in \mathcal{O}^{K \times K}}\lVert U - U_{P}R\rVert_{2,\infty} \lesssim  \log N /{N\alpha_{N}^{1/2}} $. By the triangle inequality, 
\begin{align}
    \lVert U_{P}W_{P} - U_{S} \rVert_{2,\infty} \lesssim \frac{1}{\sqrt{N\log N}} + \frac{\log^{1/2}N}{N \alpha_{N}^{1/2}} \coloneqq \beta,
\end{align}
with probability at least $1 - (T+N)^{-\nu}$. Set $r=\beta \sqrt{N_{\text{min}}/N} \asymp \beta$ (since we assume balanced communities). Let $\mathcal{B}_{1}, \mathcal{B}_{2}, \ldots, \mathcal{B}_{K}$ be $L^{2}$-balls with radii $2r$ around the $K$ distinct rows of $U_{P}W_{P}$. The balls are disjoint if $\min_{i,j : z_{i} \neq z_{j}} \lVert (U_{P}W_{P})_{i,:} - (U_{P}W_{P})_{j,:}  \rVert_{2} > 4r$. By \citetSM[][Lemma 1]{agterberg2023overview}, $\min_{i,j : z_{i} \neq z_{j}} \lVert (U_{P}W_{P})_{i,:} - (U_{P}W_{P})_{j,:}  \rVert_{2} \gtrsim 1/\sqrt{N}$. Recall Assumption \ref{ass:regularity} that $\alpha_{N} \geq c_{0}\log N /N$ for some $c_{0}>0$. We can choose $c_{0}>0$ such that $\beta \lesssim 1/\sqrt{N} \lesssim \min_{i,j : z_{i} \neq z_{j}} \lVert (U_{P}W_{P})_{i,:} - (U_{P}W_{P})_{j,:}  \rVert_{2}$ and the balls are disjoint. 

Suppose there exists $k \in [K]$ such that $\mathcal{B}_{k}$ does not contain any rows of $\hat{C}$. Then $\lVert \hat{C} - U_{P}W_{P} \rVert_{F} > 2r\sqrt{N_{\min}}$, as for each $k$, no row of $\hat{C}$ is within $2r$ of the (at least $N_{\min}$) rows of $U_{P}W_P$ in $\mathcal{B}_{k}$. By the triangle inequality, this implies that
\begin{align}
\lVert \hat{C} - U_{S} \rVert_{F} &\geq \lVert \hat{C} - U_PW_P \rVert_{F} - \lVert U_S - U_PW_P \rVert_{F} > 2r\sqrt{N_{\min}} - \beta\sqrt{N} \\
&> 2\beta\sqrt{\frac{N}{N_{\min}}} \cdot \sqrt{N_{\min}} - \beta\sqrt{N} = \beta\sqrt{N},
\end{align}
a contradiction. Therefore, $\lVert \hat{C} - U_PW_P \rVert_{2, \infty} \leq 2r$. Hence, by the pigeonhole principle, each ball $\mathcal{B}_{k}$ contains precisely one distinct row of $\hat{C}$.

If $(U_PW_P)_{i} = (U_PW_P)_{j}$, then both $\hat{C}_{i}$ and $\hat{C}_{j}$ are elements of $\mathcal{B}_{z_{i}}$, and since there is exactly one distinct row of $\hat{C}$ in $\mathcal{B}_{k}$, $\hat{C}_{i} = \hat{C}_{j}$. Conversely, if $\hat{C}_{i} \neq \hat{C}_{j}$, then $(U_PW_P)_{i}$ and $(U_PW_P)_{j}$ are in disjoint balls $\mathcal{B}_{k}$ and $\mathcal{B}_{k'}$ for some $k, k' \in [K]$, implying that $(U_PW_P)_{i} \neq (U_PW_P)_{j}$. Thus, $(U_PW_P)_{i} = (U_PW_P)_{j}$ if and only if $\hat{C}_{i} = \hat{C}_{j}$, proving the theorem. \qed 

\section{Laplacian model}\label{sec:appendix-laplacian}

\subsection{Proof of Lemma~\ref{lemma:eigs-laplacian}} \label{sec:appendix-eigs-lap} 
Following \citetSM[][Lemma 3.1]{rohe2011spectral}, we write $\mathcal{L}(P) = \alpha_{N}ZB_{L}Z^{\prime}$ where $B_{L} \coloneqq D_{B}^{-1/2}BD_{B}^{-1/2}$ and $D_{B} \coloneqq \text{diag}(BZ^{\prime}\bm{1}) \in \mathbb{R}^{K \times K}$. 

By Ostrowski's theorem \citepSM[see][Theorem 4.5.9]{horn2012matrix}, for $r \in [K]$,
    \begin{align}
        \lambda_{r}(\mathcal{L}(P)) &= \lambda_{r}(\alpha_{N}ZB_{L}Z^{\prime}) = \alpha_{N}\theta_{r} \lambda_{r}(B_{L}),
    \end{align}
    where $\theta_{r} \in [N_{\text{min}},N_{\max}]$ with $N_{\text{min}}$ and $N_{\text{max}}$ being the smallest and largest block size, respectively. Under Assumption \ref{ass:regularity-laplacian}, the blocks are balanced and therefore, 
    \begin{align}
        \lambda_{r}(\mathcal{L}(P)) \asymp \frac{  \alpha_{N} N \lambda_{r}(B_{L})}{K}.
    \end{align}
    Let $\delta \coloneqq \min_{i\in [N]} \sum_{j=1}^{N}P_{ij}$. Then $D_{B} \succeq \delta I_{K}$ implies $\lambda_{r}(B_{L}) \asymp \lambda_{r}(B)/\delta$. Furthermore, $\delta \asymp \alpha_{N}N$ implies $\lambda_{r}(\mathcal{L}(P)) \asymp \lambda_{r}(B)$. 
    By \citetSM[][Theorem 2]{chung2011spectra}, for any $s>0$, with probability at least $1 - N^{-s}$,
    \begin{align}
        |\lambda_{r}(\mathcal{L}(A)) - \lambda_{r}(\mathcal{L}(P))| \lesssim \frac{\log^{1/2}N}{\sqrt{\delta}} \asymp \frac{\log^{1/2}N}{\sqrt{N\alpha_{N}}}. 
    \end{align}
    Thus, since $\lambda_{r}(B)$ is bounded away from zero for $r \in [K]$, with probability at least $1 - N^{-s}$,
    \begin{align}
        \lambda_{r}(\mathcal{L}(A)) \asymp \lambda_{r}(B).
    \end{align}
    Therefore, with probability at least $1 - N^{-s}$,
    \begin{align}
        \lambda_{r}(\Phi_{L}) = \rho \lambda_{r}(\mathcal{L}(A)) 
        \asymp \rho \lambda_{r}(B). 
    \end{align}
    For $r = K+1,\dots,N$, since $\mathcal{L}(P)$ is rank $K$, $\lambda_{r}(\mathcal{L}(P)) = 0$. Thus with probability at least $1 - N^{-s}$, $|\lambda_{r}(\mathcal{L}(A))| \lesssim \log^{1/2}N/\sqrt{N \alpha_{N}}$ and 
    \begin{align}
        \lambda_{r}(\Phi_{L}) &\asymp \frac{\rho\log^{1/2}N}{\sqrt{N \alpha_{N}}} 
        \xrightarrow[]{N \to \infty} 0. 
    \end{align}
    Therefore, with probability at least $1 - N^{-s}$, $\lambda_{r}(\Phi_{L}) = o(1)$ for $r = K+1,\dots,N$. \qed

\subsection{Proof of Theorem \ref{thm:dist_U_S_U_P-lap}}
    The Davis-Kahan sin$\Theta$ theorem \citepSM{davis1970rotation} implies that if 
    \begin{align}\label{eq:davis-kahan-assumption-lap}
        \left\lVert S_{T} - \Gamma \right\rVert_{2} < (1 - 1/\sqrt{2})\Delta_{\Gamma},
    \end{align}
    then 
    \begin{equation}
        \text{dist} (U_{S}, U) \leq \frac{2\left\lVert S_{T} - \Gamma \right\rVert_{2} }{\Delta_{\Gamma}}, 
    \end{equation}
    where $U$ is the $N \times K$ matrix whose columns are simultaneously the eigenvectors of $\Gamma$ and $\Phi_{L}$ corresponding to the $K$ largest-in-magnitude eigenvalues.
    The assumption on $T$ specified in the statement of the theorem implies that 
    \begin{align}\label{eq:num_obs_rate-lap}
        \frac{2\left\lVert S_{T} - \Gamma \right\rVert_{2} }{\Delta_{\Gamma}} &< \frac{1}{\log T} 
    \end{align}
    and therefore \eqref{eq:davis-kahan-assumption-lap} holds under the assumptions of the theorem.  Following the proof of Lemma \ref{lemma:S_T-Gamma}, with Lemma \ref{lemma:eigs} replaced by Lemma \ref{lemma:eigs-laplacian}, we have that, for any $\nu \geq 8$, with probability exceeding $1 - (T+N)^{-\nu}$, 
    \begin{align}\label{eq:S-Gamma-rate-lap}
        \left\lVert S_{T} - \Gamma \right\rVert_{2} &\lesssim \Bigg{[} \frac{1}{(1 - \rho^{2})^{3/2}} \sqrt{\frac{K}{T}} +  \left\{1 + \frac{1}{\sqrt{ 1 - \rho^{2}}}\right\}\sqrt{\frac{N}{T}} \Bigg{]} \log^{\frac{1}{2}}(T+N).
    \end{align}
    By Assumption \ref{ass:regularity}, Lemma \ref{lemma:eigs-laplacian}, and \citetSM[][Proposition 4.1]{martin2024nirvar_arxiv},
    $1/\Delta_{\Gamma} \asymp (\kappa^{2} - \rho^{2})/\rho^{2} \lesssim 1/\rho^{2}\lambda_{K}^{2}(B)$ with high probability.
Therefore, with probability exceeding $1 - (T+N)^{-\nu}$
\begin{multline}\label{eq:dist_U_S-lap}
        \text{dist}(U_{S},U) \lesssim \left\{ \frac{1}{\rho^{2} \lambda_{K}^{2}(B) }\right\}\left\{ \frac{1}{(1 - \rho^{2})^{3/2}} \sqrt{\frac{K}{T}} +  \left(1 + \frac{1}{ \sqrt{1 - \rho^{2}}}\right)\sqrt{\frac{N}{T}}\right\} \log^{\frac{1}{2}}(T+N).
    \end{multline}
By \citetSM[][Proposition 6.3]{cape2019two}, $\lVert U - \tilde{U}_{P}R \rVert_{2} \leq \sqrt{N} \lVert U - \tilde{U}_{P}R \rVert_{2,\infty}$ for any $R \in \mathcal{O}(K)$. Under Assumptions \ref{ass:regularity} and \ref{ass:regularity-laplacian}, by \citetSM[][Theorem 4.1]{ke2025optimal}, with probability exceeding $1 - N^{-4}$,
\begin{align}
    \min_{R \in \mathcal{O}(K)}\lVert U - \tilde{U}_{P}R \rVert_{2,\infty} \lesssim \frac{\log^{1/2}N}{N \alpha_{N}^{3/4} \lambda_{K}(B)},
\end{align}
and therefore, 
\begin{align}\label{eq:dist_U_UP-lap}
    \text{dist}(U,\tilde{U}_{P}) \lesssim \frac{\log^{1/2}N}{\sqrt{N} \alpha_{N}^{3/4} \lambda_{K}(B)}.
\end{align}
Note that \citetSM[][Condition 2.1(a)]{ke2025optimal} holds by Assumption \ref{ass:regularity} of balanced communities and \citetSM[][Condition 2.1(d)]{ke2025optimal} holds since we do not consider degree-corrected nodes. 
By the triangle inequality,
\begin{align}\label{eq:triangle-lap}
    \text{dist}(U_{S},\tilde{U}_{P}) \leq \text{dist}(U_{S},U) + \text{dist}(U,\tilde{U}_{P}).
\end{align}
The result follows from substituting \eqref{eq:dist_U_S-lap} and \eqref{eq:dist_U_UP-lap} into \eqref{eq:triangle-lap}. \qed

\subsection{Proof of Theorem~\ref{thm:kmeans-lap}} \label{proof:kmeans-lap}
    We apply \citetSM[][Lemma 5.3]{lei2015consistency} to $U_{S}$ and $\tilde{U}_{P}R^{*}$ where $R^{*} = \argmin_{R \in \mathbb{R}^{K \times K}} \lVert U_{S} - \tilde{U}_{P}R \rVert_{2}$. By \citetSM[][Lemma 3.1]{rohe2011spectral} (see also the discussion of \citetSM[][Lemma 3]{agterberg2023overview}), $\tilde{U}_{P}R^{*}  = Z \Tilde{\Theta}$ where $\lVert \Tilde{\Theta}_{k,:} - \Tilde{\Theta}_{l,:} \rVert_{2} = \sqrt{1/N_{k} + 1/N_{l}}$. We choose $\delta_{k} \leq \min_{l \neq k} \lVert \Tilde{\Theta}_{k,:} - \Tilde{\Theta}_{l,:} \rVert_{2} = \sqrt{1/N_{k} + 1/\max \{N_{l} :l \neq k}\}$ and require 
    \begin{align}\label{eq:kmeans-requirement-lap}
        (16 + 8 \varepsilon)\lVert U_{S} - \tilde{U}_{P} \rVert_{F}^{2} \leq N_{k} \delta_{k}^{2}
    \end{align}
    for all $k\in[K]$. 
    A sufficient condition for \eqref{eq:kmeans-requirement-lap} to hold is 
    \begin{align}\label{eq:kmeans-condition-lap}
        (16 + 8 \varepsilon)K\lVert U_{S} - U_{P} \rVert_{2}^{2} \leq 1.
    \end{align}
    Condition \eqref{eq:kmeans-condition-lap} is implied by the assumption on $T$ in Theorem \ref{thm:dist_U_S_U_P-lap} since $K$ is a fixed constant. Therefore \citetSM[][Lemma 5.3]{lei2015consistency} implies there exists subsets $S_{k} \subset \zeta_{k}$ for $k = 1, \dots, K$ and a $K \times K$ permutation matrix $J$ such that $\hat{Z}_{\zeta,:}J = Z_{\zeta,:}$ where $\zeta = \cup_{k=1}^{K} (\zeta_{k} \setminus S_{k})$ and, upon choosing $\delta_{k} = \sqrt{1/N_{k}}$, Theorem \ref{thm:dist_U_S_U_P-lap} implies 
    \begin{align}
        \sum_{k=1}^{K}&\, \frac{|S_{k}|}{N_{k}} \leq (16 + 8 \varepsilon)K\lVert U_{S} - U_{P} \rVert_{2}^{2} \nonumber \\ \lesssim&\,  \left(\frac{1}{\rho^{2} \lambda_{K}^{2}(B) }\right)^{2} \left\{ \frac{1}{(1 - \rho^{2})^{3}} \frac{K^{2}}{T} +  \left(1 + \frac{1}{ \sqrt{1 - \rho^{2}}}\right)^{2}\frac{KN}{T}\right\} \log (T+N) +   \frac{\log N }{N \alpha_{N}^{3/2}\lambda_{K}^{2}(B)} .
    \end{align}
    with probability at least $1 - (T+N)^{-\nu}$. 
    \qed

\subsection{Proof of Theorem~\ref{thm:strong-clustering-lap}} \label{proof:strong-clustering-lap}

Following the proof of Theorem \ref{thm:strong-cov-estimation}, under the assumptions of Theorem \ref{thm:dist_U_S_U_P-lap}, there exists a $K \times K$ orthogonal matrix $W_{\tilde{U}}$ and a constant $C > 0$ depending only on $c_{\beta}$ such that, with probability $1 - (T+N)^{-\nu}$,
    \begin{align}\label{eq:strong-cov-bound-lap}
    &\|U_{S} - U W_{\tilde{U}}\|_{2 ,\infty} \lesssim   \frac{1}{\rho^{2} \lambda_{K}^{2}(B) } \Bigg{[} \left\{ \frac{1}{(1 - \rho^{2})^{3/2}} \sqrt{\frac{K^{3}}{NT}} +  \frac{1}{\sqrt{1 - \rho^{2}}}  \sqrt{\frac{K^{2}}{T}} \right\} \log^{1/2}T \\ 
    &\qquad+  \left\{ \frac{1}{(1 - \rho^{2})^{3/2}} \frac{\sqrt{NK}}{T} +   \frac{1}{ \sqrt{1 - \rho^{2}}}\frac{N}{T}\right\} \log T \log^{\frac{1}{2}}(T+N) \\
    &\qquad+  \left\{ \frac{1}{(1 - \rho^{2})^{3/2}} \sqrt{\frac{K}{T}} +   \frac{1}{ \sqrt{1 - \rho^{2}}}\sqrt{\frac{N}{T}}\right\} \log^{\frac{1}{2}}(T+N) \\
    &\qquad+ \left\{ \frac{1}{(1 - \rho^{2})^{3}} \frac{K}{T} + \frac{1}{(1 - \rho^{2})^{2}} \frac{\sqrt{NK}}{T} +   \frac{1}{ 1 - \rho^{2}}\frac{N}{T}\right\}\sqrt{\frac{K}{N}} \log (T+N) \Bigg{]}
\end{align}
for any $\nu \geq 8$. 

Under the additional assumption that $T \gtrsim N^{2}\log^{2}N$, \eqref{eq:strong-cov-bound-lap} implies that with probability at least $1 - (T+N)^{-\nu}$,
$\min_{R \in \mathcal{O}^{K \times K}}\lVert U_{S} - UR\rVert_{2,\infty} \lesssim 1/\sqrt{N\log N}$. \citetSM[][Theorem 4.1]{ke2025optimal} implies $\min_{R \in \mathcal{O}^{K \times K}}\lVert U - \tilde{U}_{P}R\rVert_{2,\infty} \lesssim \log^{1/2}N/N\alpha_{N}^{3/4}\lambda_{K}(B)$ with probability at least $1-N^{-\nu}$. By the triangle inequality, with probability at least $1 - (T+N)^{-\nu}$,
\begin{align}\label{eq:strong-proof-rate-lap}
    \min_{R \in \mathcal{O}^{K \times K}}\lVert U_{S} - \tilde{U}_{P}R \rVert_{2,\infty} &\leq \min_{R \in \mathcal{O}^{K \times K}}\lVert U_{S} - UR\rVert_{2,\infty} + \min_{R \in \mathcal{O}^{K \times K}}\lVert U - \tilde{U}_{P}R\rVert_{2,\infty} \\ &\lesssim \frac{1}{\sqrt{N\log N}} + \frac{\log^{1/2}N}{N \alpha_{N}^{3/4}\lambda_{K}(B)} \coloneqq \beta_{L}.
\end{align}

The $k$-means clustering of the rows of $U_{S}$ into $K$ blocks returns 
\begin{align}
\hat{C} \coloneqq \min_{C \in \mathcal{C}_{K}} \lVert C - U_{S} \rVert_{F}, \quad \text{where}\quad
\mathcal{C}_{K} = \{C \in \mathbb{R}^{n \times d} : C \text{ has } K \text{ distinct rows}\},
\end{align}
are the optimal cluster centroids. 
Let $\tilde{W}_{P} \coloneqq \argmin_{W \in \mathcal{O}(K \times K)}\lVert \tilde{U}_{P}W - U_{S} \rVert_{F}$. Since $\tilde{U}_{P}$ has $K$ distinct rows, it follows that
\begin{align}
\lVert \hat{C} - U_{S} \rVert_{F} \leq\lVert \tilde{U}_{P}\tilde{W}_{P} - U_{S} \rVert_{F} \leq \sqrt{N}\lVert \tilde{U}_{P}\tilde{W}_{P} - U_{S} \rVert_{2,\infty}.
\end{align}
Set $r=\beta_{L} \sqrt{N_{\text{min}}/N} \asymp \beta_{L}$ (since we assume balanced communities). Let $\mathcal{B}_{1}, \mathcal{B}_{2}, \ldots, \mathcal{B}_{K}$ be $L^{2}$-balls with radii $2r$ around the $K$ distinct rows of $\tilde{U}_{P}\tilde{W}_{P}$. The balls are disjoint if $\min_{i,j : z_{i} \neq z_{j}} \lVert (\tilde{U}_{P}\tilde{W}_{P})_{i,:} - (\tilde{U}_{P}\tilde{W}_{P})_{j,:}  \rVert_{2} > 4r$. By \citetSM[][Lemma 3]{rohe2011spectral}, $\min_{i,j : z_{i} \neq z_{j}} \lVert (\tilde{U}_{P}\tilde{W}_{P})_{i,:} - (\tilde{U}_{P}\tilde{W}_{P})_{j,:}  \rVert_{2} \gtrsim 1/\sqrt{N}$. Recall Assumption \ref{ass:regularity} that $\alpha_{N} \geq c_{0}\log N /N$ for some $c_{0}>0$. We can choose $c_{0}>0$ such that $\beta_{L} \lesssim 1/\sqrt{N} \lesssim \min_{i,j : z_{i} \neq z_{j}} \lVert (\tilde{U}_{P}\tilde{W}_{P})_{i,:} - (\tilde{U}_{P}\tilde{W}_{P})_{j,:}  \rVert_{2}$ and the balls are disjoint. 

Suppose there exists $k \in [K]$ such that $\mathcal{B}_{k}$ does not contain any rows of $\hat{C}$. Then $\lVert \hat{C} - \tilde{U}_{P}\tilde{W}_{P} \rVert_{F} > 2r\sqrt{N_{\min}}$, as for each $k$, no row of $\hat{C}$ is within $2r$ of the (at least $N_{\min}$) rows of $\tilde{U}_{P}\tilde{W}_{P}$ in $\mathcal{B}_{k}$. By the triangle inequality, this implies that
\begin{align}
\lVert \hat{C} - U_{S} \rVert_{F} &\geq \lVert \hat{C} - \tilde{U}_{P}\tilde{W}_{P} \rVert_{F} - \lVert U_S - \tilde{U}_{P}\tilde{W}_{P} \rVert_{F} > 2r\sqrt{N_{\min}} - \beta_{L}\sqrt{N} \notag \\
&> 2\beta_{L}\sqrt{\frac{N}{N_{\min}}} \cdot \sqrt{N_{\min}} - \beta_{L}\sqrt{N} = \beta_{L}\sqrt{N},
\end{align}
a contradiction. Therefore, $\lVert \hat{C} - \tilde{U}_{P}\tilde{W}_{P} \rVert_{2, \infty} \leq 2r$. Hence, by the pigeonhole principle, each ball $\mathcal{B}_{k}$ contains precisely one distinct row of $\hat{C}$.

If $(\tilde{U}_{P}\tilde{W}_{P})_{i} = (\tilde{U}_{P}\tilde{W}_{P})_{j}$, then both $\hat{C}_{i}$ and $\hat{C}_{j}$ are elements of $\mathcal{B}_{z_{i}}$, and since there is exactly one distinct row of $\hat{C}$ in $\mathcal{B}_{k}$, $\hat{C}_{i} = \hat{C}_{j}$. Conversely, if $\hat{C}_{i} \neq \hat{C}_{j}$, then $(\tilde{U}_{P}\tilde{W}_{P})_{i}$ and $(\tilde{U}_{P}\tilde{W}_{P})_{j}$ are in disjoint balls $\mathcal{B}_{k}$ and $\mathcal{B}_{k'}$ for some $k, k' \in [K]$, implying that $(\tilde{U}_{P}\tilde{W}_{P})_{i} \neq (\tilde{U}_{P}\tilde{W}_{P})_{j}$. Thus, $(\tilde{U}_{P}\tilde{W}_{P})_{i} = (\tilde{U}_{P}\tilde{W}_{P})_{j}$ if and only if $\hat{C}_{i} = \hat{C}_{j}$, proving the theorem. \qed 

\section{General error covariance, weighted directed graphs}

\begin{lemma}\label{lemma:VAR-autocovariance-bound}
Let $(\bm{f}_{t})_{t \in \mathbb{Z}}$ be any stationary VAR$(p)$ process $\bm{f}_{t} = \Phi_{1}\bm{f}_{t-1} + \cdots + \Phi_{p}\bm{f}_{t-p} + \bm{u}_{t}$ with $p \in \mathbb{N}$. The noise process is assumed to be white noise with $\bm{u}_{t}$ being continuous random vectors satisfying $\mathbb{E}(\bm{u}_{t}) = 0$, $\Sigma_{u} = \mathbb{E}(\bm{u}_{t}\bm{u}_{t}^{\prime})$ is non-singular, and $\bm{u}_{t}$ and $\bm{u}_{s}$ are independent for $t \neq s$. Then the lag-$h$ autocovariance $\Gamma_{f}(h) = \mathbb{E}(\bm{f}_{t} \bm{f}_{t+h}^{\prime})$ satisfies 
    \begin{align}
        \lVert \Gamma_{f}(h) \rVert_{2} \lesssim \lVert \Sigma_{u}\rVert_{2} \frac{\rho^{h}}{1 - \rho^{2}},
    \end{align}
    where $0<\rho<1$ is the spectral radius of the companion matrix 
    \begin{equation}
\Phi = 
\begin{bmatrix}
\Phi_1 & \Phi_2 & \Phi_3 & \cdots & \Phi_{p} \\
I_r & 0   & 0   & \cdots & 0 \\
0   & I_r & 0   & \cdots & 0 \\
\vdots & \ddots & \ddots & \ddots & \vdots \\
0 & \cdots & 0 & I_r & 0
\end{bmatrix}.
\end{equation}
\end{lemma}

\subsection{Proof of Lemma \ref{lemma:VAR-autocovariance-bound}}
We start by writing the VAR$(p)$ process as a VAR$(1)$ process in companion matrix form. Let $\bm{F}_{t} \coloneqq (\bm{f}_{t}^{\prime},\dots,\bm{f}_{t-p+1}^{\prime})^{\prime} \in \mathbb{R}^{Np \times 1}$, $\bm{U}_{t} \coloneqq (\bm{u}_{t}^{\prime},0,\dots,0)^{\prime} \in \mathbb{R}^{Np \times 1}$, $\Sigma_{U} \coloneqq \mathbb{E}(\bm{U}_{t} \bm{U}_{t}^{\prime}) \in \mathbb{R}^{Np \times Np}$.
Then $\bm{F}_{t} = \Phi\bm{F}_{t-1} + \bm{U}_{t}$, and the autocovariances are given by \citepSM[see][Equation 2.1.10, for example]{lutkepohl2005new}
\begin{equation}
    \Gamma_{F}(h) = \sum_{j=0}^{\infty} \Phi^{j+h} \Sigma_{U} (\Phi^{j})^{\prime},
\end{equation}
with $\Gamma_{f}(h) = J \Gamma_{F}(h) J^{\prime}$ where $J \coloneqq (I_{r} |0| \cdots |0) \in \mathbb{R}^{r \times Np}$. Using the triangle inequality, and noting that $\lVert J\rVert_{2} = 1$ and $\lVert \Sigma_{U}\rVert_{2} = \lVert \Sigma_{u}\rVert_{2}$,
\begin{align}
    \lVert \Gamma_{f}(h) \rVert_{2} \leq \lVert J\rVert_{2}^{2} \lVert \Gamma_{F}(h) \rVert_{2} \leq \sum_{j=0}^{\infty} \lVert \Phi^{j+h} \rVert_{2} \lVert \Sigma_{u}\rVert_{2} \lVert \Phi^{j} \rVert_{2}.
\end{align}
By assumption the spectral radius of $\Phi$ satisfies $\rho < 1$. Gelfand's formula \citepSM[see][Corollary 5.6.14, for example]{horn2012matrix} implies that $\rho = \lim_{j \to \infty} \lVert \Phi^{j}\rVert_{2}^{1/j}$. Therefore, for any $\epsilon \in (0, 1 - \rho)$, there exists $N_{0} \in \mathbb{N}$ such that for $j > N_{0}$, $\lVert \Phi^{j} \rVert_{2} < (\rho + \epsilon)^{j}$. Thus,
\begin{align}
    \lVert \Gamma_{f}(h) \rVert_{2} &\leq \lVert \Sigma_{u}\rVert_{2} \sum_{j=0 }^{N_{0}} \lVert \Phi \rVert_{2}^{2j+h} +\lVert \Sigma_{u}\rVert_{2} (\rho+\epsilon)^{h} \sum_{j=N_{0} }^{\infty} (\rho+\epsilon)^{2j} .
\end{align}
Since $\lVert \Phi \rVert_{2} < \infty$, the first term above, which is a finite sum, is finite. Since the second term is a sum of positive terms and $\rho + \epsilon < 1$,
\begin{align}
    \sum_{j=N_{0} }^{\infty} (\rho+\epsilon)^{2j} < \sum_{j=0 }^{\infty} (\rho+\epsilon)^{2j} = \frac{1}{1-(\rho + \epsilon)^{2}} = \frac{1}{1 - \rho^{2}} + O(\epsilon).
\end{align}
In conclusion,
\begin{align}\label{eq:factor-decay}
    \lVert \Gamma_{f}(h) \rVert_{2} \lesssim \lVert \Sigma_{u}\rVert_{2} \frac{\rho^{h}}{1 - \rho^{2}}.
\end{align}
\qed

\subsection{Proof of Lemma \ref{lemma:S_T-Gamma-general-errors}}
The proof is identical to that of Lemma \ref{lemma:S_T-Gamma} but with different bounds on $v^{2}$, $v_{\perp}^{2}$, and $v_{G}^{2}$. Starting with $v^{2}$, Isserlis' theorem gives 
\begin{align}
    \mathbb{E}(W_{t}W_{t+h}) &= \mathbb{E}(\bm{Y}_{t} \bm{Y}_{t}^{\prime} \bm{Y}_{t+h} \bm{Y}_{t+h}^{\prime} ) - \Lambda_{\Gamma}^{2} \\ 
    &= \mathbb{E}(\bm{Y}_{t} \bm{Y}_{t}^{\prime})\mathbb{E}(\bm{Y}_{t+h} \bm{Y}_{t+h}^{\prime}) + \mathbb{E}(\bm{Y}_{t} \bm{Y}_{t+h}^{\prime})\mathbb{E}(\bm{Y}_{t}^{\prime} \bm{Y}_{t+h}) + \mathbb{E}(\bm{Y}_{t} \bm{Y}_{t+h}^{\prime})\mathbb{E}(\bm{Y}_{t} \bm{Y}_{t+h}^{\prime}) - \Lambda_{\Gamma}^{2} \\
    &= \Lambda_{\Gamma}^{2} + U_{\Gamma}^{\prime}\Gamma(h)U_{\Gamma} \text{tr}(U_{\Gamma}^{\prime}\Gamma(h)U_{\Gamma}) + (U_{\Gamma}^{\prime}\Gamma(h)U_{\Gamma})^{2} - \Lambda_{\Gamma}^{2}. 
\end{align}
Noting that $\text{tr}(U_{\Gamma}^{\prime}\Gamma(h)U_{\Gamma}) \leq K\lVert \Gamma(h) \rVert_{2}$ then gives $\lVert \mathbb{E}(W_{t}W_{t+h}) \rVert_{2} \leq (K+1) \lVert \Gamma(h) \rVert_{2}^{2}$. By Lemma \ref{lemma:VAR-autocovariance-bound}, $\lVert \Gamma(h) \rVert_{2} \lesssim \lVert\Sigma\rVert_{2} \rho^{h}/(1-\rho^{2})$ and therefore $\lVert \mathbb{E}(W_{t}W_{t+h}) \rVert_{2} \lesssim K\lVert\Sigma\rVert_{2}^{2} \rho^{2h}/(1-\rho^{2})^{2}$. Following the same algebra as in Lemma \ref{lemma:S_T-Gamma} then gives 
\begin{align}
    v^{2} \lesssim \frac{K\lVert\Sigma\rVert_{2}^{2}}{(1 - \rho^{2})^{3}}.
\end{align}
Similarly, we have 
\begin{align}
    v_{\perp}^{2} \lesssim \frac{N\lVert\Sigma\rVert_{2}^{2}}{(1 - \rho^{2})^{3}}.
\end{align}
To bound $v_{G}^{2}$, we note that
\begin{align}
    \mathbb{E}(G_{t}G_{t+h}) 
    &= \mathbb{E}( (\bm{X}_{t}^{U_{\perp}})^{\prime}\bm{X}_{t+h}^{U_{\perp}}) \mathbb{E}(\bm{X}_{t}^{U} (\bm{X}_{t+h}^{U})^{\prime}) + \mathbb{E}( (\bm{X}_{t}^{U})^{\prime}\bm{X}_{t+h}^{U}) \mathbb{E}(\bm{X}_{t}^{U_{\perp}} (\bm{X}_{t+h}^{U_{\perp}})^{\prime}) \\
    &= \text{tr}\{U_{\Gamma,\perp}^{\prime}\Gamma(h)U_{\Gamma,\perp}\}U_{\Gamma}^{\prime}\Gamma(h)U_{\Gamma} + \text{tr}\{U_{\Gamma}^{\prime}\Gamma(h)U_{\Gamma}\}U_{\Gamma,\perp}^{\prime}\Gamma(h)U_{\Gamma,\perp},
\end{align}
and therefore $\lVert \mathbb{E}(G_{t}G_{t+h}) \rVert_{2} \leq N \lVert \Gamma(h) \rVert_{2}^{2}$. Using Lemma \ref{lemma:VAR-autocovariance-bound} and repeating the variance parameter computation from Lemma \ref{lemma:S_T-Gamma} gives
\begin{align}
    v_{G}^{2} \lesssim \frac{N\lVert\Sigma\rVert_{2}^{2}}{(1 - \rho^{2})^{3}}.
\end{align}
Repeating the proof of Lemma \ref{lemma:S_T-Gamma} with these bounds on $v^{2}$, $v_{\perp}^{2}$, and $v_{G}^{2}$ then yields the result. \qed

\subsection{Proof of Lemma \ref{lemma:singular-angle}}
We will compare $\Gamma$ to the matrix $\Xi \coloneqq (I_{N} - \Phi \Phi^{\prime})^{-1}$. Since $\Xi$ has eigenvectors $U_{\Phi,N} \coloneqq (U_{\Phi}|U_{\Phi,\perp}) \in \mathbb{R}^{N \times N}$ and eigenvalues $\Lambda_{\Xi} = \text{diag}(1/\{1-\sigma_{K+1}^{2}(\Phi)\},\dots,1/\{1-\sigma_{N}^{2}(\Phi)\})$, by the Davis-Kahan sin$\Theta$ theorem, $$\text{dist}(U_{\Gamma},U_{\Phi})~\leq~\lVert E \rVert_{2}/\{\lambda_{K}(\Xi) - \lambda_{K+1}(\Xi)\}$$ where $E \coloneqq \Gamma - \Xi$. Since $\Gamma = \Phi \Gamma \Phi^{\prime} + I_{N}$ and $\Xi = \Phi \Phi^{\prime}\Xi + I_{N}$, then $E = \Phi E \Phi^{\prime} + \Pi$ where $\Pi \coloneqq \Phi \Xi \Phi^{\prime} - \Phi \Phi^{\prime}\Xi$. This is an example of a Lyapunov matrix equation and its formal solution $E = \sum_{k=0}^{\infty} \Phi^{k} \Pi (\Phi^{\prime})^{k}$ converges since $\lVert \Phi \rVert_{2} \leq 1$ \citepSM{smith1968matrix}. By the triangle inequality, $\lVert E \rVert_{2} \leq \lVert \Pi \rVert_{2} \sum_{k=0}^{\infty} \lVert \Phi \rVert_{2}^{2k} = \lVert \Pi \rVert_{2} / \{1 - \sigma_{1}^{2}(\Phi)\}$. 

From the proof of Lemma \ref{lemma:eigs}, $\lambda_{r}(P) \asymp \alpha_{N}N\lambda_{r}(B)$ for $r \in [K]$ and $\lambda_{r}(P)=0$ for $r = K+1,\dots,N$. By \citetSM[][Corollary 7.3.5]{horn2012matrix}, $|\sigma_{r}(A) - \sigma_{r}(P)| \leq \lVert A - P\rVert_{2}$. By \citetSM[][Proposition 1]{gallagher2023spectral}, $\lVert A - P\rVert_{2} = O(N^{1/2}\log^{\alpha + 1/2}N)$ almost surely for some $\alpha >0$. Therefore, $\sigma_{r}(\Phi) \asymp \rho \lambda_{r}(B)/\lambda_{1}(B)$ for $r \in [K]$ and $\sigma_{r}(\Phi) = o(1)$ for $r = K+1,\dots N$ almost surely. This implies $\lVert E \rVert_{2} \leq \lVert \Pi \rVert_{2} / \{1 - \rho^{2}\}$ almost surely. 

Next, we bound $\lVert \Pi \rVert_{2}$. Let $V_{\Phi,N} \coloneqq (V_{\Phi}|V_{\Phi,\perp}) \in \mathbb{R}^{N \times N}$ and $\Sigma_{\Phi,N} \coloneqq (\Sigma_{\Phi}|\Sigma_{\Phi,\perp})$ so that the full singular value decomposition of $\Phi$ is $\Phi = U_{N,\Phi} \Sigma_{N,\Phi}V_{N,\Phi}^{\prime}$. Define $R_{\Phi} \coloneqq V_{N,\Phi}^{\prime}U_{N,\Phi}$ so that 
\begin{align}
\Pi &= U_{N,\Phi} \Sigma_{N,\Phi}R_{\Phi} \Lambda_{\Xi}R_{\Phi}^{\prime}\Sigma_{N,\Phi} U_{N,\Phi}^{\prime} - U_{N,\Phi} \Sigma_{N,\Phi}^{2} \Lambda_{\Xi} U_{N,\Phi}^{\prime} \\
&=  U_{N,\Phi} \Sigma_{N,\Phi}(R_{\Phi} \Lambda_{\Xi}R_{\Phi}^{\prime} -   \Lambda_{\Xi})\Sigma_{N,\Phi} U_{N,\Phi}^{\prime} \\
&= U_{N,\Phi} \Sigma_{N,\Phi}(R_{\Phi} \Lambda_{\Xi} -   \Lambda_{\Xi}R_{\Phi})R_{\Phi}^{\prime}\Sigma_{N,\Phi} U_{N,\Phi}^{\prime} \\
&= U_{N,\Phi} \Sigma_{N,\Phi}\Lambda_{\Xi}\{(1-\Sigma_{N,\Phi}^{2})R_{\Phi}  -   R_{\Phi}(1-\Sigma_{N,\Phi}^{2})\}\Lambda_{\Xi}R_{\Phi}^{\prime}\Sigma_{N,\Phi} U_{N,\Phi}^{\prime} \\
&= U_{N,\Phi} \Sigma_{N,\Phi}\Lambda_{\Xi}(R_{\Phi}\Sigma_{N,\Phi}^{2}  -   \Sigma_{N,\Phi}^{2}R_{\Phi})\Lambda_{\Xi}R_{\Phi}^{\prime}\Sigma_{N,\Phi} U_{N,\Phi}^{\prime}
\end{align}
where we used the fact that $R_{\Phi}R_{\Phi}^{\prime} = I_{N}$ and $\Lambda_{\Xi} = (1-\Sigma_{N,\Phi}^{2})^{-1}$. Since $\Phi\Phi^{\prime} = U_{N,\Phi} \Sigma_{N,\Phi}^{2} U_{N,\Phi}^{\prime}$ and $\Phi^{\prime}\Phi = V_{N,\Phi} \Sigma_{N,\Phi}^{2} V_{N,\Phi}^{\prime}$, then 
\begin{align}
    V_{N,\Phi}^{\prime}(\Phi\Phi^{\prime} - \Phi^{\prime}\Phi)V_{N,\Phi} = R_{\phi}\Sigma_{N,\Phi}^{2}R_{\phi}^{\prime} - \Sigma_{N,\Phi}^{2} = (R_{\phi}\Sigma_{N,\Phi}^{2} - \Sigma_{N,\Phi}^{2}R_{\phi})R_{\Phi}^{\prime}.
\end{align}
Since $\lVert \Lambda_{\Xi}\rVert_{2} \asymp 1/(1-\rho^{2})$ almost surely and $\lVert U_{N,\Phi} \rVert_{2} = \lVert V_{N,\Phi} \rVert_{2} = \lVert R_{\Phi} \rVert_{2} = 1$, then 
\begin{align}
    \lVert \Pi \rVert \leq \frac{\rho^{2}}{(1-\rho^{2})^{2}} \lVert \Phi\Phi^{\prime} - \Phi^{\prime}\Phi\rVert_{2}
\end{align}
almost surely. The result then follows by noting that the eigengap satisfies $1/\{\lambda_{K}(\Xi) - \lambda_{K+1}(\Xi)\} \asymp 1/\rho^{2}\lambda_{K}^{2}(B)$. \qed

\subsection{Proof of Theorem~\ref{thm:hetero-strong-cov}}
The proof is identical to that of Theorem \ref{thm:strong-cov-estimation} but with different bounds on $v^{2}$ and $v_{\perp}^{2}$. We therefore use the same definitions as in the proof of Theorem \ref{thm:strong-cov-estimation}.  Starting with $v^{2}$, by Isserlis' theorem, 
\begin{align}
    \mathbb{E}(Z_{t}Z_{t+h}) &= \mathbb{E}(X_{i,t}^{U}Y_{j,t}X_{i,t+h}^{U}Y_{j,t+h}) - (U\Lambda_{\Gamma})_{ij}^{2} \\
    &= \mathbb{E}(X_{i,t}^{U}Y_{j,t})\mathbb{E}(X_{i,t+h}^{U}Y_{j,t+h}) + \mathbb{E}(X_{i,t}^{U}X_{i,t+h}^{U})\mathbb{E}(Y_{j,t}Y_{j,t+h}) \\ &\qquad+ \mathbb{E}(X_{i,t}^{U}Y_{j,t+h})\mathbb{E}(Y_{j,t}X_{i,t+h}^{U}) - (U\Lambda_{\Gamma})_{ij}^{2}. 
\end{align}
Let $\bm{e}_{i}$ denote the standard basis vector with a $1$ in the $i$th position and zeros in all other entries. Now $\mathbb{E}(X_{i,t}^{U}Y_{j,t+h}) = \bm{e}_{i}^{\prime}UU^{\prime}\Gamma(h)U\bm{e}_{j}$ which implies $|\mathbb{E}(X_{i,t}^{U}Y_{j,t+h})| \leq \lVert \Gamma(h) \rVert_{2} \lVert U \rVert_{2,\infty}$ where we used the fact that $\lVert \bm{e}_{i}^{\prime}U \rVert_{2} \leq \lVert U \rVert_{2,\infty}$. Similarly, $|\mathbb{E}(X_{i,t}^{U}X_{i,t+h}^{U})| \leq \lVert \Gamma(h) \rVert_{2} \lVert U \rVert_{2,\infty}^{2}$, and $|\mathbb{E}(Y_{j,t}Y_{j,t+h})| \leq \lVert \Gamma(h) \rVert_{2}$. By the assumption of balanced communities, $ \lVert U \rVert_{2,\infty} \lesssim \sqrt{K/N}$. By Lemma \ref{lemma:VAR-autocovariance-bound}, $\lVert \Gamma(h) \rVert_{2} \lesssim \lVert \Sigma \rVert_{2}\rho^{h}/(1-\rho^{2})$. Therefore, $|\mathbb{E}(Z_{t}Z_{t+h})| \lesssim K\lVert \Sigma \rVert_{2}^{2}\rho^{2h}/\{(1-\rho^{2})^{2}N\}$ and thus 
\begin{align}
    v^{2} \lesssim \frac{K\lVert \Sigma \rVert_{2}^{2}}{(1-\rho^{2})^{3}N}.
\end{align}
Next, we bound $v_{\perp}^{2}$. Let $s \equiv t+h$ for $h \in \{0,\dots,T-1\}$. By Isserlis' theorem, 
\begin{multline}
    \mathbb{E}(X_{i,t}^{U_{\perp}} Y_{j,t} X_{i,s}^{U_{\perp}} Y_{j,s}) = \\ \mathbb{E}(X_{i,t}^{U_{\perp}} Y_{j,t})\mathbb{E}( X_{i,s}^{U_{\perp}} Y_{j,s}) + \mathbb{E}(X_{i,t}^{U_{\perp}} X_{i,s}^{U_{\perp}})\mathbb{E}( Y_{j,t} Y_{j,s}) + \mathbb{E}(X_{i,t}^{U_{\perp}} Y_{j,s})\mathbb{E}( X_{i,s}^{U_{\perp}} Y_{j,t}).
\end{multline}
We note that $\mathbb{E}(X_{i,t}^{U_{\perp}} Y_{j,t+h}) = \bm{e}_{i}^{\prime}U_{\perp}U_{\perp}^{\prime}\Gamma(h)U_{\perp}\bm{e}_{j}$ which implies $|\mathbb{E}(X_{i,t}^{U_{\perp}}Y_{j,t+h})| \leq \lVert \Gamma(h) \rVert_{2}$. Similarly, $|\mathbb{E}(X_{i,t}^{U_{\perp}}X_{i,t+h}^{U_{\perp}})|~\lesssim~\lVert \Gamma(h) \rVert_{2} $, and therefore $|\mathbb{E}(X_{i,t}^{U_{\perp}} Y_{j,t} X_{i,s}^{U_{\perp}} Y_{j,s})| \lesssim \lVert \Gamma(h) \rVert_{2}^{2} \lesssim \lVert \Sigma \rVert_{2}^{2}\rho^{2h}/(1-\rho^{2})^{2}$. This implies 
\begin{align}
    v_{\perp}^{2} \lesssim \frac{\lVert \Sigma \rVert_{2}^{2}}{(1-\rho^{2})^{3}}.
\end{align}
By Lemma \ref{lemma:S_T-Gamma-general-errors} and the Davis-Kahan sin$\Theta$ theorem, for any $\nu \geq 8$, with probability exceeding $1 - (T+N)^{-\nu}$, 
    \begin{align} 
        \lVert \sin \Theta (U_{S},U) \rVert_{2} &\lesssim \frac{\lVert \Sigma \rVert_{2}}{\rho^{2}\lambda_{K}^{2}(B)(1 - \rho^{2})^{3/2}} \left( \sqrt{\frac{K}{T}}  +\sqrt{\frac{N}{T}} \right)  \log^{\frac{1}{2}}(T+N) .
    \end{align}
Repeating the proof of Theorem \ref{thm:strong-cov-estimation} with these bounds on $\lVert \sin \Theta (U_{S},U) \rVert_{2}$, $v^{2}$, and $v_{\perp}^{2}$, then yields the result. \qed

\bibliographystyleSM{rss}
\bibliographySM{references}

\end{document}